\documentclass[12pt]{amsart}
\usepackage{amsfonts,amssymb,amscd,amsmath}
\usepackage{latexsym}
\usepackage{mathrsfs}
\usepackage{tocvsec2}
\usepackage[%dvipdfm,
hyperindex=true,pagebackref=true,bookmarks=true,
colorlinks=true,linkcolor=blue,citecolor=red]{hyperref}
\def\~{{\rm --}}
\begin{document}

%MACOS FOR LECTURES ON DAHA
\renewcommand{\tilde}{\widetilde}
\renewcommand{\hat}{\widehat}

\newcommand{\BR}{{\mathbb R}}
\newcommand{\BQ}{{\mathbb Q}}
\newcommand{\BC}{{\mathbb C}}
\newcommand{\BP}{{\mathbb P}}
\newcommand{\BZ}{{\mathbb Z}}
\newcommand{\BN}{{\mathbb N}}
\newcommand{\BS}{{\mathbb S}}

\newcommand{\cH}{{\mathcal H}}
\newcommand{\cA}{{\mathcal A}}
\newcommand{\cB}{{\mathcal B}}
\newcommand{\ccF}{{\mathfrak F}}
\newcommand{\cD}{{\mathcal D}}
\newcommand{\cL}{{\mathcal L}}
\newcommand{\cF}{{\mathcal F}}
\newcommand{\cP}{{\mathcal P}}
\newcommand{\cX}{{\mathcal X}}
\newcommand{\cY}{{\mathcal Y}}
\newcommand{\cS}{{\mathcal S}}
\newcommand{\cSol}{\hbox{$\mathcal Sol$}}
\newcommand{\cT}{\hbox{$\mathcal T$}}

\newcommand{\Z}{{\mathbb Z}}
\newcommand{\Q}{{\mathbb Q}}
\newcommand{\N}{{\mathbb N}}
\newcommand{\C}{{\mathbb C}}
\newcommand{\R}{{\mathbb R}}
\newcommand{\X}{{\mathbb X}}
\newcommand{\Y}{{\mathbb Y}}

\newcommand{\CH}{{\mathcal H}}
\newcommand{\CA}{{\mathcal A}}

\def\HH{\mbox{${\mathcal H}$\kern-5.2pt${\mathcal H}$}}

\newcommand{\binomial}[2]{\genfrac{(}{)}{0pt}{}{ #1 }{ #2 }}
\newcommand{\qbinomial}[2]{\genfrac{[}{]}{0pt}{}{ #1 }{ #2 }_q }
\newcommand{\qbinom}[3]{\genfrac{[}{]}{0pt}{}{ #1 }{ #2 }_{ #3 } }

%%SPECIAL SEC 1.0

\def\der{\partial}
\def\tensor{\otimes}
\def\gam{\gamma} \def\Gam{\Gamma}
\def\del{\delta} \def\Del{\Delta}
\def\kap{\kappa}
\def\lam{\lambda} \def\Lam{\Lambda}
\def\Comp{{\mathbb C}}
\def\sM{{\mathcal M}}

\newtheorem{theorem}{Theorem}[section]
\newtheorem{maintheorem}[theorem]{Main Theorem}
\newtheorem{proposition}[theorem]{Proposition}
\newtheorem{definition}[theorem]{Definition}
\newtheorem{lemma}[theorem]{Lemma}
\newtheorem{corollary}[theorem]{Corollary}
\newtheorem{notation}[theorem]{Notation}
\newtheorem{remark}[theorem]{Remark}
\newtheorem{example}[theorem]{Example}

\newtheorem{theorem }{Theorem}[section]
\newtheorem{maintheorem }[theorem]{Main Theorem}
\newtheorem{proposition }[theorem]{Proposition}
\newtheorem{definition }[theorem]{Definition}
\newtheorem{lemma }[theorem]{Lemma}
\newtheorem{corollary }[theorem]{Corollary}
\newtheorem{notation }[theorem]{Notation}
\newtheorem{remark }[theorem]{Remark}
\newtheorem{example }[theorem]{Example}

\newtheorem{ maintheorem }[theorem]{Main Theorem}
\newtheorem{ theorem}{Theorem}[section]
\newtheorem{ proposition}[theorem]{Proposition}
\newtheorem{ definition}[theorem]{Definition}
\newtheorem{ lemma}[theorem]{Lemma}
\newtheorem{ corollary}[theorem]{Corollary}
\newtheorem{ notation}[theorem]{Notation}
\newtheorem{ remark}[theorem]{Remark}
\newtheorem{ example}[theorem]{Example}

\newtheorem{thm}{Theorem}[section]
\newtheorem{prop}[thm]{Proposition}
\newtheorem{lem}[thm]{Lemma}
\newtheorem{cor}[thm]{Corollary}
\newtheorem{conj}[thm]{Conjecture}
\newtheorem{con}[thm]{Conjecture}
\newtheorem{dfn}[thm]{Definition}
\newtheorem{df}[thm]{Definition}
 \newcommand{\rem}{{\bf Comment.\ }}
 \newcommand{\rmk}{{\bf Comment.\ }}
 \newcommand{\exmp}{{\bf Example.\ }}
 \newcommand{\ex}{{\bf Example.\ }}
 \newcommand{\prob}{{\bf Problem.\ }}

\newtheorem{note}{Note} 
\renewcommand{\thenote}{}
\newtheorem*{acka}{Acknowledgments}
\newtheorem{ack}{Acknowledgments}
\renewcommand{\theack}{}
\renewcommand{\appendixname}{\bf Appendix}
\renewcommand{\proof}{{\em Proof.\ }}

\hyphenation{
ap-pen-dix as-ymp-tot-ic at-trib-uted at-trib-ut-able
Bry-li-n-sky com-mu-ta-tion de-ge-ne-rate
de-riv-a-tive dis-trib-ute equi-vari-ant ex-tra-or-di-nary  
geo-met-ric griev-ance griev-ous grad-ed ho-lo-no-my ho-mo-thetic
in-fin-ite-ly in-fin-i-tes-i-mal Ha-rish Cha-n-dra mul-ti-plic-able 
non-euclid-ean non-iso-mor-phic non-smooth par-a-digm 
par-a-bol-ic pa-rab-o-loid pa-ram-e-trize phe-nom-e-non 
post-script pseu-do-dif-fer-en-tial pseu-do-fi-nite 
qua-drat-ics quad-ra-ture Han-kel rec-tan-gle semi-def-i-nite 
set-up wide-spread Euler-ian Feb-ru-ary Gauss-ian Grothen-dieck 
Hamil-ton-ian Her-mi-t-ian her-mi-t-ian Jan-u-ary 
Japan-ese Ka-shi-wa-ra Kor-te-weg Le-gendre No-vem-ber Rie-mann-ian 
Sep-tem-ber Za-mo-lo-d-chi-kov Kni-zh-nik quan-tum Op-dam
Mac-do-nald Ca-lo-ge-ro Su-ther-land Mo-ser 
Ol-sha-net-sky  Pe-re-lo-mov in-de-pen-dent ope-ra-tors 
cy-clo-to-mic ra-tio-nal de-gen-er-a-tion 
in-ter-est-ing de-for-ma-tions de-for-ma-tion pro-ce-dure 
fol-lows ope-ra-tors  pre-serve suf-fices ap-proach 
for-mu-las con-sider its com-ple-tion cor-re-spond-ing 
au-to-mor-phism be-cause pro-por-tional fi-nal-ly let-ting 
equi-v-a-lence ge-n-er-al-ized Mac-do-nald iden-ti-ties 
cor-re-s-pond sub-dia-grams par-ti-tion na-t-u-ral-ly 
or-dered stan-dard de-for-ma-tion ar-gu-ment com-bined 
sphe-r-i-cal rep-re-sen-ta-tions tri-go-no-me-t-ric
ge-n-er-al-ly speak-ing pri-m-it-ive ir-re-du-cible 
sum-ma-tion  rep-re-sen-ta-tives pro-por-ti-o-na-li-ty
ultra-sphe-ri-cal Ro-gers}

\def\ffor{\quad\hbox{ for }\quad}
\def\wwhen{\quad\hbox{ when }\quad}
\def\wwhere{\quad\hbox{ where }\quad}
\def\aand{\quad\hbox{ and }\quad}
\def\for{\  \hbox{ for } \ }
\def\iif{ \ \hbox{ if } \ }
\def\when{ \ \hbox{ when } \ }
\def\where{\  \hbox{ where } \ }
\def\and{\  \hbox{ and } \ }
\def\and{\  \hbox{ and } \ }
\def\oor{\  \hbox{ or } \ }
\def\proof{{\em Proof. \  }}

\def\equal{\stackrel{\,\mathbf{def}}{= \kern-3pt =}}

\def\la{\lambda}
\def\La{\Lambda}
\def\om{\omega}
\def\Om{\Omega}
\def\Th{\Theta}
\def\th{\theta}
\def\al{\alpha}
\def\be{\beta}
\def\ga{\gamma}
\def\ep{\epsilon}
\def\up{\upsilon}
\def\Up{\Upsilon}
\def\de{\delta}
\def\De{\Delta}
\def\ka{\kappa}
\def\kapp{\hbox{\bf \ae}}
\def\si{\sigma}
\def\Si{\Sigma}
\def\Ga{\Gamma}
\def\ze{\zeta}
\def\io{\iota}
\def\bio{b^\iota}
\def\aio{a^\iota}
\def\twio{\tilde{w}^\iota}
\def\hwio{\hat{w}^\iota}
\def\gio{\g^\iota}
\def\Bio{B^\iota}

\def\del{\delta}
\def\pa{\partial}
\def\vp{\varphi}
\def\ve{\varepsilon}
\def\inf{\infty}

\def\vph{\varphi}
\def\vps{\varpsi}
\def\vPh{\varPhi}
\def\vep{\varepsilon}
\def\vpi{{\varpi}}
\def\vth{{\vartheta}}
\def\vsi{{\varsigma}}
\def\vrh{{\varrho}}

\def\bph{\bar{\phi}}
\def\bsi{\bar{\si}}
\def\bvp{\bar{\varphi}}

\newcommand{\bS}{{\mathbf S}}
\newcommand{\bH}{{\mathbf H}}
\newcommand{\bF}{{\mathbf F}}
\newcommand{\bE}{{\mathbf E}}

\def\tal{\tilde{\alpha}}
\def\tbe{\tilde{\beta}}
\def\tde{\tilde{\delta}}
\def\tpi{\tilde{\pi}}
\def\txi{\tilde{\xi}}
\def\tPi{\tilde{\Pi}}
\def\tPhi{\tilde{\Phi}}
\def\tV{\tilde{V}}
\def\tJ{\tilde{J}}
\def\tla{\tilde{\lambda}}
\def\tga{\tilde{\gamma}}
\def\tGa{\tilde{\Gamma}}
\def\tvs{\tilde{{\varsigma}}}
\def\tu{\tilde{u}}
\def\tU{\tilde{U}}
\def\tw{\widetilde w}
\def\tW{\widetilde W}
\def\tB{\tilde B}
\def\tv{\tilde v}
\def\tV{\tilde V}
\def\tz{\tilde z}
\def\tb{\tilde b}
\def\ta{\tilde a}
\def\tih{\tilde h}
\def\trh{\tilde {\rho}}
\def\tx{\tilde x}
\def\tf{\tilde f}
\def\tg{\tilde g}
\def\tG{\tilde G}
\def\tk{\tilde k}
\def\tl{\tilde l}
\def\tL{\tilde L}
\def\tD{\tilde D}
\def\tR{\tilde R}
\def\tP{\tilde P}
\def\tH{\tilde H}
\def\tp{\tilde p}

\def\hH{\hat{H}}
\def\hh{\hat{h}}
\def\hR{\hat{R}}
\def\hY{\hat{Y}}
\def\hX{\hat{X}}
\def\hP{\hat{P}}
\def\hT{\hat{T}}
\def\hV{\hat{V}}
\def\hG{\hat{G}}
\def\hF{\hat{F}}
\def\hw{\widehat{w}}
\def\hW{\widehat{W}}
\def\hu{\hat{u}}
\def\hs{\hat{s}}
\def\hv{\hat{v}}
\def\hb{\hat{b}}
\def\hB{\widehat{B}}
\def\hze{\hat{\zeta}}
\def\hsi{\hat{\sigma}}
\def\hrh{\hat{\rho}}
\def\hth{\hat{\theta}}
\def\hy{\hat{y}}
\def\hx{\hat{x}}
\def\hz{\hat{z}}
\def\hg{\hat{g}}
\def\he{\hat{e}}
\def\hE{\widehat{E}}

\def\B{\mathbf{B}}
\def\I{\mathbf{I}}
\def\P{\mathbf{P}}
\def\G{\mathbf{G}}
\def\S{\mathbf{S}}
\def\F{\mathbf{F}}
\def\one{\mathbf{1}}
\def\Sn{\mathbf{S}_n}
\def\0{\mathbf{0}}
\def\H{\mathbf{H}}
\def\V{\mathbf{V}}

\def\f{\mathcal{F}}
\def\çF{\mathcal{F}}
\def\o{\mathcal{O}}
\def\t{\mathcal{T}}
\def\r{\mathcal{R}}
\def\l{\mathcal{L}}
\def\m{\mathcal{M}}
\def\k{\mathcal{K}}
\def\n{\mathcal{N}}
\def\d{\mathcal{D}}
\def\p{\mathcal{P}}
\def\cP{\mathcal{P}}
\def\a{\mathcal{A}}
\def\h{\mathcal{H}}
\def\c{\mathcal{C}}
\def\y{\mathcal{Y}}
\def\e{\mathcal{E}}
\def\v{\mathcal{V}}
\def\z{\mathcal{Z}}
\def\x{\mathcal{X}}
\def\s{\mathcal{S}}
\def\g{\mathcal{G}}
\def\u{\mathcal{U}}
\def\w{\mathcal{W}}
\def\i{\mathcal{I}}
\def\j{\mathcal{J}}
\def\b{\mathcal{B}}

\def\lan{\langle}
\def\llb{(\!(}
\def\ran{\rangle}
\def\rrb{)\!)}
 \def\dim{{\hbox{\rm dim}}_{\mathbb C}\,}
\def\lng{\hbox{\rm{\tiny lng}}}
\def\sht{\hbox{\rm{\tiny sht}}}
\def\sph{\hbox{\rm{\tiny sph}}}
\def\inv{\hbox{\rm{\tiny inv}}}

\def\br#1{\langle #1 \rangle}

\def\rank{\hbox{rank}}
\def\gl{\mathfrak{gl}_N}
%\def\sgn{\hbox{sgn}}
%\font\germ=eufb10 %at 12pt 
%\def\mathfrak#1{\hbox{\germ #1}}

\newcommand{\Aut}{\operatorname{Aut}}
\newcommand{\Hom}{\operatorname{Hom}}
\newcommand{\End}{\operatorname{End}}
\newcommand{\Ind}{\operatorname{Ind}}
\newcommand{\ad}{\operatorname{ad}}
\newcommand{\pr}{\operatorname{pr}}
\newcommand{\aweyl}{\tilde{\mathbb S}_n}
\newcommand{\hec}{{\mathcal H}^t_n}
\newcommand{\Func}{{\mathcal F}({\mathbb C}^n,{\mathcal H}^t_n)}
\newcommand{\tr}{\operatorname{tr}}
\newcommand{\Out}{\operatorname{Out}}
\newcommand{\Rad}{\operatorname{Rad}}
\newcommand{\Spec}{\operatorname{Spec}}
\newcommand{\id}{\operatorname{id}}
\newcommand{\Int}{\operatorname{Int}}
\newcommand{\ct} {\operatorname{ct}}

\newcommand{\rat}{{\mathbb Q}}
\newcommand{\real}{{\mathbb R}}
\newcommand{\cplx}{{\mathbb C}}
\newcommand{\zint}{{\mathbb Z}}

\newcommand{\sq}{\phantom{1}\hfill$\qed$}
\newcommand{\Rea}{\Re}
\newcommand{\Ima}{\Im}

\newcommand{\st}{\bowtie}
\newcommand{\modd}{\mbox{\,mod\,}}
\newcommand{\lr}{\langle}
\newcommand{\rr}{\rangle}
\newcommand{\eps}{\varepsilon}
\newcommand{\phk}{\phi^{(k)}}
\newcommand{\psk}{\psi^{(k)}}
\newcommand{\Res}{\mbox{Res}\;}
\newcommand{\sgn}{\mbox{sgn}}
\newcommand{\mn} {\left\{ \begin{array}{c}m\\
n\end{array}\right\}}

\def\sX{\mathscr{X}}
\def\sH{\mathscr{H}}
\def\sY{\mathscr{Y}}
\def\TT{\mathfrak{T}}
\def\JJ{\mathfrak{J}}
\def\HH{\mathfrak{H}}
\def\FF{\mathfrak{F}}
\def\GG{\mathfrak{G}}
\def\CC{\mathfrak{C}}
\def\LL{\mathfrak{L}}

\def\BB{\mathfrak{B}}
\def\AA{\mathfrak{A}}
\def\ZZ{\mathfrak{Z}}
\def\HH{\hbox{${\mathcal H}$\kern-5.2pt${\mathcal H}$}}
\def\HHH{\hbox{${\mathbb H}$\kern-4.2pt${\mathbb H}$}}
\def\tHH{\widetilde{\HH\ }}

\font\smm=msbm10 at 12pt 
\def\symbol#1{\hbox{\smm #1}}
\def\lsmash{{\symbol n}}
\def\rsmash{{\symbol o}}
\def\#{\sharp}

\font\tenbf=cmbx10
\font\tenrm=cmr10
\font\tenit=cmti10
\font\ninebf=cmbx9
\font\ninerm=cmr9
\font\nineit=cmti9
\font\eightbf=cmbx8
\font\eightrm=cmr8
\font\eightit=cmti8
\font\sevenrm=cmr7
\font\sevenbf=cmbx7

%END MACROS

\title [DAHA-Jones polynomials of torus knots]
{DAHA-Jones polynomials of torus knots}
\author[Ivan Cherednik]{Ivan Cherednik $^\dag$}
%\date{February 2, 2014}

\begin{abstract}
DAHA-Jones polynomials of torus knots $T(r,s)$ are 
studied systematically for reduced root systems
and in the case of $C^\vee C_1$. We prove the 
polynomiality and evaluation conjectures from the 
author's previous paper on torus knots and
extend the theory by the color exchange and further 
symmetries. DAHA-Jones polynomials for $C^\vee C_1$ 
depend on $5$ parameters. Their surprising connection 
to the DAHA-superpolynomials (type $A$) for the 
knots $T(2p+1,2)$ is obtained, a remarkable combination 
of the color exchange conditions and the author's duality 
conjecture (justified by Gorsky and Negut). 
The DAHA-superpolynomials for symmetric and wedge powers
(and torus knots) conjecturally coincide with the 
Khovanov-Rozansky stable polynomials, those originated in 
the theory of BPS states and the superpolynomials defined via 
rational DAHA in connection with certain Hilbert schemes, 
though not much is known about such connections beyond the 
HOMFLYPT and Kauffman polynomials.
We also define certain arithmetic counterparts of DAHA-Jones 
polynomials for the absolute Galois group instead of torus 
knots in the case of $C^\vee C_1$.
\end{abstract}

\thanks{$^\dag$ \today.
\ \ \ Partially supported by NSF grant
DMS--1101535 and the Simons Foundation}

\address[I. Cherednik]{Department of Mathematics, UNC
Chapel Hill, North Carolina 27599, USA\\
chered@email.unc.edu}

 \def\sht{\raisebox{0.4ex}{\hbox{\rm{\tiny sht}}}}
 \def\bysame{{\bf --- }}
 \def\~{{\bf --}}
 \def\rr{{\mathsf r}}
 \def\ss{{\mathsf s}}
 \def\mm{{\mathsf m}}
 \def\pp{{\mathsf p}}
 \def\ll{{\mathsf l}}
 \def\NS{\hbox{\tiny\sf ns}}
\newcommand{\comment}[1]{}
\renewcommand{\tilde}{\widetilde}
\renewcommand{\hat}{\widehat}
\renewcommand{\V}{\mathbb{V}}
\renewcommand{\F}{\mathbb{F}}
\newcommand{\dagx}{\hbox{\tiny\mathversion{bold}$\dag$}}
\newcommand{\ddagx}{\hbox{\tiny\mathversion{bold}$\ddag$}}
\newtheorem{conjecture}[theorem]{Conjecture}
\newcommand*\toeq{
\raisebox{-0.15 em}{\,\ensuremath{
\xrightarrow{\raisebox{-0.3 em}{\ensuremath{\sim}}}}\,}
}

\vskip -0.0cm
%\par
%{\centering
%\medskip
%\par}
%\vskip -0.0cm
\maketitle
\smallskip

{\tiny
\centerline{{\bf MSC} (2010): 14G32, 17B22, 17B45, 20C08,
20F36, 33D52, 30F10, 57M25, 57R58,}
}
\smallskip

\vskip -0.0cm
\renewcommand{\baselinestretch}{1.2}
{\textmd
\tableofcontents
}
\renewcommand{\baselinestretch}{1.0}
%}
%$\mathfrak{a,b,c,d,e,f,g,h,i,j,k,l,m,n,o,p,q,r,s,t,e^u,e^v,e^w}$
\vfill\eject

\renewcommand{\natural}{\wr}

\setcounter{section}{-1}
\setcounter{equation}{0}
\section{\sc Introduction}
This paper is a systematic exposition of the
new theory of {\em DAHA-Jones polynomials\,} of torus 
knots $T(\rr,\ss)$ for reduced root systems and 
in the case of $C^\vee C_1$. We prove the (key) 
Polynomiality Conjecture from \cite{CJ} and the 
Evaluation Conjecture there, also extending 
the theory by the important Color Exchange, 
the symmetry $(\rr,\ss)\!\mapsto\!(\ss,\rr)$ and
further properties of DAHA-Jones polynomials
and {\em superpolynomials\,}. We mention that the 
justification of the Stabilization Conjecture in type $A$ 
is provided in \cite{GoN}. It was announced as proven
in \cite{CJ} (with a link to \cite{SV}); it
justifies the existence of DAHA-superpolynomials.
\smallskip

{\sf Superpolynomials.}
The DAHA-superpolynomials of $T(\rr,\ss)$ conjecturally 
coincide with the Poincar\`e polynomials of stable
triply graded {\em Khovanov-Rozansky homology\,} of torus knots 
for $sl_{n+1}$ as $n\to \infty$ \cite{KhR1,KhR2,Ras}; see
also \cite{Web,Rou} and references there concerning the 
{\em categorification\,}.  Only {\em reduced\, }
Jones and other polynomials (with value $1$ at the unknot) 
are considered in this paper, as well as in \cite{CJ}. 
Such Connection Conjecture was checked for 
{\em HOMFLYPT and Kauffman polynomials\,} \cite{CJ,Ste}
and when the Khovanov-Rozansky polynomials are available; 
there is also a link to the {\em Heegard- Floer homology\,}. 
Our color exchange has something in common with 
the so-called colored differentials (see \cite{GGS}), 
but there are significant differences
and the existence of the latter is a conjecture.
\smallskip 

Relying on such a connection,
DAHA provide a self-consistent powerful
approach to torus knots (arbitrary root systems and weights). 
It is relatively simple. Only the positivity for the rectangles 
and the stabilization-duality in 
types $B\!-\!\!C\!\!-\!\!D$ remain open among the intrinsic properties 
of the DAHA-polynomials conjectured in \cite{CJ}. The $B\!-\!C\!-\!D$
systems are more natural to address within the $C^\vee C_n$\~theory, 
to be considered in author's further works. The system $C^\vee C_n$,
the most general classical one, can be managed following the
present paper. We note that counterparts of the properties 
of DAHA-superpolynomials are generally difficult challenges 
in the Khovanov-Rozansky theory.
\smallskip

Let us briefly mention other approaches to 
superpolynomials. The (physics) 
theory of the {\em BPS states\,} is used, which is not mathematically 
rigorous but fruitful; see e.g. \cite{DGR,AS,FGS,DMMSS}. We
mention that the duality for symmetric and wedge powers
was conjectured in \cite{GS}.
%by Gukov and Sto$\breve{\hbox{s}}$i\'{c}. 
The {\em rational DAHA\,} can be employed here (in connection 
with certain {\em Hilbert schemes\,}); see \cite{GORS,GoN} and 
references there. The superpolynomials obtained via rational 
DAHA perfectly match our DAHA-superpolynomials, but this approach is
restricted to symmetric or wedge powers so far. 
We do not discuss this direction and related important geometric 
developments; the link to our superpolynomials has now
no explanation (in spite of using DAHA in both theories).
\medskip

{\sf The $C^\vee C_1$\~case.} The
DAHA-Jones polynomials in this case depend on 
$5$ parameters and are related to those of 
type $A_1$ in nontrivial ways. Their surprising 
connection to the DAHA-superpolynomials (type $A$) 
for $T(2p+1,2)$ is discovered; its proof is a 
remarkable combination of the color exchange and 
the Duality Conjecture from \cite{CJ},
justified in \cite{GoN}. We will outline 
our approach to the proof of the latter
based on the new $q,t$\~level-rank 
duality, which deserves a separate paper.

The relation of $C^\vee C_1$ to superpolynomials requires 
special parameters. Geometric meaning of the whole set of
parameters is unclear at the moment. 
This restricts potential value of the general
$C^\vee C_n$\~theory, though there are solid reasons to 
follow this avenue.
\smallskip

We note that the $C^\vee C_1$\~theory somewhat deviates from 
that in \cite{CJ} (any reduced root systems). For instance, the 
symmetry $(\rr,\ss)\mapsto (\ss,\rr)$  holds in this case only 
if $\,\rr\ss\,$ is even. Also,  the DAHA-Jones polynomials 
for $(2\pp+1,1)$ become nontrivial for odd $\pp>1$, 
which are exactly the ones resulting in the $A$\~type
superpolynomials for $T(2\pp+1,2)$. For any $\ss\,$ and
odd $\rr\,$, the DAHA-Jones 
polynomials for $(\rr,\ss)$ in type $C^\vee C_1$  
coincide with those of type $A_1$ for $T(\rr,2\ss)$ under certain 
specialization of the parameters, which is, by the way, related to 
the connection between the continued fraction of $2\ss/\rr$ 
and  $\ss/\rr$ 
(Hurwitz).
\smallskip

Our formula for the superpolynomials 
of torus knots $T(2\pp+1,2)$ (conjectured in \cite{DMMSS,FGS})
is of independent interest due to the connection to 
(half-differentials of) {\em modular forms\,} via the 
{\em WRT and Kashaev invariants\,}; see \cite{Hi}.  
Hopefully it will clarify and generalize 
the modularity properties of the Jones polynomials.
Its proof is a convincing demonstration of the efficiency 
of the DAHA-based approach.

\smallskip

At the end of the paper, counterparts
of DAHA-Jones polynomials in type $C^\vee C_1$ are 
defined for the {\em absolute Galois group\,} instead of
torus knots, which is based on the rigidity of (``major")
{\em perfect DAHA-modules\,} at roots of unity
in the terminology of \cite{C101}.
We follow \cite{CG}, though do not 
have the classification of such perfect modules for $C^\vee C_1$. 
See \cite{ObS} concerning the rigidity for generic $\,q$
in the case of $C^\vee C_1$.
\medskip

We would like to mention paper \cite{BeS}, which approaches
the {\em skein modules\,} (for $sl_2$) via DAHA of type 
$C^\vee C_1$. 
The knots $T(\rr,2)$ for odd $\rr\,$ are among the main examples 
there, but not only torus knots are considered in this paper.
We see no connection with our work at 
the moment. Actually, such a connection is more likely at the 
level of {\em $A$\~polynomials\,} (see e.g., \cite{FGS}), which 
are beyond the present work. 
\medskip

\setcounter{section}{0}
\setcounter{equation}{0}
\section{\sc Reduced root systems}
\subsection{\bf Affine Weyl group}
Let $R=\{\al\}   \subset \R^n$ be a root system of type
$A,\!B,\!C,\!D,\!E,\!F,\!G$
with respect to a euclidean form $(z,z')$ on $\R^n
\ni z,z'$,
$W$ the Weyl group 
generated by the reflections $s_\al$,
$R_{+}$ the set of positive  roots
corresponding to fixed simple 
roots $\al_1,...,\al_n;$ $R_-=-R_+$. 
The form is normalized
by the condition  $(\al,\al)=2$ for 
{\em short\,} roots. 
The root lattice and the weight lattice are:
\begin{align}
& Q=\oplus^n_{i=1}\Z \al_i \subset P=\oplus^n_{i=1}\Z \om_i,
\notag
\end{align}
where $\{\om_i\}$ are fundamental weights:
$ (\om_i,\al_j^\vee)=\de_{ij}$ for the
coroots $\al^\vee=2\al/(\al,\al).$
Replacing $\Z$ by $\Z_{+}=\{m\in\Z, m\ge 0\}$, we obtain
$Q_+, P_+.$
Here and further see  e.g., \cite{Bo} or \cite{C101}. 

Setting 
$\nu_\al\equal (\al,\al)/2$,
the vectors $\ \tal=[\al,\nu_\al j] \in
\R^n\times \R \subset \R^{n+1}$
for $\al \in R, j \in \Z $ form the
{\em twisted affine root system\,}
$\tR \supset R$ ($z\in \R^n$ are identified with $ [z,0]$).
We add $\al_0 \equal [-\vth,1]$ to the simple
roots for the {\em maximal short root\,} $\vth\in R_+$.
The corresponding set
$\tR_+$ of positive roots is %%%% tR_+BOOK!!!
$R_+\cup \{[\al,\nu_\al j],\ \al\in R, \ j > 0\}$.

The set of the indices of the images of $\al_0$ by all
automorphisms of the affine Dynkin diagram will be denoted by 
$O$ ($O=\{0\} \for E_8,F_4,G_2$). Let $O'\equal\{r\in O, r\neq 0\}$.
The elements $\om_r$ for $r\in O'$ are  
{\em minuscule weights\,}, defined by the inequalities 
$(\om_r,\al^\vee)\le 1$ for all $\al \in R_+$. We set $\om_0=0$
for the sake of uniformity.

Given $\tal=[\al,\nu_\al j]\in \tR,  \ b \in P$, let
\begin{align}
&s_{\tal}(\tz)\ =\  \tz-(z,\al^\vee)\tal,\
\ b'(\tz)\ =\ [z,\ze-(z,b)]
\label{ondon}
\end{align}
for $\tz=[z,\ze] \in \R^{n+1}$.
The
{\em affine Weyl group\,} $\tW=\lan s_{\tal}, \tal\in \tR_+\ran)$ is
the semidirect product $W\lsmash Q$ of
its subgroups $W=$ $\lan s_\al,
\al \in R_+\ran$ and $Q$, where $\al$ is identified with
\begin{align*}
& s_{\al}s_{[\al,\,\nu_{\al}]}=\
s_{[-\al,\,\nu_\al]}s_{\al}\for
\al\in R.
\end{align*}

The {\em extended Weyl group\,} $ \hW$ is $W\lsmash P$, where
the corresponding action in $\R^{n+1}$ is 
\begin{align}
&(wb)([z,\ze])\ =\ [w(z),\ze-(z,b)] \for w\in W, b\in P.
\label{ondthr}
\end{align}
It is isomorphic to $\tW\lsmash \Pi$ for $\Pi\equal P/Q$. 
The latter group consists of $\pi_0=$id\, and the images $\pi_r$
of minuscule $\om_r$ in $P/Q$. 
%Note that 
%$\pi_r^{-1}$ is $\pi_{r^*}$, where $^*$ is the standard involution 
%(sometimes trivial) of the {\em nonaffine\,} Dynkin diagram.

The group $\Pi$
is naturally identified with the subgroup of $\hW$ of the
elements of the length zero; the {\em length\, } is defined as 
follows:
\begin{align*}
&l(\hw)=|\la(\hw)| \for \la(\hw)\equal\tR_+\cap \hw^{-1}(-\tR_+).
\end{align*}
One has $\om_r=\pi_r u_r$ for $r\in O'$, where $u_r$ is the 
element $u\in W$ of minimal length such that $u(\om_r)\in -P_+$. 
\smallskip

Setting $\hw = \pi_r\tw \in \hW$ for $\pi_r\in \Pi,\, \tw\in \tW,$
\,$l(\hw)$ coincides with the length of any reduced decomposition
of $\tw$ in terms of the simple reflections
$s_i, 0\le i\le n.$ We will also use the {\em partial lenghts\,}
$l_{\sht},l_{\lng}$, which count correspondingly short and long
$s_i$ in reduced decompositions.
\medskip

\subsection{\bf Definition of DAHA}
We follow \cite{CJ,C101}.
Let $\mm,$ be the least natural number
such that  $(P,P)=(1/\mm)\Z.$  Thus
$\mm=|\Pi|$ unless 
$\mm=2 \for D_{2k}$ and $\ \mm=1 \for B_{2k},C_{k}.$ 

The double affine Hecke algebra, {\em DAHA\,}, depends
on the parameters
$q, t_\nu\, (\nu\in \{\nu_\al\})\,$ and will be defined
over the ring
$Z_{q,t}\equal\Z[q^{\pm 1/\mm},t_\nu^{\pm 1/2}]$
formed by
polynomials in terms of $q^{\pm 1/\mm}$ and
$\{t_\nu^{1/2}\}.$ Note that the coefficients of the 
Macdonald polynomials will belong to 
$
\Q(q,t_\nu).
$

For $\tal=[\al,\nu_\al j] \in \tR,\ 0\le i\le n$, we set
\begin{align*}
&   t_{\tal} =t_{\al}=t_{\nu_\al}=q_\al^{k_\nu} ,\ \, 
q_{\tal}=q^{\nu_\al}, \ \, t_i = t_{\al_i},  q_i=q_{\al_i},
\end{align*}

Also, using here and below {\em\small sht,\ lng\,} instead 
of $\nu$, let
\begin{align*}
\rho_k\equal \frac{1}{2}\!\sum_{\al>0} k_\al \al=
k_{\sht}\rho_{\sht}\!+\!k_{\lng}\rho_{\lng},\ \,
\rho_\nu=\frac{1}{2}\!\sum_{\nu_\al=\nu} \al=
\!\!\sum_{\nu_i=\nu,i>0}  \om_i.
\end{align*}

For pairwise commutative $X_1,\ldots,X_n,$
\begin{align}
& X_{\tb}\ \equal\ \prod_{i=1}^nX_i^{l_i} q^{ j}
\iif \tb=[b,j],\ \hw(X_{\tb})\ =\ X_{\hw(\tb)},
\label{Xdex}\\
&\hbox{where\ } b=\sum_{i=1}^n l_i \om_i\in P,\ j \in
\frac{1}{ m}\Z,\ \hw\in \hW.
\notag \end{align}
For instance, $X_0\equal X_{\al_0}=qX_\vth^{-1}$.
\medskip

We use that $\pi_r^{-1}$ is $\pi_{\iota(i)}$, where
$\iota$ is the standard involution  of the nonaffine Dynkin diagram,
induced by $\al_i\mapsto -w_0(\al_i)$;
$w_0$ is the longest element in $W$. We set $m_{ij}=2,3,4,6$
when the number of links between $\al_i$ and $\al_j$ in the affine 
Dynkin diagram is $0,1,2,3$. Recall that 
$\om_r=\pi_r u_r$ for $r\in O'$ (see above).

\begin{definition}
The double affine Hecke algebra $\HH\ $
is generated over $\Z_{q,t}$ by
the elements $\{ T_i,\ 0\le i\le n\}$,
pairwise commutative $\{X_b, \ b\in P\}$ satisfying
(\ref{Xdex})
and the group $\Pi,$ where the following relations are imposed:

(o)\ \  $ (T_i-t_i^{1/2})(T_i+t_i^{-1/2})\ =\
0,\ 0\ \le\ i\ \le\ n$;

(i)\ \ \ $ T_iT_jT_i...\ =\ T_jT_iT_j...,\ m_{ij}$
factors on each side;

(ii)\ \   $ \pi_rT_i\pi_r^{-1}\ =\ T_j \iif
\pi_r(\al_i)=\al_j$;

(iii)\  $T_iX_b \ =\ X_b X_{\al_i}^{-1} T_i^{-1} \iif
(b,\al^\vee_i)=1,\
0 \le i\le  n$;

(iv)\ $T_iX_b\ =\ X_b T_i\ $ if $\ (b,\al^\vee_i)=0
\for 0 \le i\le  n$;

(v)\ \ $\pi_rX_b \pi_r^{-1}\ =\ X_{\pi_r(b)}\ =\
X_{ u^{-1}_r(b)}
 q^{(\om_{\iota(r)},b)},\  r\in O'$.
\label{double}
\end{definition}

Given $\tw \in \tW, r\in O,\ $ the product
\begin{align}
&T_{\pi_r\tw}\equal \pi_r T_{i_l}\cdots T_{i_1},\where
\tw=s_{i_l}\cdots s_{i_1} \for l=l(\tw),
\label{Twx}
\end{align}
does not depend on the choice of the reduced decomposition
Moreover,
\begin{align}
&T_{\hv}T_{\hw}\ =\ T_{\hv\hw}\  \hbox{ whenever\,}\
 l(\hv\hw)=l(\hv)+l(\hw) \for
\hv,\hw \in \hW. \label{TTx}
\end{align}
In particular, we arrive at the pairwise
commutative elements: 
\begin{align}
& Y_{b}\equal
\prod_{i=1}^nY_i^{l_i} \iif
b=\sum_{i=1}^n l_i\om_i\in P,\ 
Y_i\equal T_{\om_i},b\in P.
\label{Ybx}
\end{align}
\smallskip

\subsection{\bf The automorphisms}\label{sect:Aut}
We will begin with the anti-involution
\begin{align}
  & X_b^\star\ =\  X_b^{-1},\ \  Y_b^\star\ =\
 Y_b^{-1},\  \
 T_i^\star \ =\  T_i^{-1}, \notag\\
&t_\nu^{\upsilon} \mapsto t_\nu^{-\upsilon},\
 q^{\upsilon}\mapsto  q^{-\upsilon},\ 0\le i\le n,
\label{star}
\end{align}
where $\upsilon\in \Q$ ($\upsilon\in\frac{1}{2\mm}\Z$ will, 
actually, be sufficient). 
See \cite{C101},(3.2.18).
This is the group inversion for  
any products of the generators $X_b,Y_b,T_i,$, as
well as for $q^{1/\mm},t_\nu^{1/2}$, and 
therefore commutes with all automorphisms and anti-automorphisms
below.

The following maps can be (uniquely) extended to
an automorphism of $\HH\,$, where 
$q^{1/(2\mm)}$ must be added to $\Z_{q,t}$
(see \cite{C101}, (3.2.10)-(3.2.15)):
\begin{align}\label{tauplus}
& \tau_+:\  X_b \mapsto X_b, \ T_i\mapsto T_i\, (i>0),\
\ Y_r \mapsto X_rY_r q^{-\frac{(\om_r,\om_r)}{2}}\,,
\\
& \tau_+:\ T_0\mapsto  q^{-1}\,X_\vth T_0^{-1},\
\pi_r \mapsto q^{-\frac{(\om_r,\om_r)}{2}}X_r\pi_r\
(r\in O'),\notag\\
& \label{taumin}
\tau_-:\ Y_b \mapsto \,Y_b, \ T_i\mapsto T_i\, (i\ge 0),\
\ X_r \mapsto Y_r X_r q^\frac{(\om_r,\om_r)}{ 2},\\
&\tau_-(X_{\vth})= 
q T_0 X_\vth^{-1} T_{s_{\vth}}^{-1};\ \
\si\equal \tau_+\tau_-^{-1}\tau_+\, =\,
\tau_-^{-1}\tau_+\tau_-^{-1},\notag\\
&\si(X_b)=Y_b^{-1},\   \si(Y_b)=
T_{w_0}^{-1}X_{b^\iota}^{-1}T_{w_0},\ \si(T_i)=T_i (i>0).
\label{taux}
\end{align}
These automorphisms fix $\ t_\nu,\ q$
and their fractional powers, as well as the
following {\em anti-involution\,}:
\begin{align}
&\phi:\ 
X_b\mapsto Y_b^{-1},\, Y_b\mapsto X_b^{-1},\
T_i\mapsto T_i\ (1\le i\le n).\label{starphi}
%\\
%&\phi(\tau_+)\equal \phi\circ \tau_+\circ \phi\ =\ \tau_-\,,\ 
%\ \, \phi(\tau_-)\ =\ \tau_+.\notag 
\end{align}

Extending the standard involution $\iota$ of the
nonaffine Dynkin diagram used above, let 
$\iota(b)=b^\iota\equal -w_0(b)$ for $b\in P$ and
the longest element $w_0\in W$.
The relations 
\begin{align}\label{iotaXY}
\iota(X_b)=X_{\iota(b)},\ 
\iota(Y_b)=Y_{\iota(b)},\ T_i^\iota=T_{\iota(i)}
\end{align}
can be uniquely extended the automorphism
$H\mapsto \iota(H)\in \HH$. Proposition 3.2.2 from
\cite{C101} states that
\begin{align}\label{sisquare}
\si^2(H)\ =\ T_{w_0}^{-1}\,H^\iota\,T_{w_0} \hbox{\, for\, } 
H\in \HH.
\end{align}

We will also need the involutions 
$\vep\equal\vph\cdot \star=\star\cdot\vph$ (here
$\,\cdot\,$ is used for the composition) and $\eta 
\equal\vep\si=\si^{-1}\vep$:
\begin{align}
\vep:\ &X_b \mapsto Y_b,\ \  Y_b \mapsto X_b,\
\ T_i \mapsto T_i^{-1}\, (1\le i\le n),
\label{vep}\\
\eta:\ &T_i\mapsto T_i^{-1},\ X_b\mapsto X_b^{-1},\
\pi_r\mapsto \pi_r \,(0\le i\le n),
\label{etatxpi}
\end{align}
where $b\in P,\ r\in O'$.

Both ``conjugate" $t,q$; namely, 
$t^{\upsilon}_\nu \mapsto t_\nu^{-\upsilon}, 
q^{\upsilon}\mapsto  q^{-\upsilon}$.
The involution $\eta\,$ extends the
{\em Kazhdan\~Lusztig involution\,} in the affine Hecke theory;
see \cite{C101}, (3.2.19-22). Note that 
$$
\vep\tau_\pm\vep=\tau_\mp=
\vph\tau_\pm\vph=\si\tau_\pm^{-1}\si^{-1},\ 
\eta\tau_{\pm}\eta=\tau_{\pm}^{-1},\ \vph\si\vph=\si^{-1}=
\eta\si\eta.
$$
Also, $\vph \vep=\vph\vep$ and $\vph\eta\vph=\eta\si^{-2}=
\si^2\eta$.
Let us list the matrices corresponding to the automorphisms and
anti-automorphisms above upon the natural projection 
onto $GL_2(\Z)$, corresponding to  $t^\upsilon_\nu=1=q^\upsilon$ for 
any rational $\upsilon$ ($\upsilon\in\frac{1}{2\mm}\Z$ is sufficient). 
The matrix {\tiny 
$\begin{pmatrix} \al & \be \\ \ga & \de\\ \end{pmatrix}$}
will represent the map $X_b\mapsto X_b^\al Y_b^\ga,
Y_b\mapsto X_b^\be Y_b^\de$ for $b\in P$. One has:
\smallskip

\ \ \ \ \ $\tau_+\mapsto$ 
{\tiny 
$\begin{pmatrix}1 & 1 \\0 & 1 \\ \end{pmatrix}$},\ 
$\tau_-\mapsto$ 
{\tiny 
$\begin{pmatrix}1 & 0 \\1 & 1 \\ \end{pmatrix}$},\
$\si\mapsto$ 
{\tiny 
$\begin{pmatrix}0 & 1 \\-1 & 0 \\ \end{pmatrix}$},\
$\si^2\mapsto$ 
{\tiny 
$\begin{pmatrix}-1 & 0 \\0 & -1 \\ \end{pmatrix}$},\

\ \ \ \ \ \ $\vep\mapsto$ 
{\tiny 
$\begin{pmatrix}0 & 1 \\1 & 0 \\ \end{pmatrix}$},\
$\vph\mapsto$ 
{\tiny 
$\begin{pmatrix}0 & -1 \\-1 & 0 \\ \end{pmatrix}$},\
$\eta\mapsto$ 
{\tiny 
$\begin{pmatrix}-1 & 0 \\0 & 1 \\ \end{pmatrix}$},\
$\phi\si^2\mapsto$ 
{\tiny 
$\begin{pmatrix}0 & 1 \\1 & 0 \\ \end{pmatrix}$}.
\smallskip

We define the {\em projective\,} $GL_{\,2}^\wedge(\Z)\,$ as
the group generated by $\tau_{\pm}, \eta$ subject to the
relations $\tau_+\tau_-^{-1}\tau_+\, =\,
\tau_-^{-1}\tau_+\tau_-^{-1},\ \eta^2=1$ and 
$\eta\tau_{\pm}\eta=\tau_{\pm}^{-1}.$ The span 
of $\tau_{\pm}$ is the projective $PSL_{\,2}^\wedge(\Z)\,$
(due to Steinberg), which is isomorphic to the braid 
group $B_3$.

\medskip

\subsection{\bf Macdonald polynomials}
Following \cite{C101}, 
we use the PBW theorem to express any $H\in \HH$ in the form 
\,$\sum_{a,w,b} c_{a,w,b}\, X_a T_{w} Y_b$\, for $w\in W$,
$a,b\in P$ (this presentation is unique). Then we substitute:
\begin{align}\label{evfunct}
\{\,\}_{ev}:\ X_a \ \mapsto\  q^{-(\rho_k,a)},\ 
Y_b \ \mapsto\  q^{(\rho_k,b)},\ 
T_i \ \mapsto\  t_i^{1/2}. 
\end{align}
By construction, the resulting functional 
$\HH\ni H\mapsto \{H\}_{ev}$ acts via the projection 
$H\mapsto H(1)$ of $\HH\,$
onto the {\em polynomial representation \,}$\v$, which is
the $\HH$\~module induced from the one-dimensional
character $T_i(1)=t_i^{-1/2}=Y_i(1)$ for $1\le i\le n$ and
$T_0(1)=t_0^{-1/2}$; recall that $t_0=t_{\sht}$.
Explicitly,
$\{\,H\,\}_{ev}=H(1)(q^{-\rho_k})$; see \cite{C101,CJ}.

In detail, the polynomial representation $\v$
is isomorphic to $\Z_{q,t}[X_b]$
as a vector space and the action of $T_i(0\le i\le n)$ there 
is given by 
the {\em Demazure-Lusztig operators\,}:
\begin{align}
&T_i\  = \  t_i^{1/2} s_i\ +\
(t_i^{1/2}-t_i^{-1/2})(X_{\al_i}-1)^{-1}(s_i-1),
\ 0\le i\le n.
\label{Demazx}
\end{align}
The elements $X_b$ become the multiplication operators 
and  $\pi_r (r\in O')$ act via the general formula
$\hw(X_b)=X_{\hw(b)}$ for $\hw\in \hW$. Note that $\tau_-$ 
and $\eta$ naturally act in the polynomial representation. 
For the latter,
\begin{align}\label{eta-conj}
\eta(f)\!=\!f^\star, \hbox{\, where\, } X_b^\star\!=\!X_{-b}, 
(q^\upsilon)^\star\!=\!q^{-\upsilon}, (t^v)^\star\!=\!t^{-v} 
\hbox{\, for\, } \upsilon\in \Q.
\end{align}

%\smallskip

One has the following relations:
\begin{align}\label{evsym}
&\{\,\vph(H)\,\}_{ev}=\{\,H\,\}_{ev},\ \{\,\iota(H)\,\}_{ev}=
\{\,H\,\}_{ev},\\ 
&\{\,\eta(H)\,\}_{ev}=\{\,H\,\}_{ev}^\star,\ 
\{\,\si^2(H)\,\}_{ev}=\{\,H\,\}_{ev}.\notag 
\end{align}

The Macdonald polynomials $P_b(X)$ (due to Kadell
for the classical root systems) are uniquely defined
as follows. For $b\in P_+$,
\begin{align*}
&P_b\! -\!\!\!\!\sum_{b'\in W(b)}\!\!\! X_{b'} 
\in\, \oplus_{b_+\neq c_+\!\in b+Q_+}\,\Q_{q,t}' X_c, 
\hbox{ and\, }
CT\bigl(P_b X_{c}\,\mu(X;q,t)\bigr) = 0,
\\
&c_+\in W(c)\cap P_+,\ \, \mu(X;q,t)\!\equal\!\!\prod_{\al \in R_+}
\prod_{j=0}^\infty \frac{(1\!-\!X_\al q_\al^{j})
(1\!-\!X_\al^{-1}q_\al^{j+1})
}{
(1\!-\!X_\al t_\al q_\al^{j})
(1\!-\!X_\al^{-1}t_\al^{}q_\al^{j+1})}\,.
\end{align*}
Here $CT$ is the constant term; 
$\mu$ is considered
a Laurent series of $X_b$ with 
the coefficients expanded in terms of
positive powers of $q$. The coefficients of
$P_b$ belong to the $\Q(q,t_\nu)$ (see also below).
One has:
\begin{align}\label{macdsym}
&P_b(X^{-1})\,=\,P_{b^\iota}(X)\,=\,
P_b^\star(X),\ \, P_{b}(q^{-\rho_k})=
P_{b}(q^{\rho_k})\\
\label{macdeval}
=&(P_{b}(q^{-\rho_k}))^\star=
q^{-(\rho_k,b)}
\prod_{\al>0}\prod_{j=0}^{(\al^{\!\vee},b)-1}
\Bigl(
\frac{
1- q_\al^{j}t_\al X_\al(q^{\rho_k})
 }{
1- q_\al^{j}X_\al(q^{\rho_k})
}
\Bigr).
\end{align}

DAHA provides an important alternative (operator) 
approach to $P$\~polynomials; namely, they satisfy
the defining relations
\begin{align}\label{macdopers}
L_f(P_b)=f(q^{\rho_k+b})P_b,\  L_f\equal f(X_a\mapsto Y_a)
\end{align}
for any symmetric ($W$\~invariant) polynomial 
$f\in \C[X_a,a\in P]^W$.

These polynomials are $t$-symmetrizations of the
nonsymmetric Macdonald polynomials $E_b\in \v$, which
are defined for  $b\in P_+$ from the
relations $\,Y_a(E_b)=q^{-(a,b+w_b(\rho_k))}E_b\,$ for all $a$.
Here $w_b$ is the element of {\em maximal\,} length in the 
centralizer of $b\,$ in $W$. The normalization here is by 
the condition that the coefficient of $X_b$ in $E_b$ is $1$. 
More exactly, 
$E_b-X_b\in \oplus_{b_+\neq c_+\!\in b+Q_+}\,\Q(q,t_\nu) X_c$.
See \cite{Mac} and (6.14) from \cite{C103} or (3.3.14) 
from \cite{C101} (the differential version is due to Opdam):
\begin{align}\label{macde-eval}
&E_{b}(q^{-\rho_k})=
q^{-(\rho_k,b)}
\prod_{\al>0}\prod_{j=1}^{(\al^{\!\vee},b)-1}
\Bigl(
\frac{
1- q_\al^{j}t_\al X_\al(q^{\rho_k})
 }{
1- q_\al^{j}X_\al(q^{\rho_k})
}
\Bigr) \for b\in P_+.
\end{align}
For any $b\in P$, we define $E$\~polynomials as follows:
\begin{align*}
&Y_a(E_b)\,=\,q^{-(a,b+w_b(\rho_k))}E_b, \hbox{\,\, where the 
coefficient of \,} X_b\in E_b \hbox{\, is\, } 1,\\
&\hbox{and\, } w_b\in W \hbox{\, is of 
maximal possible length such that\, } b\in w_b(P_+).
\end{align*}

The technique of intertwining operators for the $E$\~polynomials
results in the following {\,\em existence criterion\,}.
Assuming that $q$ is not a root of unity,
the products $\,den_{\,b}\cdot E_b,\,\, \,den_{\,b}\cdot P_b\,$ and 
$\,num_{\,b}\cdot E_b/E_b(q^{-\rho_k})\,$ for $b\in P_+$
are well defined for the numerator $\,num_{\,b}\,$ and the denominator
$\,den_{\,b}\,$ of the expression for $E_{b}(q^{-\rho_k})$ from
(\ref{macde-eval}). Note that all binomials in $num_{\,b},\, 
den_{\,b}\,$ involve $q$ and $t$ (both).  Also : 
$$
P_{b}(q^{-\rho_k})=\Pi_R^b\,E_{b}(q^{-\rho_k}) \for b\in P_+,\  
\Pi_R^b=\!\!\!\prod_{\al>0,(\al,b)>0}
\frac{
1- t_\al X_\al(q^{\rho_k})
 }{
1- X_\al(q^{\rho_k})
}.
$$
The product $\Pi_R^b$ becomes the {\em Poincar\'e polynomial\,} 
$\Pi_R$ for generic $b\in P_+$. See formulas (6.33), (7.15) from 
\cite{C103} and (3.3.45) from \cite{C101}.

Setting $E_b^\circ=E_b/E_b(q^{-\rho_k})$ for any $b$ and 
$P_b^\circ=P_b/P_b(q^{-\rho_k})$ for $b\in P_+$, the following
relation holds for any $b\in P$:
\begin{align}\label{nontosym}
P_{b_+}^\circ\!=\!(\Pi_R)^{-1}\mathscr{P}
(E^\circ_b)\, \hbox{\, for\, }\, 
\mathscr{P}_+\!\equal\!\sum_{w\in W}
t_{\sht}^{l_{\sht}(w)/2}
t_{\lng}^{l_{\lng}(w)/2}T_w,
\end{align}
where $l_{\sht}, l_{\lng}$ count correspondingly the number 
of $s_i$ for short and long $\al_i$ in any reduced decomposition 
$w=s_{i_l}\cdots s_{i_1}$. We check that the right-hand side here
is proportional to $P_{b_+}$ and then evaluate at $\,q^{-\rho_k}$ by
applying the general formula
$\{T_i (f)\}_{ev}=t_i^{1/2}\{f\}_{ev}\,$ for $i>0$.

The following particular case of (\ref{nontosym})
will be mainly used:
$$
\prod_{\al>0,(\al,b)=0}
\frac{
1- t_\al X_\al(q^{\rho_k})
 }{
1- X_\al(q^{\rho_k})
}
\,P_{b}=\mathscr{P}_+(E_{b}) \for b\in P_+.
$$  
\smallskip

\subsection{\bf DAHA-Jones polynomials}
The following theorem is the first part of Conjecture 2.1 from
\cite{CJ} extended as follows.
We enlarge $PSL_{\,2}^{\wedge}(\Z)$ used
there to $GL_{\,2}^{\wedge}(\Z)$, consider in detail 
the symmetry $(\rr,\ss)\mapsto (\ss,\rr)$ and the
mirror images in (\ref{jones-sym}), and 
add important color-exchange relations (\ref{jones-bc});
last but not least, the justifications are provided.

We represent a given torus knot $T(\rr,\ss)$  
by a matrix $\ga_{\rr,\ss}\in GL_{\,2}(\Z)$ such that its 
first column 
is $(\rr,\ss)^{tr}$ ($tr$ is the transposition),
assuming that \,gcd$(\rr,\ss)=1$; $\rr,\ss$ can be $0$ or negative. 
Then $\hat{\ga}_{\rr,\ss}\in GL_{\,2}^{\wedge}(\Z)$
is by definition a pullback of $\ga_{\rr,\ss}$.

Note that any $\,(\rr,\ss)\,$ can be obviously lifted
to $\,\ga\,$ of determinant $1$ and, accordingly, 
to $\hat{\ga}\,$ from the subgroup $PSL^\wedge_{\,2}$
generated by $\{\tau_{\pm}\}$, i.e. without using
$\eta$. This is sufficient for the construction of
the DAHA-Jones polynomials below; using $GL_2(\Z)$ 
instead of $PSL_2(\Z)\subset GL_2(\Z)$ 
is equivalent to the conjugation relation 
from (\ref{jones-sym}).
\medskip

For a polynomial $F$ in terms of positive and negative 
fractional powers of $q$ and $t_\nu$, 
the {\em tilde-normalization}
$\tilde{F}$ will be the result of division by the lowest 
$q,t_\nu$\~monomial (assuming that it is well defined);
$\tilde{F}$ contains only non-negative powers of 
$q,t_\nu$. In the following theorem, the lowest monomials always 
exist and the powers are integral.

\medskip
\begin{theorem} \label{MAINTHM} 
(i) {\sf Polynomiality}.
Given $T(\rr,\ss)$, a root system $R$ and a weight 
$b\in P_+\,$, we set 
\begin{align}\label{jones-d}
J\!D_{\rr,\ss}^{R}(b\,;\,q,t)=
J\!D_{\rr,\ss}(b\,;\,q,t)\equal 
\{\,\hat{\ga}_{\rr,\ss}(P_b)/P_b(q^{-\rho_k})\,\}_{ev}.
\end{align}
Then $\tilde{J\!D}_{\rr,\ss}(b\,;\,q,t)$, the tilde-normalization
of the (reduced) DAHA-Jones polynomial,  does not depend on the 
particular choice of $\ga_{\rr,\ss}\in GL_2(\Z)$
and $\hat{\ga}_{\rr,\ss}\in GL_{\,2}^{\wedge}(\Z)$.
It is a polynomial in terms of $q,\{t_\nu\}$. 
Moreover, using the conjugation involution $\star$ from 
(\ref{eta-conj}),
\begin{align}\label{jones-nonsym} 
&J\!D_{\rr,\ss}^{R}(b\,;\,q,t)\!=\! 
\{\,\hat{\ga}_{\rr,\ss}
(E_c)/E_c(q^{-\rho_k})\,\}_{ev} \hbox{\, for\, }
c\in W(b),\\
\label{jones-sym} 
&\tilde{J\!D}_{\rr,1}(b\,;\,q,t)=1,\ \
J\!D_{\rr,-\ss}(b\,;\,q,t)=\bigr(J\!D_{\rr,\ss}
(b\,;\,q,t)\bigl)^\star,
\notag\\
&J\!D_{\rr,\ss}(b\,;\,q,t)\,\ =\ 
\,J\!D_{\ss,\rr}(b\,;\,q,t)\,\ =\ \,
J\!D_{-\ss,-\rr}(b\,;\,q,t).
\end{align}

(ii) {\sf Color Exchange.} Let us assume that 
for $b,c\in P_+$ and certain \,id\,$\neq w\in W$,  
\begin{align}\label{bcw-rel}
q^{\,(b+\rho_k-w(\rho_k)-w(c)\,,\,\al)\,}=1=
q_{\al}^{\,(b-w(c)\,,\,\al^{\!\vee})}\,
t_{\sht}^{\,(\rho^w_{\sht}\,,\,\al^{\!\vee})}\,
t_{\lng}^{\,(\rho^w_{\lng}\,,\,\al^{\!\vee})}
\end{align}
for any $\al\in R_+$, where $q$ is not a root of unity and
$\rho_{\nu}^w\equal w(\rho_{\nu})-\rho_{\nu}$.
For instance, taking $P_+\ni b=w(\rho_k)-\rho_k+w(c)$ 
for a given $c\in P_+$ is sufficient. 
The product in the right-hand side is in terms of
integral powers of $q_\nu$ and $t_{\nu}^{l_{\nu}}$, where
$\rho_{\nu}^w\in l_\nu P$ for $l_\nu\in \Z_+$.
Then 
\begin{align}\label{jones-bc}
J\!D_{\rr,\ss}(b\,;\,q,t)=
J\!D_{\rr,\ss}(c\,;\,q,t)
\hbox{\, for such \, } q,\{t_\nu\} \hbox{\, and any\, } \rr,\ss.
\end{align}
In particular, this holds for $w=s_i(i>0), c\in P_+$ and for
\begin{align}\label{jonbcsi}
&b=c-\bigl(k_i+(c,\al_i^{\!\vee})\bigr)\al_i  
\hbox{\, provided \, }
2k_i+(c,\al_i^{\!\vee})\in -\Z_+ \hbox{\, and}\\
(c,\al_j^{\!\vee})&+\nu(k_i+(c,\al_i^{\!\vee}))\in \Z_+
\hbox{\,\, for all\, } j>0 \hbox{ such that } (\al_i,\al_j)<0,\notag
\end{align}
where $\nu=\nu_i$ for short $\al_j$ and 
$\nu=1$ otherwise
(thus $k_i\in -\Z_+$ unless for $A_1$ and if $\,\al_i\,$ is long
for $\,C_{n\ge 2}$). Another particular case of (\ref{bcw-rel}) 
is $w=w_0$ with $\,2k_\nu\in -\Z_+$, where 
\begin{align}\label{jonbcwo}
b=-2\rho_k+w_0(c) \hbox{\, provided \, }
(c,\al_i^{\!\vee})\le -2k_i \hbox{\, for all\, } i>0.
\end{align}

(iii) {\sf Evaluation}. Relation (\ref{bcw-rel}) 
results in (\ref{jones-bc})
for any $\,q$, including roots of unity, due to the polynomiality
of $J\!D_{\rr,\ss}(b\,;\,q,t)$. Upon $\,q\to 1$, one has 
$\,E_{w_0(b)}=P_{b}$ and $\,P_{b+c}=P_{b} P_{c}$ 
for $\,b,c\in P_+$. Accordingly, %for $b=\sum_{i=1}^n b_i \om_i$,
\begin{align}\label{jones-eval}
\tilde{J\!D}_{\rr,\ss}\bigl(\sum_{i=1}^n b_i \om_i\,;\,q\!=\!1,t
\bigr) =\!
\prod_{i=1}^n\tilde{J\!D}_{\rr,\ss}(\om_i\,;\,q\!=\!1,t)^{b_i}
\hbox{\, for any \,} \rr,\ss.
\end{align}
 \end{theorem}

{\it Proof.}
By construction,
$\tilde{J\!D}_{\rr,\ss}^{R}(b\,;\,q,t)$ depends only 
on the first column 
of the matrix $\ga=\ga_{\rr,\ss}$. The switch from
$P_b^\circ =P_b/P_b(q^{-\rho_k})$ to 
$E_b^\circ =E_b/E_b(q^{-\rho_k})$ is straightforward, 
since one can perform the $t$\~symmetrization inside 
(\ref{jones-nonsym}) using formula (\ref{nontosym}).
Let us deduce the symmetries in  (\ref{jones-nonsym})
from (\ref{evsym}). First of all, 
(\ref{iotaXY}) and (\ref{evsym}) result in 
$\,\iota\cdot \hat{\ga}=
\hat{\ga}\cdot \iota\,$ for any $\hat{\ga}\in GL_{\,2}^{\wedge}(\Z)$
and therefore provide 
\begin{align}\label{jdiota}
J\!D_{\rr,\ss}^R(\iota(b)\,;\,q,t)=
J\!D_{\rr,\ss}^{R}(b\,;\,q,t).
\end{align}
Then, using the $\vph$\~invariance of $\,\{\cdot\}_{ev}$,
\begin{align}\label{symminus}
\{ \vph \bigl(\cdots \tau_+^\be\tau_-^\al(P^\circ_b)\bigr)\}_{ev} =
\{ \cdots \tau_-^\be\tau_+^\al \bigl(P^\circ_b(Y^{-1})\bigr)\}_{ev}
\\
=
\{ \cdots \tau_-^\be\tau_+^\al\si\bigl(P^\circ_b(X)\bigr)\}_{ev}=
\{ \bigl(\cdots \tau_-^\be\tau_+^\al\si\bigr)(P^\circ_b)\}_{ev},
\notag
\end{align}
which proves that 
$J\!D_{\rr,\ss}(b\,;\,q,t)=
J\!D_{-\ss,-\rr}(b\,;\,q,t).$
Similarly,
\begin{align*}
&\{ \eta \bigl(\cdots \tau_+^\be\tau_-^\al(P^\circ_b)\bigr)\}_{ev}=
\{ \cdots \tau_+^{-\be}\tau_-^{-\al}
\bigl((P^\circ_b)^\star\bigr)\}_{ev}\\
=
&\{ \cdots \tau_+^{-\be}\tau_-^{-\al}
\bigl(P^\circ_{\iota(b)})\bigr)\}_{ev}=
\{ \cdots \tau_+^{-\be}\tau_-^{-\al}\bigl(P^\circ_{b})\bigr)\}_{ev}=
\{\cdots \tau_+^{\be}\tau_-^\al(P^\circ_b)\}_{ev}^\star\,,
\end{align*}
which results in $J\!D_{\rr,\ss}(b\,;\,q,t)=
J\!D_{\rr,-\ss}(b\,;\,q,t)^\star.$ Combining this symmetry with the
previous one,
\begin{align*}
&J\!D_{\rr,\ss}(b\,;\,q,t)=J\!D_{-s,r}(b\,;\,q,t)^\star=
J\!D_{-\rr,\ss}(b\,;\,q,t)^\star\\
&= J\!D_{-\rr,-\ss}(b\,;\,q,t)=
J\!D_{\ss,\rr}(b\,;\,q,t).
\end{align*}
The latter symmetry can be 
established directly following (\ref{symminus}):
\begin{align*}
&\{ \vph \bigl(\cdots \tau_+^{\be}\tau_-^{\al}
(P^\circ_b)\bigr)\}_{ev}=
\{ \cdots \tau_-^{\be}\tau_+^\al
\bigl(P^\circ_b(Y^{-1})\bigr)\}_{ev}\\
=&
\{ \cdots \tau_-^{\be}\tau_+^{\al}\si^{-1}\bigl
(P^\circ_{\iota(b)}(X)\bigr)\}_{ev}=
\{ \cdots \tau_-^{\be}\tau_+^{\al}\si^{-1}\bigl
(P^\circ_{b}(X)\bigr)\}_{ev}\,,
\end{align*}
where we use (\ref{sisquare}) and (\ref{jdiota}).
\smallskip

Finally for $b\in P_+$,
\begin{align*}
&J\!D_{\rr,1}(b\,;\,q,t)=J\!D_{1,\rr}(b\,;\,q,t)
\!=\!\{\,\tau_-^{\rr}(P_b^\circ)\,\}_{ev}\!=\!
\{\,\tau_-^{\rr}\cdot P_b^\circ\cdot \tau_-^{-\rr}(1)\mid_{{}_\v})
\,\}_{ev}
\\
&=\,\{\,\tau_-^{\rr}(P^\circ_b)\mid_{{}_\v}\,\}_{ev}
\,=\,\{\,q^{-\rr(b,b)/2-\rr(b,\rho_k)}P_b^\circ\,\}_{ev}\,=\,
q^{-\rr(b,b)/2-\rr(b,\rho_k)},
\end{align*}
where $\tau^{\pm \,r}_-(\cdot)\mid_{{}_\v}$ means
the action in $\v$; see 
formula (\ref{tau-min-e}) below.
\medskip

{\sf Polynomiality of $J\!D$}.
We will use the {\em radical\,}, $Rad$, of the
{\em evaluation pairing\,}  defined as follows: 
$$
\{E,F\}_{ev}\,=\, E(Y^{-1})(F(X))(t^{-\rho}),\ E,F \in \v\,.
$$
Theorem 11.8 from \cite{C103} gives the necessary and
sufficient conditions for  $Rad\neq \{0\}$
for generic $q$. For instance (though we do not really need this),
the radical is nonzero in the case of $ADE$ if and only if
$$
t\!=\! q^{-l -\!\frac{j}{m_i\!+\!1}}\,\ze_i^{j\,'}
\hbox{\, for \,} 1\le i\le n,\ 
0\le j,j\,'\le m_i,\, j\!+\!j\,'\!>\!0, \,\, l\! \in \!\Z_+, 
$$ 
where  $t=t_{\sht}$, $\ze_i$ are primitive $(m_i+1)${\small th}\,
roots of unity for the classical {\em exponents\,} 
$m_i$ (see \cite{Bo}).  This is 
Theorem 11.1 from \cite{C103}.

\comment{
Note that the exponent 
$n-1$ appears twice in the list of $\{m_i\}$ for
the root systems $D_{\hbox{\tiny even}}$. We will exclude 
this case from our further considerations. The polynomiality
of  $J\!D$ for $D_{n}$ will follow from that in the
case of $C_n$ upon the reduction $t_{\lng}=1$;
see formula (4.7) from \cite{CJ}. 
}

\smallskip
If  $J\!D_{\rr,\ss}(b;q,t)$ is not a polynomial for $b\in P_+$
and admissible  $\,\rr,\ss$, then $E_b^\circ=E_b/E_b(q^{-\rho_k})$ 
has a pole at $\ep=0$ of order $l>0$ for $\ep$ equal 
to one of the binomials $(1-q_\al^j t_\al^p)$ in the numerator
$num_{\,b}\,$ of (\ref{macde-eval}) for certain $\,j>0,p>0$. 
Recall that
the inequality $\,j>0$ here is due to using  
$E$\~polynomials; see (\ref{jones-nonsym}). Here and further
we will localize and complete the ring of constants 
$\Z_{q,t}=\Z[q^{\pm 1/\mm},t_\nu^{\pm 1/2}]$ with respect to $\ep$,
i.e. for the principle ideal $(1-q_\al^j t_\al^p)$; the notations 
will be $\Z_{q,t}^{(\ep)}, \v^{(\ep)}$. For the sake of definiteness, 
we will take a maximal binomial here, i.e. such that  
$(1-q_\al^{jv} t_\al^{pv})$ is not in $num_{\,b}\,$ for any 
$\Z\ni v>1$.
Accordingly, we introduce the following filtration of submodules
of $\v$:
\begin{align}\label{radep}
Rad_{\ep,\ell}=\bigl\{F\in \v\,\mid\, \{F,\v\}_{ev}\in 
\ep^\ell\Z_{q,t}\bigr\}
\for \ell\in \Z_+\, 
\end{align}
setting $Rad(\ep=0)=Rad_{\ep,\infty}$; we 
add $q^{1/(2\mm)}$ to $\Z_{q,t}^{(\ep)}$ here and below.

Note that it is generally not impossible that such 
$(1-q_\al^j t_\al^p)$ coincides with one of the binomials in 
the denominator $den_{\,b}\,$ of (\ref{macde-eval}). 
As a matter of fact, this is not the case due to the consideration 
below, but we do not really need this fact.
\smallskip

The polynomials  $E_b'=\ep^l E_b^\circ,\, 
P_b'=\ep^l P_b^\circ \in \v^{(\ep)}$
are eigenfunctions respectively
for $\{Y_a\}$ and $\{L_f\}$ from (\ref{macdopers}). Moreover, 
$E_b'(q^{-\rho_k})\in \ep^l\Z_{q,t}^{(\ep)}\ni P_b'(q^{-\rho_k})$.   
%We impose no other relations on $q,t_\nu$ at 
%$\ep=0=(1-q_\al^j t_\al^p)$.  
Therefore $E_b' \in Rad_{\ep,l}\ni P_b'$; see e.g., 
Lemmas 11.4-5 from \cite{C103}.
Thus $Rad(\ep=0)\neq \{0\}$ and, for instance, 
Theorem 11.8 there implies that $(1-q_\al^j t_\al^p)$
does not coincide with the binomials in $den_{\,b}\,$
from (\ref{macde-eval}) (as it was claimed above).
\smallskip

Let us now switch from $Rad\,$ to  
\begin{align}\label{RADep}
R\!A\!D_{\ep,\ell}\equal\bigl\{ H\in \HH\, \mid\, 
\{\HH\, H\, \HH\}_{ev}\in \ep^\ell\Z_{q,t}\bigr\}
\for \ell\in \Z_+.
\end{align}
We set $R\!A\!D(\ep=0)=R\!A\!D_{\ep,\infty}$.
Equivalently, one has:
$$
R\!A\!D_{\ep,\ell}=\bigr\{H\in \HH \mid 
H(\v)\subset Rad_{\ep,\ell}\bigl\},
$$ 
since $Rad_{\ep,\ell}=\bigl\{F\in \v \,\mid\, \{\HH(F)\}_{ev}\in 
\ep^\ell\Z_{q,t}^{(\ep)}\bigr\}$; see Lemma 11.3, \cite{C103}.
%\smallskip

Here $q$ is not a root of unity. Therefore Proposition 3.2
from \cite{C100} states that any $Y$\~invariant
submodule of $\v$ is invariant with respect to the
natural action of $\tau_-$ in $\v$. For instance,
\begin{align}\label{tau-min-e}
&\tau_-(E_b)=q^{-(b,b)/2-(b,\rho_k)}E_b \hbox{\, and\, }
\tau_-(P_b)=q^{-(b,b)/2-(b,\rho_k)}P_b,
\end{align}
assuming that $E_b$ for $b\in P_+$ is well defined. See also
Proposition 3.3.4 from \cite{C101}. 

Therefore
$\psi$ and $\tau_-$ preserve $R\!A\!D_{\ep,r}$ for any $r\in \Z_+$
(and generic $q$), as well as $\eta$. Thus the whole 
$GL_{\,2}^{\wedge}(\Z)$ fixes each $R\!A\!D_{\ep,r}$. This implies that
$\hat{\ga}(P'_b)\in R\!A\!D_{\ep,l}$ and 
$\{\hat{\ga}(P_b^\circ)\}_{ev}$
is divisible by $\ep^l$. Hence 
$\tilde{J\!D}_{\rr,\ss}(b\,;\,q,t)$ has no singularity
at $\ep=0$, which contradiction is sufficient to claim 
polynomiality of $\tilde{J\!D}_{\rr,\ss}(b\,;\,q,t)$
for any $b\in P_+$ and $\rr,\ss$.
\medskip

{\sf Parts $(ii,iii)$}.
The justification of Part $(ii)$ of the theorem is quite similar.
Let $\ep$ be $1$ minus the right-hand side of 
(\ref{bcw-rel}) and $l\in \Z_+$ is minimal such that
$\ep^l\,P_b^\circ$ and $\ep^l\,P_c^\circ$ are regular in 
$\v^{(\ep)}$; as above,   $\v^{(\ep)}$ is defined over 
$\Z_{q,t}^{(\ep)}$.  
The conditions from (\ref{bcw-rel}) are necessary
and sufficient for these polynomials to have coinciding 
eigenvalues for any $L_f$ from (\ref{macdopers}). 
This coincidence results in the relations 
$\ep^l(P_b^\circ-P_c^\circ)\in Rad_{\ep,l+1}$, which 
gives the required. 

Note that given $c,q,t_\nu$, 
there always exists $b=\in P_+$ satisfying (\ref{bcw-rel})
such that the polynomial $P_b^\circ$ is regular at $\ep=0$. 
The regularity of $J\!D_{\rr,\ss}(b\,;q,t)$ from $(i)$
is granted automatically for such $b$. Moreover one such $b$
can be canonically constructed in terms of the 
{\em right Bruhat ordering\,} 
defined for the root subsystem of $\tR$ associated
with the stabilizer of $\rho_k$ in $\hW$. It was called
{\em primary\,} in  \cite{C103}; see formula (9.2) 
from Section 9.1 there.
\smallskip

It is not really necessary to assume in (\ref{bcw-rel}) that
$q$ is not a root of unity, since we have the polynomiality
of the $J\!D$\~polynomials.
Let us make $q=1$. Then the technique of intertwining
operators of $E$\~polynomials readily results in the formula
$E_{w_0(b)}=P_b$ for $b\in P_+$ and the multiplicative
property $P_b P_c=P_{b+c}$ for $b,c\in P_+$. We use that
$E_{s_i(b)}$ is $(s_i+1)E_b$ for $b\in P$ assuming that
$(\al_i,b)>0$ for $i>0$. Also,  $E_{b+u\vth}=
X_{\vth}s_\vth(E_b)+E_b$ if $u=1-(\vth,b)>0$ and
$E_{b+\om_r}=X_{\om_r}\pi_r(E_b)$ for $r\in O'$.
See e.g., Proposition 3.3.5 from \cite{C101}. Therefore 
\begin{align}\label{gaPqone}
\hat{\ga}(P_b)=\prod_{i=1}^n \hat{\ga}(P_{\om_i})^{b_i} \for
b=\sum_{i=1}^n b_i\om_i\in P_+.
\end{align}
We need to apply this operator to $1\in \v$. 
Firstly, $L_f$ is a constant, which is
$\{f\}=f(t_{\sht}^{\rho_{\sht}}t_{\lng}^{\rho_{\lng}})$, 
when acting on symmetric polynomial $F\in \v^W$
for any symmetric $f$ due to $q=1$; see (\ref{macdopers}).
Secondly, the elements $f(X)$ and $L_f$ are {\em central\, }
in $\HH$ if $q^{1/\mm}=1$, as well as any operators
$\hat{\ga}(f(X))$ for symmetric $f$ and any 
$\ga\in GL_{\,2}^\wedge(\Z)$. See Theorem 3.8.5 from \cite{C101}.
The action of $GL_{\,2}^\wedge(\Z)$ in the center of $\HH$ at
$q=1$ and at the roots of unity is important in DAHA theory
and have geometric applications.
For instance, we have that $f(Y)=\{f\}$ in the whole $\v$ for
symmetric $f$ and $q^{1/\mm}=1$. 
Thirdly and finally, 
\begin{align*}
&\hat{\ga}(P_b)(F\in \v^W)\ =\
\prod_{i=1}^n \,\bigl(\hat{\ga}(P_{\om_i})(F)\bigr)^{b_i}
\for \ga=\ga_{\rr,\ss} \hbox{\,\, and}\\
&J\!D_{\rr,\ss}(b\,;q\!=\!1,t)=
\bigl\{\bigl(\hat{\ga}(P_b)(1\!\in\! \v^W)\bigr\}=
\prod_{i=1}^n \bigl\{\hat{\ga}(P_{\om_i})(1)\bigr\}^{b_i}.
\end{align*}
\vskip -1.5cm
\sq
\vskip 0.2cm

We mention that Part $(ii)$ leads to 
{\em resonance $\tilde{J\!D}$\~polynomials\,}, which are
the limits of linear combinations of 
$\tilde{J\!D}_{\rr,\ss}(b\,;\,q,t)$ for the same 
$q,t_\nu$ and for $b$ from $(ii)$ corresponding to
a given $c$. We divide such linear combinations
by the leading powers of $\,\ep\,$ before taking the 
limit $\ep\to 0$; $\ep$  is $1$ minus the right-hand side of 
(\ref{bcw-rel}). They are related to the {\em nonsemisimple
Macdonald polynomials\,} from \cite{C103} in the spherical
normalization, but we will not touch this upon in the present
paper.
\medskip

\subsection{\bf Superpolynomials}
Theorem \ref{MAINTHM}  has the following extension
to the {\em DAHA- superpolynomials\,}
in type $A_n$. The following stabilization theorem
was announced in \cite{CJ}; it was mentioned there
(Section 2.4.1, ``Confirmations")
that its justification is similar to Lemma 4.3 and 
formula (4.1) from \cite{SV}. The complete proof
is published in \cite{GoN} (following \cite{SV}). 
We switch to $A_n$ and naturally set $t=t_{\sht}=q^k$.

\begin{theorem}\label{STABILIZ}
(i) For the root system $A_n$, let
us consider $P_+\ni b=$ $\sum_{i=1}^n b_i \om_i$ as 
a (dominant) weight for any $A_N$ with $N\ge n-1$,
where we set $\om_{n}=0$ upon the restriction to $A_{n-1}$.
Then given $T(\rr,\ss)$, there exists $\h_{\rr,\ss}(b\,;\,q,t,a)$,
a polynomial in terms of $a,q,t^{\pm 1}$, such that
its coefficient of $a^0$ is tilde-normalized 
(i.e. in the form $\sum_{u,v\ge 0}C_{u,v}q^u t^v$ 
with $C_{0,0}=1$) and the following specializations hold:
\begin{align}\label{jones-sup}
\h_{\rr,\ss}(b\,;\,q,t,a=-t^{N+1})=
\tilde{J\!D}_{\rr,\ss}^{A_{N}}(b\,;\,q,t)
\hbox{\, for any\, } N\ge n-1. 
\end{align}

(ii) Imposing the Color Exchange relation (\ref{bcw-rel}), 
we will consider $\,w\,$ there
%from Part (ii) of Theorem  \ref{MAINTHM} 
as an element of $\S_{N+1}$ for every
$N\ge n$ (the Weyl group for $A_N$) naturally acting in 
the corresponding $P$. Then given \, $\rr,\ss$\, and
up to a proper power $\,q^{\bullet}t^{\bullet}$,
\begin{align}\label{jones-bca}
\h_{\rr,\ss}(b\,;\,q,t,a) = q^{\bullet}t^{\bullet}\,
\h_{\rr,\ss}(c\,;\,q,t,a)
\hbox{\, for such \, } q,\{t_\nu\}.
\end{align}
In particular, let $w=s_i=(i,i+1)$ with $i<n$. Then for
a dominant $c$\, and $\,b=c-\bigl(k+(c,\al_i)\bigr)\al_i$,
the components of $\,b$ are 
\begin{align}\label{jonacsi}
b_i\!=\!-2k\!-\!c_i,\ 
b_j\!=\!c_j\!+\!c_i\!+\!k \hbox{ for } 
j\!=\!i\!\pm\!1\!>\!0,\ b_j\!=\!c_j \hbox{ otherwise.}
\end{align}
Here $k\!\in\! -\Z_+$ satisfies the 
relations $\,c_i/2\!\le\! -k\!\le\! c_i+\min\{c_{i\pm 1>0}\}\,$,
which are necessary and sufficient for $b\in P_+$.

(iii) Making $q=1$, one has:
\begin{align}\label{jones-seval}
\h_{\rr,\ss}(b\,;\,q\!=\!1,t,a)\!=\!
\prod_{i=1}^n\h_{\rr,\ss}(\om_i\,;\,q\!=\!1,t,a)^{b_i}
\hbox{\, for\, } b\!=\!\sum_{i=1}^n b_i \om_i.
\end{align}
This gives that the $\,a$\~degree\,
deg${}_a(\h_{\rr,\ss}(b\,;\,q,t,a)\,$ of 
$\,\h_{\rr,\ss}(b\,;\,q,t,a)$\,
equals $\,\min(|\rr|,|\ss|)\,$ times
the number of boxes in the Young diagram $\la_{b}$
associated with $b\in P_+$. These claims are from
Conjecture 2.6 of \cite{CJ}. 
\end{theorem} 
\vskip -1.25cm \sq
\vskip 0.2cm

The conjectural relation of DAHA-superpolynomials 
to the stable {\em 
Khovanov- Rozansky polynomials\,} for $sl_{n\!+1}$
is as follows: $a\mapsto\! t^{n\!+1}\sqrt{t/q}$.
Note that this is not the substitution from (\ref{jones-sup}),
but they are linked. See Section 2.3 in  \cite{CJ} for
discussion and some references; in particular, 
formula (2.12) there connects our (DAHA -based)
parameters with the standard ones. Such relations match 
the conjectural links between the Khovanov- Rozansky 
homology \cite{KhR1,KhR2,Ras} and the superpolynomials 
introduced via the {\em BPS states\,} \cite{DGR,AS,FGS,GGS} as 
well as those obtained in terms of {\em rational DAHA\,},
which are deeply related to certain {\em Hilbert schemes\,} 
\cite{GORS,GoN}.

The coincidence was confirmed when the stable Khovanov-Rozansky
polynomials are available and
for $\h_{\rr,\ss}(b\,;\,q,q,-a)$, which do coincide 
with the corresponding {\em HOMFLYPT polynomial\,}.
See Proposition 2.3 from \cite{CJ},\cite{Ste} and
\cite{RJ} (concerning 
the Jones polynomials and Quantum groups) and references 
there.
When $n\!+\!1\!=\!0$, the relation to the 
{\em Heegard- Floer homology\,} of torus knots is 
expected.

\medskip

{\sf Color exchange combinatorially.}
We associate with
$c=\sum_{i=1}^n c_i\om_i$ in Part $(ii)$ the Young diagram
$$
\la_c=\{m_1\!=\!c_1+\ldots+c_n,m_2\!=\!c_1+\ldots+c_{n-1},
\ldots, m_n\!=\!c_n,0,0,\ldots\}.
$$
Then we switch to $\la_c'=\{m'_i=m_i-k(i-1)\}$, apply $w\in W$ to
$\la_c'$  and finally obtain 
\begin{align}\label{lamw}
\la_b=\{m'_{w(i)}+k(i-1),\}=\{m_{w(i)}+k(i-w(i))\}.
\end{align} 
Here 
$w$ transforms the rows of $\la$ and we set 
$w\{m_1,m_2,\ldots,m_n\}=\{m_{w(1)},m_{w(2)},\ldots,m_{w(n)}\}$.
Given $k<0$ (it can be fractional), $\la_b$ 
{\em must be a Young diagram\,},
which determines the set of all $\,w\,$ that can be used.
Note that when $k=+1$ (formally), the procedure $\la_c\mapsto \la_b$ 
is actually from the Jacobi-Trudi formula. 

Let us briefly discuss the Duality Conjectures 2.5 from \cite{CJ}.  
It states that 
\begin{align}\label{dualityc}
\h_{\rr,\ss}(\la\,;q,t,a)=q^{\bullet}t^{\bullet}
\h_{\rr,\ss}(\la^{tr}\,;t^{-1},q^{-1},a),
\end{align}
where we switch from weights to the corresponding Young diagrams
and $\la^{tr}$ is the transposition of $\la$.
This was justified in \cite{GoN} using the {\em modified 
Macdonald polynomials\,}. 
Let us outline a direct proof (which is expected
to work for $C^\vee C_n$). We use the 
{\em perfect repsentation\,} at $t=q^{-(s+1)/(n+1)}$,  
which is defined as $\v/Rad$ for $\,s\in \Z_+$ provided\,
gcd$(s+1,n+1)=1$. More generally, 
the Coxeter number must be used here instead of $(n+1)$. 
See \cite{C101}. Then $\,a=-t^{n+1}=-(q^{-1})^{s+1}\,$ 
and we can identify the left-hand side of (\ref{dualityc}) 
for $A_n$ with the right-hand side for $A_s$. This identification 
goes via the theory at roots of unity; we compare the zeros in both 
sides of (\ref{dualityc}).
\smallskip

{\sf More on Part $(iii)$.} Formula (\ref{jones-seval}) follows
directly from (\ref{jones-eval}). Let us mention that the 
formulas for $\h_{\rr,\ss}(\om_i\,;\,q\!=\!1,t,a)$
are not difficult to calculate for simple knots, 
but generally they are involved (even for $q=1$).
Note that (\ref{jones-eval}) 
combined with the duality results in
\begin{align}\label{jones-sdeval}
\h_{\rr,\ss}(\la\,;\,q,t\!=\!1,a)=
\prod_{i=1}^n\h_{\rr,\ss}(m_i\om_1\,;\,q,t\!=\!1,a)^{m_i},
\end{align}
where $m_i$ is the number of boxes in the $i${\small th} row of
$\la$. A direct justification of (\ref{jones-sdeval}) 
without using the duality is not clear at the moment.

Concerning the $a$\~degrees there, it is relatively
straightforward to justify that 
deg${}_a(J\!D_{\rr,\ss}(\la\,;\,q,t,a)\le 
\min(|\rr|,|\ss|)\deg(\la)$;
see e.g. \cite{GoN}. Then the evaluation
formula at $q=1$  (or even that for $\,q=1=t$) gives that we
actually have the equality. 
\smallskip

%So the terms in the (formal) determinant of the Jacobi-Trudi matrix 
%that are Young diagrams have superpolynomials coinciding (for $qt=1$)
%with that for the diagram $\la_c$ (treated as the diagonal of the 
%Jacobi-Trudi matrix). 

\subsection{\bf Examples of color exchange}
To begin with, let $i=1,k=-2,$ $c=\om_1+\om_2, b=3\om_1$ 
in Part $(ii)$ of the theorem. 
According to Section 3.1 from \cite{CJ},
\renewcommand{\baselinestretch}{0.5} 
$$
\h_{3,2}(\om_1+\om_2;\,q,t,a)\ =
$$
{\small
\noindent
\(
1+\frac{a^3 q^6}{t}+2 q t-q t^2+2 q^2 t^2+q^3 t^2-q^2 t^3
+2 q^3 t^3-q^3 t^4+2 q^4 t^4+q^5 t^5
\)

\noindent
\(
+a \bigl(2 q^2+q^3+\frac{q}{t}-q^2 t+3 q^3 t+q^4 t-q^3 t^2
+3 q^4 t^2+q^5 t^2-q^4 t^3+2 q^5 t^3+q^6 t^4\bigr)
\)

\noindent
\(
+a^2 \bigl(q^4
+q^5+\frac{q^3}{t}+\frac{q^4}{t}-q^4 t+q^5 t+q^6 t+q^6 t^2\bigr)
\)\,.

}
\renewcommand{\baselinestretch}{1.2} 
\smallskip

This is, by the way, the simplest example of a
DAHA-superpolynomial with negative coefficients.
Conjecturally, the negative terms are not present only 
for rectangle Young diagrams; see \cite{CJ}. This is in 
obvious contrast with formula (4.14) from \cite{GGS} for 
the {\em 3-hook\,}. It resembles our one (both
have the same number of terms), but 
has no negative terms. The  $t$\~powers are odd and 
even there (they are even in our one upon using the standard 
parameters $q,t,a$). Formula (4.14) is of course a suggestion,
not the result of a formal calculation.
\smallskip   

The second superpolynomial is
\renewcommand{\baselinestretch}{0.5} 
$$
\h_{3,2}(3\om_1;\,q,t,a)\ =
$$
{\small
\noindent
\(
1 + a^3 q^{12} + q^3 t + q^4 t + q^5 t + q^6 t^2 + q^7 t^2 + 
q^8 t^2 + q^9 t^3
\)

\noindent
\(
+a \bigl(q^3 + q^4 + q^5 + q^6 t + 2 q^7 t + 2 q^8 t + q^9 t + 
q^9 t^2 + q^{10} t^2 + q^{11} t^2\bigr)
\)

\noindent
\(
+a^2 \bigl(q^7 + q^8 + q^9 + q^{10} t + q^{11} t + q^{12} t\bigr)
\)\,.

}
\renewcommand{\baselinestretch}{1.2} 
\smallskip

Combinatorially, the procedure is as follows:

\vskip 0.2cm
\begin{thicklines}
\raisebox{-1.0ex}{\hbox{\ $\la_c=$\ }}
{\framebox(10,10){}}$\kern -2.95pt$
{\framebox(10,10){}} 
\raisebox{-1.0ex}{\hbox{\ $\overset{\!\!+\!k\rho'}{\longmapsto}$\ }}
\ {\framebox(10,10){}}$\kern -2.95pt$
  {\framebox(10,10){}} 
\ \ \ \ \raisebox{-1.0ex}{\hbox{\ $\overset{w}\longmapsto$\ }}
{\framebox(10,10){}}$\kern -2.95pt$
{\framebox(10,10){+}}$\kern -2.95pt$
{\framebox(10,10){+}} \ \ 
\raisebox{-1.0ex}{\hbox{\ $\overset{\!\!-\!k\rho'}\longmapsto$\ }}
{\framebox(10,10){}}$\kern -2.95pt$
{\framebox(10,10){}}$\kern -2.95pt$
{\framebox(10,10){}} 
\raisebox{-1.0ex}{\hbox{\ $=\la_b.$\ }}\hfill

%\vskip -0.6cm
\vskip -0.25cm
%\thicklines
{%$\kern -1.4pt$ 
$\kern +35.4pt$\framebox(10,10){}}
\ \ \ \ \ \ \ \ \ \
\ $\kern -2.0pt$ 
  $\kern -2.95pt$
{\framebox(10,10){}}$\kern -2.95pt$
{\framebox(10,10){+}}$\kern -2.95pt$
{\framebox(10,10){+}}%$\kern -2.95pt$ 
$\kern 0.5pt$ 
\ \ \ \ \ \ \ \ \ {\framebox(10,10){}}$\kern -2.95pt$
{\framebox(10,10){}} \hfill
\end{thicklines}
\vskip 0.2cm

%$T_\mathsf{r,s},T_\mathit{r,s}, T_\mathtt{r,s}$

Here $-\rho'$ is a stable
re-normalization of  $-\rho$ 
given by the diagram $\{0,1,2,3,\cdots\}$. The boxes
added due to the translation by $k\rho'=-2\rho'$ 
are marked by $+$\,. Generally, we can add fractional
boxes here if $k$ is fractional (recall that it is always negative). 
\smallskip

We obtain that
\begin{align}\label{21to3}
q^4\,\h_{3,2}(\om_1+\om_2;\,q,t,a) =
\h_{3,2}(3\om_1;\,q,t,a) \for t=q^{-2}.
\end{align}
The corresponding {\em resonance $J\!D$\~polynomial\,} 
is  
$$
\lim_{\ep\to 0}\bigl(tq^6\,\h_{3,2}(\om_1\!+\!\om_2;\,q,t,a)\!-\!
\h_{3,2}(3\om_1;\,q,t,a)\bigr)/\ep \for \ep=1-tq^2.
$$ 

\smallskip
Since $\h_{3,2}(\om_1+\om_2;\,q,t,a)$
contains negative terms,  a remarkable
cancelation of monomials occurs in this polynomials 
upon $t=q^{-2}$.
Relation (\ref{21to3}) holds for 
any $\rr,\ss$ with a proper coefficient of proportionality, 
which is $q^{12}$ for the knot $T(3,4)$. 
\smallskip

Let us take now $c=\om_1+\om_2, w=s_1 s_2 s_1, k=-1/2$. Then
$b=\om_3$ and we obtain that\, 
$t^{\bullet}\,\h_{\rr,\ss}(\om_1+\om_2;\,q,t,a) =
\h_{\rr,\ss}(\om_3;\,q,t,a)$ for $q=t^{-2}$; the coefficient
of proportionality is $t^4$ for $T(3,2)$.

The corresponding combinatorial transformation is as follows:
\vskip 0.2cm
\begin{thicklines}
\raisebox{-1.5ex}{\hbox{\ $\la_c=$\ }}
{\framebox(10,10){}}$\kern -2.95pt$
{\framebox(10,10){}} 
\raisebox{-1.5ex}{\hbox{\ $\overset{\!\!+\!k\rho'}{\longmapsto}$\ }}
\ {\framebox(10,10){}}$\kern -2.95pt$
  {\framebox(10,10){}} 
\ \ \raisebox{-1.5ex}{\hbox{\ $\overset{w}\longmapsto$\ }}
{\framebox(10,10){+}}
\ \ \ \ 
\raisebox{-1.5ex}{\hbox{\ $\overset{\!\!-\!k\rho'}\longmapsto$\ }}
\ $\kern -0.7pt$
{\framebox(10,10){+}}
\raisebox{-1.5ex}{\hbox{\ $=\la_b.$\ }}

%\vskip -0.6cm
\vskip -0.35cm
%\thicklines
{%$\kern -1.4pt$ 
$\kern +35.4pt$\framebox(10,10){}}
\ \ \ \ \ \ \ \ \ \
\ $\kern -2.0pt$ 
  $\kern -2.95pt$
{\framebox(10,10){}}$\kern -2.95pt$
{\framebox(5,10){{\tiny +}}}
$\kern 12.5pt$ 
\ \ \ \ \ \ \ {\framebox(10,10){}}$\kern -2.95pt$
{\framebox(5,10){{\tiny +}}} 
\ \ \ \ \ \ \ \ \ \ \ \ $\kern 1pt${\framebox(10,10){}} 

%\vskip -0.6cm
\vskip -0.1cm
%\thicklines
{%$\kern -1.4pt$ 
$\kern +35.40pt$\framebox(-0.02,10){}}
\ \ \ \ \ \ \ \ \ \ \ \ 
\ $\kern -0.0pt$ 
  $\kern -2.935pt$
{\framebox(10,10){+}}
$\kern 10.7pt$  
\ \ \ \ \ \ \ \ \ {\framebox(10,10){}}$\kern -2.95pt$
{\framebox(10,10){}} 
\ \ \ \ \ \ \ \ \ \ \ $\kern -4.1pt$ {\framebox(10,10){}} 
\end{thicklines}
\vskip 0.2cm

The output coincides with (\ref{21to3}) if the duality 
(\ref{dualityc}) is employed.
\medskip

The following two examples deal with more involved $\,w$.
Let  us take $\,w^{(1)}=s_1s_2\cdots s_p$\, and
$\,w^{(2)}=s_p\cdots s_2s_1$\, for
$k=-\ell$ and the corresponding weights
$c^{(1)}=\ell\om_{p+1},\, c^{(2)}=\ell p\,\om_{p+1}$. 
Then one finds that 
$$
b^{(1)}=\ell (p+1)\om_1 \hbox{\, and\, }
b^{(2)}=\ell (p+1)\, \om_p.
$$
We apply here the general relation
$\rho-w(\rho)=\sum_{\al\in R_+\,\cap \,w(R_-)}\al$ and 
formulas $\al_i=2\om_i-\om_{i-1}-\om_{i+1}$,
which hold in any $A_n$ if $n>i$ ($\om_0=0$).
See below the combinatorial interpretation of this
calculation.

Thus we obtain that 
\begin{align}\label{s1-sp}
\h_{\rr,\ss}(c^{(i)};\,q,t\!=\!q^{-\ell},a) = q^{\bullet}\,
\h_{\rr,\ss}(b^{(i)};\,q,t\!=\!q^{-\ell},a) \hbox{\, for\, } i=1,2.
\end{align}

When $\ell=1$, the Young diagrams for 
$\{\om_{p+1},(p+1)\om_1\}$  and those for 
$\{p\om_{p+1},(p+1)\om_p\}$ are transposed to each other. 
Therefore (\ref{s1-sp}) follow from (\ref{dualityc}) 
for $\ell=1$; the relation $tq=1$ is obviously preserved by 
the duality transformation $t\leftrightarrow q^{-1}$. 

Combining (\ref{s1-sp}) for $\ell=1$
and $i=1$ with its restriction to
$A_{n}$ for $n=p+1$, where $\om_{n}$ and $\om_1$ result in
coinciding $\tilde{J\!D}$\~polynomials,
\begin{align*}
\h_{\rr,\ss}(\om_1;\,q,t=q^{-1},a=-t^{n+1}) = q^{\bullet}\,
\h_{\rr,\ss}(n\om_1;\,q,t=q^{-1},a=-t^{n+1}).
\end{align*}

Here we see some links to the so-called 
{\em colored differentials\,}; 
see e.g., \cite{DGR,GS,FGS,GGS}. However our  
color exchange from (\ref{jones-bca})
preserves the number of boxes of Young diagrams,
in contrast to the colored differentials from \cite{GGS}
and other papers. Also, the differentials there do not seem 
mathematically rigorous (the differentials
from \cite{KhR1,KhR2,Ras} are) and it is surprising to us 
that the colored differentials are introduced
only for symmetric and wedge powers.

\comment{
Combinatorially,
applying the differential associated with $t^l q^m=1$
(or a similar ones) is essentially replacing the monomials
$t^{i+l}q^{j+m}$ by $t^aq^b$ (in some manner)
in the expression for $\h$. 
Compare with Section 3.6 from \cite{CJ}.
There will be further reductions of topological 
nature if one wants to establish a connection with 
unstable homology. 
However, adding weights (colors) to the Khovanov-Rozansky theory 
is generally a difficult problem and we do not see any
mathematically rigorous rationale for the existence of
colored differentials.  
We emphasize that the color exchange among different
$i\om_1$ (or different $\om_i$) 
can hold only for DAHA-Jones polynomials, not for
superpolynomials (in contrast to \cite{GGS}). 
The coincidence formulas (\ref{jones-bca}) are 
only among Young diagrams with the same number of boxes. 
}

%the theory of {\em perfect modules\,} of DAHA
%from \cite{C101}, Section 3.10.2. 

\medskip
{\sf Fat hooks.}
Generalizing, let us consider
a {\em fat hook\,} corresponding to
$c=v_1\om_{u_1+u_2}+v_2\om_{u_1}$ for $u_1,u_2,v_1,v_2\in \Z_+$.
Let $k=-\ell\in -\Z_+$ (can be fractional) subject to the following 
relations:
\begin{align}\label{fathook}
v_1+v_2\,\ge\, u_2\ell\,\in\, \Z_+\, \ni\, (u_1+u_2)\ell\,\ge\, v_2.
\end{align}
We set $b=(v_1+v_2-u_2\ell)\om_{u_1+u_2}+
((u_1+u_2)\ell-v_2)\om_{u_2}$ and claim that
\begin{align}\label{symqt}
\h_{\rr,\ss}(c;\,q,t\!=\!q^{-\ell},a) = q^{\bullet}\,
\h_{\rr,\ss}(b;\,q,t\!=\!q^{-\ell},a) \hbox{\, for any\, }
\rr,\ss. 
\end{align}

Combinatorially, $\,u_1,u_2,v_1,v_2\,$ are 
the number of lines and columns in the Young diagram
for $c$. This diagram and the corresponding
Young block-diagram for $b\,$ are as follows: 

\vskip 0.2cm
\begin{thicklines} 
\ \ \ \ \,\,  \ \ \ \ \ \ \ \ \ \ 
$v_1$  $v_2$ \ \ \ \ \ \ \ \ \ \ \  \ \ \ \ \ \ \ \ \ \ \ 
$v_1\!+\!v_2\!-\!u_2\ell$ \ \ \ $(u_1\!+\!u_2)\ell\!-\!v_2$ 

\ \ \raisebox{-1ex}{\hbox{\ $\la_c=$\ }}\ \ 
$u_1$\ {\framebox(10,10){}}$\kern -2.95pt$
{\framebox(10,10){}}\ \raisebox{-1ex}{,} 
\ \ \ \ \raisebox{-1ex}{\hbox{\ $\la_b=$\ }}\ \ $u_2$  
{\framebox(80,10){}}$\kern -2.95pt$
{\framebox(80,10){}} \ \raisebox{-1ex}{.} 

\vskip -0.27cm
\ \  \ \ \ \ \ \ \ \ \ \ $u_2$ $\kern -4.1pt$ 
{\framebox(10,10){}}\ 
\ \ \ \ \ \ \ \  \ \ \ \ \ \ \ \ \ \ $\kern -0.8pt$ 
$u_1$ $\kern -4.2pt$ {\framebox(80,10){}}
\end{thicklines}
\vskip 0.2cm

The permutation $w$ transposes
the two block-lines of the diagram $\la_c$, i.e. it
moves the first $u_1$ lines of $\la_c$ down and the last 
$u_2$ lines up (without changing their relative positions 
within the corresponding block-lines).
The integrality condition for $k=\ell$ is\,
gcd$(u_1,u_2)\ell\in \Z_+$. For instance, $\ell=1/u$ is
possible if $u_1=u=u_2$; in this case, the inequalities from
(\ref{fathook}) become  $\,v_1+v_2\ge 1,\, v_2\le 2\,$ and
the pair $v_1,v_2$ will be transformed to \,$v_1+v_2-1, 2-v_2$.

Note that when $\ell=1$ (i.e. $k=-1$) and $u_1=v_2,u_2=v_1$,
the diagram $\la_b$ is the transposition of $\la_c$. Then
(\ref{symqt}) follows from the duality. 
\smallskip

{\sf Rectangles.}
This construction has the following application to 
the Young diagrams that
are rectangles. Setting
$v_1=0=v_1+v_2-u_2\ell=v_2-u_2\ell$, one has
$c=v_2\om_{u_1}=u_2\ell\om_{u_1}$ and
$b=u_1\ell\om_{u_2}$. Thus the rectangles
\,$\la_c=\{u_2\ell\times u_1\}$ and
\,$\la_b=\{u_1\ell\times u_2\}\,$ satisfy 
(\ref{symqt}). Here $u_1$ and $u_2$ are
arbitrary positive, \,gcd$(u_1,u_2)\ell\,$ must be
integral. Note that one can set here
$v_2=0=v_1+v_2-u_2\ell=v_1-u_2\ell$.
The corresponding weights become \,$c=u_2\ell\om_{u_1+u_2}$ and 
$b=(u_1+u_2)\ell\om_{u_2}$; thus we arrive at a particular case
of the previous relation (use $u_1+u_2> u_2$). Finally:

%Thus we obtain the duality we obtain that for any $u,v\in \Z_+$,
\begin{align}\label{symrect}
&q^{\bullet}\,\h_{\rr,\ss}(u\om_v;\,q,q^{-u},a) =
q^{\bullet}\,\h_{\rr,\ss}(uv\om_1;\,q,q^{-u},a),\\ 
%=q^{\bullet}\,\h_{\rr,\ss}(v\om_u;\,q^u,q^{-1},a)  
%=q^{\bullet}\,\h_{\rr,\ss}(\om_{uv};\,q^u,q^{-1},a)
&t^{\bullet}\,\h_{\rr,\ss}(u\om_v;\,t^{-v},t,a) =
t^{\bullet}\,\h_{\rr,\ss}(\om_{uv};\,t^{-v},t,a),\notag 
\end{align}
where the second formula is due to the duality (\ref{dualityc}). 
Actually it can be justified directly  (similarly to the proof 
of the first formula). 
\smallskip

For instance, $t^8\h_{3,2}(2\om_2;\,t^{-2},t,a)\ =$
$\h_{3,2}(\om_4;\,t^{-2},t,a)$ for 

\centerline{$\h_{3,2}(2\om_2;\,q,t,a)\ =$}
\vskip 0.1cm

\renewcommand{\baselinestretch}{0.5} 
{\small

\noindent
\(
1+\frac{a^4 q^{10}}{t^2}+q^2 t+q^3 t+q^2 t^2+q^3 t^2+q^4 t^2
+q^4 t^3+q^5 t^3+q^4 t^4+q^5 t^4+q^6 t^4+q^6 t^5+q^7 t^5
+q^6 t^6+q^7 t^6+q^8 t^8+a^3 \bigl(q^9+q^{10}+\frac{q^7}{t^2}
+\frac{q^8}{t^2}+\frac{q^7}{t}+\frac{q^8}{t}+q^9 t
+q^{10} t\bigr)
\)

\noindent
\(
+a^2 \bigl(q^5+q^6+2 q^7+q^8+\frac{q^5}{t^2}
+\frac{q^4}{t}+\frac{2 q^5}{t}+\frac{q^6}{t}+q^6 t+3 q^7 t
+2 q^8 t+q^7 t^2+q^8 t^2+q^9 t^2+q^8 t^3+2 q^9 t^3
+q^{10} t^3+q^9 t^4\bigr)
\)

\noindent
\(
+a \bigl(q^2+q^3+q^4+q^5
+\frac{q^2}{t}+\frac{q^3}{t}+2 q^4 t+3 q^5 t+q^6 t+q^4 t^2
+2 q^5 t^2+2 q^6 t^2+q^7 t^2+2 q^6 t^3+3 q^7 t^3+q^8 t^3
+q^6 t^4+2 q^7 t^4+q^8 t^4+q^8 t^5+q^9 t^5+q^8 t^6
+q^9 t^6\bigr)\,,
\)

}
\renewcommand{\baselinestretch}{1.2} 
\smallskip

\centerline{$\h_{3,2}(\om_4;\,q,t,a)\ =$}
\vskip 0.1cm

\renewcommand{\baselinestretch}{0.5} 
{\small

\noindent
\(
1 + \frac{a^4 q^4}{t^6} + q t + q t^2 + q t^3 + q t^4 + q^2 t^4 + 
q^2 t^5 + 2 q^2 t^6 + q^2 t^7 + q^2 t^8 + q^3 t^9 + q^3 t^{10} + 
q^3 t^{11} + q^3 t^{12} +  q^4 t^{16} + 
a^3 \bigl(q^4 + \frac{q^3}{t^6} + \frac{q^3}{t^5} + \frac{q^3}{t^4} + 
\frac{q^3}{t^3} + 
\frac{q^4}{t^2} + 
   \frac{q^4}{t} + q^4 t\bigr)
\)

\noindent
\(
+a^2 \bigl(3 q^3 + \frac{q^2}{t^5} + \frac{q^2}{t^4} + 
\frac{2 q^2}{t^3} + \frac{q^2}{t^2} + 
   \frac{q^3}{t^2} + \frac{q^2}{t} + \frac{2 q^3}{t} + 3 q^3 t + 
2 q^3 t^2 + q^3 t^3 + 
q^4 t^3 + q^4 t^4 + 2 q^4 t^5 + q^4 t^6 + q^4 t^7\bigr)
\)

\noindent
\(
+a \bigl(q + 2 q^2 + \frac{q}{t^3} + \frac{q}{t^2} + \frac{q}{t} + 
\frac{q^2}{t} + 3 q^2 t + 
3 q^2 t^2 + 
   2 q^2 t^3 + q^3 t^3 + q^2 t^4 + 2 q^3 t^4 + 3 q^3 t^5 + 
3 q^3 t^6 + 
   2 q^3 t^7 + q^3 t^8 + q^4 t^9 + q^4 t^{10}+ 
q^4 t^{11} + q^4 {t^{12}}\bigr)\,.
\)
}
\renewcommand{\baselinestretch}{1.2} 
\smallskip

{\sf Color exchange as a recovery tool.}
The positivity conjecture for rectangles (the last 
from \cite{CJ} in type $A$ that remains open) states
that the coefficients of $\h_{\rr,\ss}(u\om_v;\,q,t,a)$ 
are all positive. In spite of the color-exchange 
connection, the number of monomials
in terms of $q,t,a$ (ignoring their numerical coefficients)
is different for  $\h_{3,2}(2\om_2;\,q,t,a)$ and 
$\h_{3,2}(\om_4;\,q,t,a)$ ($66$ and $60$). Both
result in $35$ different monomials after the 
substitution $q=t^{-2}$ (a significant reduction).

We claim that the positivity of $\h_{3,2}(2\om_2;\,q,t,a)$,
its self-duality, knowing the above specialization at 
$q=t^{-2}$ and the evaluation at $q=1$ uniquely determine this
polynomial (up to a multiplier $q^\bullet t^\bullet$). 

\smallskip
Let us provide some details.
Assuming that we know the $q,t$\~monomials that occur
in $\h_{3,2}(2\om_2;\,q,t,a)$, the duality restricts 
the $66$-dimensional space  
of the corresponding undetermined coefficients 
to the space $V$ of dimension $36$. Decomposing such $V$ 
with respect to $a^i$, one has $V=\oplus_{i=0}^4 V_i$, where
the corresponding dimensions 
are $\{\hbox{dim}_0\!=\!9, 12, 10, 4,\hbox{dim}_4\!=\!1\}$. 
Then the dimensions of the kernels of the specialization
map $q=t^{-2}$ are $\{1,3,2,0,0\}$.  The evaluation 
$q=1$ has the corresponding kernels of  
dimensions $\{2,4,3,1,0\}$. There will be no 
common kernel if these $2$ maps are combined.
%\smallskip

Here we do not use the positivity, but assume that the
monomials are those from the actual $\h_{3,2}(2\om_2;\,q,t,a)$.
%Without imposing the positivity, the recovery of the 
%coefficients can be unique only
%modulo the ideal $I=\bigl((1-q)(1-t)(t-q^{-2})(q-t^{-2})\bigr)$. 
The positivity gives that any other set of coefficients
will be obtained from the actual one when at least 
$8$ monomials (for one of $V_i$) are replaced by some 
other $8$ (with positive coefficients) and the
corresponding difference belongs to the
ideal $\bigl((1-q)(1-t)(t-q^{-2})(q-t^{-2})\bigr)$. 
This appeared impossible.

Moreover, one can involve here the reduction $t\!=\!q$
(which gives essentially the HOMFLYPT polynomial of $T(3,2)$
for the weight $2\om_2$ and can be assumed known),
as well as  
$\h_{3,2}(2\om_2;\,q,t,a\!=\!-t^2)\!=\!1$. The resulting system
of (liner) equations for the coefficients of 
$\h_{3,2}(2\om_2;\,q,t,a)$ will then become significantly  
overdetermined.   

\smallskip
However the ``next" polynomial $\h_{3,2}(2\om_3;\,q,t,a)$ 
can not be recovered from the color-exchange relations  
(\ref{symrect}) combined with the evaluations 
at $t=1$ and $q=1$; here $t=q^{-2}$ for the reduction
to $6\om_1$ and $q=t^{-3}$ for that to
$\om_6$. As above, we do not use the positivity,
but assume that the monomials can be only those
actually present in  $\h_{3,2}(2\om_3;\,q,t,a)$
(with undetermined coefficients).

\comment{
The total space $V=\oplus_{i=0}^6 V_i$
of possible coefficients is now of dimension $284$; 
the corresponding dimensions are $\{47,67,69,54,34,12,1\}$. 
The dimensions of the total kernel
for the $4$ restriction maps above are $\{0,4,7,6,1,0,0\}$. 
Omitting $\,t=1=q\,$, such dimensions jump to
$\{10,24,27,18,7,0,0\}$. Adding $\,t=q\,$ (the reduction to
the HOMFLYPT polynomial) makes the recovery of coefficients 
unique in this case; the kernel becomes zero.

Note that the number of monomials
in $u\om_v$ upon the reduction to $uv\om_1$
and $\om_{uv}$ seems approximately $(1/u)$ and $(1/v)$ of
the total number of monomials (exactly
$124$,  $106$ for $2\om_3$)
and one can expect the dimension of the common
kernel of these $2$ reduction to be approximately 
$(uv-u-v)/uv$ of the total number of monomials if these 
reductions are ``in a general position". However this is not the 
case for $2\om_3$; the dimension of such a kernel is $87$, which is 
significantly greater than $286/6\simeq 48$ (or $286-124-106=56$
if we use the exact dimensions of the images). Obviously the 
unique recovery cannot be expected using such and similar 
specializations for general rectangles, even assuming (as we did) 
that we know the $q-t$\~monomials involved. 
}

\setcounter{equation}{0}
\section{\sc DAHA of type 
\texorpdfstring{{\mathversion{normal}$C^\vee C_1$}}{C-check-C-1}}

\subsection{\bf Main definitions}
Double affine Hecke algebra of type $C^\vee C_1$, denoted
by $\sH=\sH_{q,u,v}$ in this paper, is generated by 
$U_1,U_0,V_1,V_0$ subject 
to the relations
\begin{align}\label{CCdoub}
&(U_i\!-\!u_i^{1/2})(U_i\!+\!u_i^{-1/2})=0,\ \,
(V_i\!-\!v_i^{1/2})(V_i\!+\!v_i^{-1/2})=0,\notag\\ 
&q^{1/4}V_1V_0U_0U_1= 1, \ \,i = 0,1
\hbox{\, here and below;\, we set :}\\
&\X \!\equal\! V_1^{-1}\,U_1^{-1}\!=\!
q^{1/4}V_0 U_0,\ 
\Y \!\equal\! U_0\, U_1\!=\!q^{-1/4}V_0^{-1}V_1^{-1}.\notag
\end{align}
The natural definition ring is 
$\Z_{q,u,v}\equal \Z[q^{\pm 1/4},u_i^{\pm 1/2}, v_i^{\pm 1/2}]$,
though we will mainly use its field of fractions $\Q_{q,u,v}$
when the Askey-Wilson polynomials are needed. Note that
$V_1,\X,\Y$ obviously generate $\sH_{q,u,v}$.
Here and below we closely follow \cite{NS}, modifying
the generators and parameters from 2.22 there as follows:
\begin{align}\label{NoStokman} 
&V_i\!\mapsto\! U_i^{-1},\, V_i^\vee\!\mapsto\! V_i^{-1} (i=0,1),\, 
q_{\NS}\!\mapsto\! q^{-\frac{1}{2}},\, 
k_i\!\mapsto\! u_i^{-\frac{1}{2}},\, 
u_i\!\mapsto\! v_i^{-\frac{1}{2}},\\ 
%&T_i\mapsto \hat{U}_i^{-1},\ T_i^\vee\mapsto \hat{V}_i^{-1},
&X\equal V_1^{-1}(V_1^\vee)^{-1}\!=\!
q_{\NS}^{1/2}\,V_0V_0^\vee  \mapsto\  
\X^{-1}\!=\!U_1V_1\!=\!q^{-1/4}\,U_0^{-1}V_0^{-1},\notag\\
&Y\equal V_1V_0\!=\!q_{\NS}^{-1/2}(V_1^\vee)^{-1}(V_0^\vee)^{-1}
\mapsto \Y^{-1}\!=\!U_1^{-1}U_0^{-1}\!=\!q^{1/4}V_1V_0.\notag
\end{align}
These arrows mean that $V_i$ from \cite{NS}
is our $U_i^{-1}$, their $q_{\NS}$ is our $q^{-1/2}$, 
their $X$ is our $\X$ and so on.  The
relation $q_{\NS}^{1/2}V_1^\vee V_1 V_0 V_0^\vee=1$
becomes that from (\ref{CCdoub}) in our notations.

Such changes are
convenient to establish the relations to DAHA of type
$A_1$ and to the superpolynomials for $T(2\pp+1,2)$ 
(a surprising application of the $C^\vee C_1$-theory).
A convenient way of obtaining our relations and formulas
from their ones (or those from \cite{Ob})
is simply by transposing all terms ($AB\cdots C\mapsto C\cdots BA$)
followed by the substitutions $V\mapsto U, V^\vee\mapsto V, q_{\NS}
\mapsto q^{1/2}, X\mapsto \X, Y\mapsto \Y$.
\smallskip

{\sf Automorphisms.}
We redefine the automorphisms $\si, \tau, \eta$
from \cite{Ob}, Proposition 1.3 as follows:
$\si\mapsto \si^{-1}$, $\tau\mapsto \tau_+^{-1};$
$\eta$ remains unchanged (this is 
an extension of the Kazhdan-Lusztig involution).    

%Also the {\em anti-isomorphism\, $\nu$ 
%from 8.5 in \cite{Nos} will be denoted by $\psi$. 

In full detail, the following maps can
be uniquely extended to automorphisms of the
whole $\sH_{q,u,v}$:
\begin{align} \label{tau+cc}
\tau_+:\ 
&U_0\mapsto V_0=q^{-1/4}\X\, U_0^{-1},\  
V_0\mapsto V_0 U_0^{-1} V_0^{-1},\  
U_1\mapsto U_1,\notag\\
&V_1 \mapsto V_1,\ \X\mapsto \X,\  
\Y\mapsto\, V_0\,U_1 =
q^{-1/4}\,V_1^{-1}\Y^{-1} U_1,\ u_0\leftrightarrow v_0,
\\ \label{tau-cc}
\tau_-: \ 
&V_1\mapsto V_0=q^{-1/4}V_1^{-1}\Y^{-1},\,
V_0\,\mapsto \,V_0^{-1} V_1 V_0,\,  
U_1\mapsto\, U_1,\notag\\ 
&U_0 \mapsto U_0,\ \Y \mapsto \Y, \
\X \mapsto V_0^{-1}U_1^{-1}\!=\!
q^{1/4} U_0\X^{-1} U_1^{-1},\ v_0\leftrightarrow v_1,
\\ \label{etacc}
\eta: \ 
&U_1\!\mapsto U_1^{-1},\, U_0\!\mapsto U_0^{-1},\,
V_1\!\mapsto U_1 V_1^{-1} U_1^{-1},\, 
V_0\!\mapsto U_0^{-1} V_0^{-1} U_0,\notag\\
&\X \mapsto \X^{-1}, \
\Y \mapsto U_1 \Y^{-1} U_1^{-1},
\, u_i\mapsto u_i^{-1}, q\mapsto q^{-1}, 
v_{i}\mapsto v_{i}^{-1}.
\end{align}
The parameter $q$ and respectively
$u_1,v_1$ and $u_0,u_1$ remain unchanged under the action
of $\tau_{\pm}$. We set $\si\equal
\tau_+\tau_-^{-1}\tau_+= \tau_-^{-1}\tau_+\tau_-^{-1}$;
this automorphism transposes $t_0$ and $v_1$ and 
corresponds to $\mu^{-1}$ from 
\cite{NS}, 8.5. One has:
\begin{align} \label{sicc}
\si:\ &U_0\mapsto V_1,\ \, V_0\mapsto V_1V_0 V_1^{-1},\ \,
V_1\mapsto U_1^{-1} U_0 U_1,\ \, U_1\mapsto U_1,\\
&\X\!\mapsto \Y^{-1},\  \Y\!\mapsto U_1^{-1} \X^{-1} U_1,
\ u_0\leftrightarrow v_1, u_1\!\mapsto u_1, v_0\!\mapsto v_0.\notag
\end{align}

In the notations from \cite{NS}, formulas 
(\ref{tau+cc})-(\ref{etacc}) are as follows:
\begin{align*}
\tau_+:\ 
&V_0\mapsto V_0^\vee\,=\,q_{\NS}^{-1/2}V_0^{-1}X,\ \, 
V_0^\vee\mapsto (V_0^\vee)^{-1}V_0 V_0^\vee,\  
V_1\mapsto V_1,\\
&V_1^\vee \mapsto V_1^\vee,\ X\mapsto X,\  
Y\mapsto V_1V_0^\vee\!=\!
q_{\NS}^{-1/2}V_1Y^{-1}(V_1^\vee)^{-1},\ k_0\leftrightarrow u_0,
\\
\tau_-: \ 
&V_1^\vee\mapsto V_0^\vee=q_{\NS}^{-1/2}V_0^{-1}X,\ 
V_0^\vee\mapsto V_0^\vee V_1^\vee (V_0^\vee)^{-1},\  
V_1\mapsto V_1,\\ 
&V_0 \mapsto V_0,\ Y \mapsto Y, \
X \mapsto V_1^{-1}(V_0^\vee)^{-1}\!=\!
q_{\NS}^{1/2} V_1^{-1}X^{-1} V_0,\ u_0\leftrightarrow u_1,
\\
\eta: \ 
&V_1\!\mapsto V_1^{-1},\, V_0\!\mapsto V_0^{-1},\,
V_1^\vee\!\mapsto V_1^{-1}(V_1^\vee)^{-1} V_1,\, 
V_0^\vee\!\mapsto V_0 (V_0^\vee)^{-1} V_0^{-1},\\
&X \mapsto X^{-1}, \ \,
Y \mapsto V_1^{-1} Y^{-1} V_1,\,
\ k_i,u_i,q\mapsto k_i^{-1},u_{i}^{-1},q^{-1}
\,(i=1,2).
\end{align*}

Next, we will need the 
{\em anti-involution\,} $\nu=\vph_{\NS}$ 
from 8.5 in \cite{NS}, which corresponds to
our $\vph$:
\begin{align} \label{vphcc}
\vph\,:\ 
&U_0\leftrightarrow V_1,\ \,
U_1\mapsto U_1,\ \, V_0 \mapsto V_0,\ \  
\X\leftrightarrow \Y^{-1},\ \,u_0 \leftrightarrow v_1,
\\ 
\vph_{\NS}:\ 
&V_0\leftrightarrow V_1^\vee,\ 
V_1\mapsto V_1,\  V_0^\vee \mapsto V_0^\vee,\   
X\leftrightarrow Y^{-1},\ k_0 \leftrightarrow u_1.
\notag
\end{align}
All relations of the reduced theory in Section \ref{sect:Aut}
holds for $\tau_{\pm},\si,\vph$. For instance,  
\begin{align}\label{vphsirel}
&\tau_\mp=
\vph\tau_\pm\vph=\si\tau_\pm^{-1}\si^{-1},\ 
\eta\tau_{\pm}\eta=\tau_{\pm}^{-1},\ 
\vph\si\vph=\si^{-1}=\eta\si\eta, \notag\\ 
& 
\vph\eta\vph=\eta\si^{-2}=
\si^2\eta,\ \, \si^2(H)=U_1^{-1}(H)U_1 \for H\in \sH_{q,t,u}.
\end{align}

Also, due to the group nature of the definition
of $\sH_{q,u,v}$, we have the {\em inversion 
anti-involution\,}
$\star$\,, sending all generators  of
$\sH_{q,u,v}$ and parameters $q,u_i,v_i$ (and
their products) to the corresponding inversions.
\smallskip

Finally, let us define the {\em sign-automorphisms\,} 
of $\sH_{q,t,u}$:
\begin{align}\label{chngsign}
&\varsigma_{x}: V_0\mapsto -V_0, v_0^{1/2}\mapsto -v_0^{1/2},
V_1\mapsto -V_1, v_1^{1/2}\mapsto -v_1^{1/2},\\
&\varsigma_{y}: U_0\mapsto -U_0, u_0^{1/2}\mapsto -u_0^{1/2},
V_0\mapsto -V_0, v_0^{1/2}\mapsto -v_0^{1/2},\notag\\
&\varsigma_{\ast}: U_i\mapsto -U_i, u_i^{1/2}\mapsto -u_i^{1/2},
V_i\mapsto -V_i, v_i^{1/2}\mapsto -v_i^{1/2},\notag\\
&\varsigma_{q}: U_1\mapsto -U_1, u_1^{1/2}\mapsto -u_1^{1/2}, 
q^{1/4}\mapsto -q^{1/4},\notag
\end{align}
where the remaining generators and parameters are unchanged
under the action of the corresponding $\varsigma$.

The group $GL_{\,2}^\wedge(\Z)$ fixes $\varsigma_q,\varsigma_\ast$
and acts in the group $\F_{\!2}^2$ generated by 
$\varsigma_x,\varsigma_y$ through its projection 
$\varpi$ onto $GL_2(\F_{\!2})$. One has:
\begin{align}\label{varsixy}
&\varsigma_x: \X, \Y\mapsto -\X, \Y, \
\varsigma_y: \X, \Y\mapsto \X, -\Y,\\ 
&\varsigma_q: \X,\Y\mapsto -\X,-\Y,\  
\varsigma_\ast: \X,\Y\mapsto \X,\Y.\notag
\end{align} 

\subsection{\bf Polynomial representation}
The following presentation of the algebra
$\sH_{q,u,v}$ by
Demazure-Lusztig operators was found
by Noumi for $C^\vee C_n$; see also \cite{Sa}, Section 2.3.
Using the notations from \cite{NS},
\begin{align*}%\label{UVops}
&\hat{V}_i=k_i s_i+\frac{(k_i-k_i^{-1})+(u_i-u_i^{-1})X_i}
{1-X_i^2}(1-s_i) \for i=0,1, \\
&X_1=X,\ X_0=q^{\frac{1}{2}}/X, \
s_1(X^m)=X^{-m},\ s_0(X^m)=q^{m} X^{-m},\\
&\hat{V}_i^{-1}=k_i^{-1} s_i+
\frac{(k_i-k_i^{-1})X_i^2+(u_i-u_i^{-1})X_i}
{1-X_i^2}\,(1-s_i),\\
&\hat{V}^\vee_1(f)=\hat{V}_i^{-1}(f)/X,\
\hat{V}^\vee_0(f)=\hat{V}_0^{-1}(f/X_0)=q^{-\frac{1}{2}}
\hat{V}_0^{-1}(Xf),\\
&(\hat{V}^\vee_1)^{-1}(f)=\hat{V}_1(Xf),\
(\hat{V}^\vee_0)^{-1}(f)\!=X_0\hat{V}_0(f)\!=\!
q^{\frac{1}{2}}X^{-1}\hat{V}_0(f).
\end{align*}
In our notations,
\begin{align}\label{Uopers}
&\hat{U}_i\,=\,u_i^{1/2} s_i+\frac{(u_i^{1/2}-u_i^{-1/2})+
(v_i^{1/2}-v_i^{-1/2})\X_i}
{1-\X_i^2}(1-s_i), \\
&\X_1=\X,\ \X_0=q^{\frac{1}{4}}/\X, \
s_1(\X^m)=\X^{-m},\ s_0(\X^m)=\!q^{\frac{m}{2}} \X^{-m}\!,\notag\\
&\hat{U}_i^{-1}\!=\!u_i^{-1/2} s_i\!+\!
\frac{(u_i^{1/2}\!-\!u_i^{-1/2})\X_i^2\!+\!
(v_i^{1/2}\!-\!v_i^{-1/2})\X_i}
{1-\X_i^2}(1\!-\!s_i),\notag\\
\label{Vopers}
&\hat{V}_1(f)=\hat{U}_1^{-1}(f/\X),\
\hat{V}_0(f)=\hat{U}_0^{-1}(f)/\X_0=q^{-\frac{1}{4}}
\X\hat{U}_0^{-1}(f),\\
&\hat{V}_1^{-1}(f)\,=\,\X\hat{U}_1(f),\ \,
\hat{V}_0^{-1}(f)\,=\,\hat{U}_0(\X_0\,f)\,=
\, q^{\frac{1}{4}}\,\hat{U}_0(f/\X).\notag
\end{align}

Here $f\in \v=\v_{{}_{\X}}\equal\Z_{q,u,v}[\X^{\pm 1}]$, which we
will call the {\em polynomial representation\, };
it is supplied with the action 
$\sH_{q,u,v}\ni H\mapsto \hat{H}$ defined
via (\ref{Uopers}),(\ref{Vopers}).
Thus (following Noumi) we claim that the relations
from (\ref{CCdoub}) are satisfied for 
$\{\hat{U}_i,\hat{V}_i\}$. The operator $\X$ acts
as the multiplication by $\X$.
Indeed, (\ref{Vopers}) gives that $\hat{\X}(f)=
 \hat{V}_1^{-1}\hat{U}_1^{-1}(f)=
\X\hat{U}_1\bigl(\hat{U}_1(f)\bigr)=f$.

The existence of this representations provides the
{\em PBW Theorem\,}, which states that the following
decomposition is unique: 
\begin{align}\label{ccpbw}
H\!=\!\!\sum_{n,\ep,m}C_{n,\ep,m}  \X^n V_1^\ep \Y^m 
\hbox{\, for any\,} 
H\in \sH_{q,u,v}, \  n,m\in\Z, \ep\!=\!0,1. 
\end{align}
Using this theorem (which can be proved directly),
$\v$ is the $\sH_{q,u,v}$\~module
induced from the one-dimensional {\em
evaluation character\,} $\chi$ on the subalgebra
$\mathscr{U}$ generated by $U_0,U_1$:
\begin{align}\label{ccpolind}
\v\!=\!\hbox{Ind}_{\mathscr{U}}^\sH(\chi),\ 
\chi:\, U_1\!\mapsto\! u_1^{1/2}, U_0\!\mapsto\! u_0^{1/2}, 
\Y\!=\!U_0U_1\!\mapsto\! (u_0 u_1)^{1/2}.
\end{align}
This readily gives that $\tau_-$ sends $\v$ to its image
under $v_0\leftrightarrow v_1$.
\smallskip

The {\em difference Dunkl operator\,} for $C^\vee C_1$
is $\hat{\Y}=\hat{U}_0\,\hat{U}_1$. 
The {\em nonsymmetric Askey-Wilson\,} 
polynomials $\e_n\,(n\in \Z)$ are defined from 
the relations
\begin{align}\label{Eccheck}
&\hat{\Y}(\e_n)\, =\, (u_0u_1)^{-\hbox{\tiny sgn}(n)/2}
q^{-n/2}\,\e_n \for
n\in \Z, \hbox{\, where : }\\
&\e_n\!=\!\X^n+C^n_{-n}\X^{-n}\!+\!\!\sum_{|m|<|n|}\! C^n_m\X^m, \,
C^n_{-n}=0 \hbox{\, for\,\ } n>0.\notag
\end{align}
Here $\hbox{\small sgn}(n\le 0)=-1$ and $+1$ otherwise,
i.e. $0$ is treated as negative. Note the formula
\begin{align}\label{tauminpol}
\tau_-(\e_n)\,=\,q^{-n^2/4}
(u_1u_0)^{-|n|/2}\,
\e_n\mid_{v_0\leftrightarrow v_1} \for n\in \Z.
\end{align}

The {\em symmetric Askey-Wilson polynomials\, }
$\p_n\, (n\in \Z_+)$ from \cite{AW} can be defined
from the relation
\begin{align}\label{Pccheck}
&(\hat{\Y}\!+\!\hat{\Y}^{-1})(\p_n) = 
\bigl((u_0u_1)^{\frac{1}{2}}q^{\frac{n}{2}}+
(u_0u_1)^{-\frac{1}{2}}q^{-\frac{n}{2}}\bigr) \p_n,\\
&\hbox{where\, } \p_n\,=\,\X^n+\X^{-n}+ \!\sum_{0\le m <n}\! C_m
(\X^m+\X^{-m}).\notag
\end{align}

Also, $\p_n$  are the $u_1$\~symmetrizations 
of $\e_n$, i.e. they can be defined using relation 
(\ref{nontosym}) for $A_1$ with $\,T_1,t_1\,$ there replaced by
$U_1,u_1$. They are formally real; namely, $\p_n^\star=\p_n$, where
$\star\,$ sends $\X\mapsto \X^{-1}$ and all parameters to
their reciprocals.
\smallskip

One has: 
\begin{align}\label{cc-e1}
\e_0=1,\ \e_1\!=\!\X+\frac{(q^{1/2} u_0/v_0)^{1/2}(1\!-\!v_0)\!+\!
u_0 (q u_1/v_1)^{1/2}(1\!-\!v_1)}{1-q^{1/2}u_0u_1},
\end{align}

\comment{
\begin{align*}
&(-u_1*v_0 - q^2*u_0^2*u_1*v_0 - q*u_0*v_1 - q^3*u_0*u_1^2*v_1 
+ q*u_0*v_0^2*v_1\\ 
&+ q^3*u_0*u_1^2*v_0^2*v_1 + u_1*v_0*v_1^2 + 
q^2*u_0^2*u_1*v_0*v_1^2)/ \\
&((-1 + q^2*u_0*u_1)*(1 + q^2*u_0*u_1)*v_0*v_1)\\
&+ X^{-1} + ((-1 + u_1^2 - q^2*u_1^2 + q^2*u_0^2*u_1^2)*X)\\
&/((-1 + q^2*u_0*u_1)*(1 + q^2*u_0*u_1)).
\end{align*}
\begin{align*}
\frac{ u_1^{\frac{1}{2}} v_0^{\frac{1}{2}} + 
q u_0 u_1^{\frac{1}{2}} v_0^{\frac{1}{2}} + 
q^{\frac{1}{2}} u_0^{\frac{1}{2}} v_1^{\frac{1}{2}}
+ q^{\frac{3}{2}} u_0^{\frac{1}{2}} u_1 v_1^{\frac{1}{2}} } 
{ v_0^{\frac{1}{2}} v_1^{\frac{1}{2}}
(1 - q u_0 u_1)}\\
 -\frac{ q^{\frac{1}{2}} u_0^{\frac{1}{2}} v_0 v_1^{\frac{1}{2}} + 
q^{\frac{3}{2}} u_0^{\frac{1}{2}} u_1 v_0 v_1^{\frac{1}{2}} 
+ u_1^{\frac{1}{2}} v_0^{\frac{1}{2}} v_1
+ q u_0 u_1^{\frac{1}{2}} v_0^{\frac{1}{2}} v_1 }
{ v_0^{\frac{1}{2}} v_1^{\frac{1}{2}}
(1 - q u_0 u_1) } \\
+ \X^{-1} + 
 \frac{1 - u_1 + q u_1 - q u_0 u_1}
{1 - q u_0 u_1}\X.
\end{align*}
}

\vskip -0.4cm
\begin{align}\label{cc-e-1}
\e_{-1}=\X^{-1}+ \frac{(1 - u_1) + q^{1/2} u_1(1 - u_0)}
{1 - q^{1/2} u_0 u_1}\,\X \hskip 2.cm&\\
+\,\frac{ (q^{1/2} u_0/v_0)^{\frac{1}{2}}
(1  + q^{1/2} u_1)(1 - v_0)+ (u_1/v_1)^{\frac{1}{2}}
(1 + q^{1/2} u_0) (1-v_1)}
{ 1 - q^{1/2} u_0 u_1 }&.\notag
\end{align}

The corresponding symmetric Askey-Wilson polynomial is
\begin{align}\label{cc-p-1}
&\p_{1}=\p_1=(1+u_1^{1/2}V_1)(\e_1)\ =\ 
\X +\X^{-1}\\
+\, 
&\frac{ (q^{1/2} u_0/v_0)^{\frac{1}{2}}
(1\!+\! u_1)(1\! -\! v_0) +  (u_1/v_1)^{\frac{1}{2}}
(1 \!+\! q^{1/2} u_0) (1\!-\!v_1)}
{ 1 - q^{1/2} u_0 u_1 }.\notag
\end{align}

Note that $\e_1$ is not $X$ as it is for $A_1$ and
$\e_n$ can have the terms $X^m$ in odd and even degrees.
However, the following symmetries for the sign-automorphisms
from (\ref{chngsign}) and (\ref{varsixy}) hold:
\begin{align}\label{varsen}
&\varsigma_x(\e_n)\!=\!(-1)^n\e_n\!=\!\varsigma_q(\e_n),\ 
\varsigma_y(\e_n)\!=\!\e_n\!=\!\varsigma_q(\e_n) 
\hbox{\, for any\, } n\in\Z.
\end{align} 
We treat here $\e_n(\X)$ as elements of $\sH_{q,u,v}$. 
These relations control the appearance of square roots of the
parameters in $\e_n$; the coefficients of $E_n$ are in terms 
of $q,t$ in the case of $A_1$, i.e. free of square roots. 
\medskip

The  $\e$\~polynomials can be obtained using
the intertwining operators for $\sH_{q,u,v}$ 
from \cite{Sa,NS}:
\begin{align}\label{Sinter}
&\S_1\equal U_1^{-1}\Y^{-1}-\Y^{-1}U_1^{-1},\ 
\S_0\equal\Y^{-1} U_0^{-1}-U_0^{-1}\Y^{-1},\\
&\hbox{intertwining\, }\Y\,:
\S_1\Y\ =\ \Y^{-1}\S_1,\ \ q^{-1}\S_0\Y\ =\ \Y^{-1}\S_0.\notag
\end{align}
They result in the following recurrence relations for 
$m=0,1,2\ldots$\, :
\comment{
\begin{align}\label{Sintere}
\frac{1-q^{-\frac{2m+1}{2}}u_0^{-1}u_1^{-1}}
{q^{-\frac{m+1}{2}}u_1^{-1}u_0^{-1/2}}
\e_{m+1}= \hat{\S}_0(\e_{-m}),\ 
\frac{1-q^{-m-1}u_0^{-1}u_1^{-1}}{q^{-\frac{m+1}{2}}u_0^{-1/2}}
\e_{-m-1}=\hat{\S}_1(\e_{m+1}).
\end{align}
}
\begin{align}\label{Sintere}
\hat{\S}_0(\e_{-m})\!=\!
\frac{q^{\frac{2m\!+\!1}{2}}u_0u_1\!-\!1}{q^{\frac{m}{2}}u_0^{1/2}}
\e_{m\!+\!1},\ \,
\hat{\S}_1(\e_{m\!+\!1})\!=\!
\frac{q^{\frac{2m\!+\!2}{2}}u_0u_1\!-\!1}
{q^{\frac{m\!+\!1}{2}}u_0^{1/2}u_1}
\e_{-m\!-\!1}.
\end{align}

See \cite{NS} for explicit formulas for the
$\e$\~polynomials, their norms and evaluations.
We will need below the evaluation formulas
(recalculated to our notations):
\comment{
\begin{align*}
&\e_{m+1}(u_1^{-1/2}v_1^{-1/2})\!=\!
\frac{u_1^{(m-1)/2}(v_0v_1)^{1/2\!+m/2}
(1\!+\!q^{-1}u_1^{-1})\cdots(1\!+\!q^{-m}u_1^{-1})}
{(1-q^{-m-1}u_0^{-1}u_1^{-1})\cdots
(1-q^{-2m-1}u_0^{-1}u_1^{-1})}\\
&\times\prod_{i=0}^{m}(v_0^{-1/2}\!+\!
q^{-i-1/2}u_0^{-1/2}u_1^{-1/2}v_1^{-1/2})
(1\!-\!q^{-i-1/2}u_0^{-1/2}u_1^{-1/2}v_0^{-1/2}v_1^{-1/2}),
\end{align*}
}
\begin{align}\label{cceval-1}
\e_{m+1}&\bigl((u_1v_1)^{-\frac{1}{2}}\bigr)=
\frac{(u_1v_0v_1)^{-(m+1)/2}
(1+q^{\frac{1}{2}} u_1)\cdots(1+q^{\frac{m}{2}} u_1)}
{(1-q^{\frac{m+1}{2}}u_0 u_1)\cdots
(1-q^{\frac{2m+1}{2}}u_0 u_1 )}\notag\\
&\times\,\prod_{i=0}^{m}\,(v_0^{1/2}+
q^{\frac{2i+1}{4}}u_0^{1/2}u_1^{1/2}v_1^{1/2})\,
(1-q^{\frac{2i+1}{4}}u_0^{1/2}u_1^{1/2}v_0^{1/2}v_1^{1/2}),
\end{align}
\begin{align}\label{cceval-2}
\e_{-m-1}&\bigl((u_1v_1)^{-\frac{1}{2}}\bigr)=
\frac{(u_1v_0v_1)^{-(m+1)/2}
(1+q^{\frac{1}{2}} u_1)\cdots(1+q^{\frac{m+1}{2}}u_1)}
{(1-q^{\frac{m+2}{2}}u_0 u_1)\cdots
(1-q^{\frac{2m+2}{2}}u_0 u_1 )}\notag\\
&\times\,\prod_{i=0}^{m}\,(v_0^{1/2}+
q^{\frac{2i+1}{4}}u_0^{1/2}u_1^{1/2}v_1^{1/2})\,
(1-q^{\frac{2i+1}{4}}u_0^{1/2}u_1^{1/2}v_0^{1/2}v_1^{1/2}).
\end{align}
Also, $\p_{m+1}\bigl((u_1v_1)^{\pm\frac{1}{2}}\bigr)
=(1+u_1)\e_{m+1}\bigl((u_1v_1)^{-\frac{1}{2}}\bigr)$\, for\, 
$m\in \Z_+$.

\comment{
\begin{align*}
\e_{-m-1}(u_1^{-1/2}v1^{-1/2})=
\frac{
(1+q^{-1}u_1^{-1})\cdots(1+q^{-m-1}u_1^{-1})
(v_0^{-1/2}+q^{-1/2}u_0^{-1/2}u_1^{-1/2}v_1^{-1/2})
(v_0^{-1/2}+q^{-3/2}u_0^{-1/2}u_1^{-1/2}v_1^{-1/2})\cdots
(v_0^{-1/2}+q^{-(2m+1)/2}u_0^{-1/2}u_1^{-1/2}v_1^{-1/2})
(1-q^{-1/2}u_0^{-1/2}u_1^{-1/2}v_0^{-1/2}v_1^{-1/2})
(1-q^{-1/2}u_0^{-3/2}u_1^{-1/2}v_0^{-1/2}v_1^{-1/2})\cdots
(1-q^{-(2m+1)/2}u_0^{-1/2}u_1^{-1/2}v_0^{-1/2}v_1^{-1/2})}
{u_1^{-m/2}
(t_0t_1)^{-m/2}
(1-q^{-m-2}u_0^{-1/2}u_1^{-1/2})\cdots
(1-q^{-2m-2}u_0^{-1/2}u_1^{-1/2})}
\end{align*}
}
 
Concerning the orthogonality relations and other
formulas for the $\e$\~polynomials, we refer to \cite{NS},
where $X,q_{\NS},k_0,k_1,u_0,u_1$ must be replaced
by $\X^{-1}, q^{-1/2}, u_0^{-1/2},u_1^{-1/2},$
$v_0^{-1/2},v_1^{-1/2}$ to match our setting.
See \cite{Ko,Mac,Sa,Sto} for general theory.

\subsection{\bf Relations to 
\texorpdfstring{{\mathversion{bold}$A_1$}}{A1}}
Let us first reduce Definition \ref{double}
to the case of DAHA of type  $A_1$, denoted  
by $\HH^{A_1}_{q,t}$. For $A_1$, let $\al=\al_1$, $s=s_1$ 
and $\om=\om_1$ be the fundamental weight; 
then $\alpha=2\omega$ and $\rho=\om$. The extended affine 
Weyl group $\widehat{W}=<\!s,\om\!>\,$ in the $A_1$\~case
is a free group generated by the involutions $\,s\,$ and
$\,\pi\equal \om s$.

The generators of
$\HH^{A_1}_{q,t}$ are
$$
Y=Y_{\om_1}=\pi T,\ \,T=T_1,\ \,X=X_{\om_1}
$$
subject to the quadratic relation $\,(T-t^{1/2})(T+t^{-1/2})=0\,$
and the cross-relations:
\begin{align}\label{dahaone}
&TXT=X^{-1},\ T^{-1}YT^{-1}=Y^{-1},\ Y^{-1}X^{-1}YXT^2q^{1/2}=1.
\end{align}
Using $\pi=YT^{-1}$, the second relation becomes $\pi^2=1$.
This algebra is defined over
$$
\Z_{q;t}\equal\Z[q^{\pm 1/2},t^{\pm 1/2}].
$$

The following maps can be extended to automorphisms
of $\HH^{A_1}_{q,t}$:
\begin{align}\label{tau-taum-def}
&\tau_+(X)\!=\!X,\ \tau_+(T)\!=\!T,\ \tau_+(Y)\!=\!q^{-1/4}XY,\
\tau_+(\pi)\!=\!q^{-1/4}X\pi,\\
\notag
&\tau_-(Y)=Y,\ \ \tau_-(T)=T,\ \ \tau_-(X)=q^{1/4}YX,\ \
\tau_-(\pi)\,=\,\pi,\\
&\si(X)\,=\, Y^{-1},\ \ \si(T)\,=\,T,\ \ \si(Y)\,=\,XT^2,\ \ 
\si(\pi)\,=\,XT,\notag\\
&\eta(X)\!=\! X^{-1},\,
\eta(Y)\!=\!T Y^{-1} T^{-1},\, \eta(T)\!=\!T^{-1},\ 
q,t \mapsto q^{-1},t^{-1}.
\label{etaa-one}
\end{align}
Recall that 
$\si= \tau_+\tau_-^{-1}\tau_+\ =\ \tau_-^{-1}\tau_+\tau_-^{-1}$
and we add $q^{\pm 1/4}$ to the ring of definition
of $\HH^{A_1}_{q,t}$.

We also have two anti-involutions:
\begin{align}\label{eta-a-one}
&\vph(X)=Y^{-1},\ \vph(Y)=X^{-1},\ \vph(T)=T,\ 
q,t \mapsto q,t,\\ 
&X^\star=X^{-1},\  Y^\star=Y^{-1},\, T^\star=T^{-1},\ \, 
q,t\mapsto q^{-1},t^{-1}.\notag
\end{align}
\smallskip

\begin{theorem}\label{CCTOONE}
(i) Let us send $u_1^{1/2}\!\!\mapsto\! t^{1/2}$ and make 
$u_0^{1/2}\!\!=\!1\!=\!v_0^{1/2}\!\!=\!v_1^{1/2}$. Then the map
\begin{align}\label{cc-a-map}
U_1\mapsto T,\ &U_0\mapsto \pi=YT^{-1}, V_0\mapsto 
\tilde{\pi}\equal q^{1/4}YT^{-1}X^{-1}, V_1\mapsto
\overline{\pi}\equal XT,\notag\\ 
&\X\!=V_1^{-1}U_1^{-1}\mapsto (XT)T^{-1}=X, \
\Y\!=U_0 U_1\mapsto (YT^{-1})T=Y 
\end{align}
can be extended to a homomorphism of algebras
$\xi:\sH_{q,u,v}\to \HH^{A_1}_{q,t}$ (with the same $q$
in both), compatible with $\tau_{\pm},\eta,\vph$
and the inversion $\star$. The nonsymmetric
Askey-Wilson polynomials $\e_n(\X)$ 
become $E_n(X)$ for $A_1$ for $n\in \Z$, where
the evaluation formulas from
(\ref{cceval-1}),(\ref{cceval-2}) become those in
(\ref{macde-eval}) for $A_1$.

(ii) Let us make now $u_1^{1/2}=t^{1/2}=u_0^{1/2}$
and consider $\sH_{q^2,u,v}$ (i.e. replace 
$\,q^{1/4}$ by $\,q^{1/2}\,$ in its definition). Then 
\begin{align}\label{cc-a-new}
&U_1\mapsto T,\ U_0\mapsto T_0\equal Y^2T^{-1}, 
V_0\mapsto \tilde{\pi}, V_1\mapsto
\breve{\pi}\equal q^{1/4}XTY,\\ 
\X\!= &V_1^{-1}U_1^{-1}\!\mapsto\! 
q^{-1/4}Y^{-1}X\!=\!\tau_-^{-1}(X),\ 
\Y\!=\!U_0 U_1\!\mapsto\! (Y^2 T^{-1})T=Y^2\notag
\end{align}
give a homomorphism $\ze:\sH_{q^2,u,v}\to \HH^{A_1}_{q,t}$,
which induces the following group homomorphism from 
the group $\Ga^\wedge_{\,0}(2)\equal <\tau_+^2,\tau_-,\eta>\subset
GL^\wedge_{\,2}(\Z)$ acting in $\sH_{q^2,u,v}$ to
the group $GL^\wedge_{\,2}(\Z)$ acting in  $\HH^{A_1}_{q,t}$:
\begin{align}\label{cc-new-aut}
&\ze:\ \tau_+^2\mapsto \tau_-^{-1}\tau_+\tau_-,\ 
\tau_-\mapsto \tau_-^2,\ \eta\mapsto \tau_-^{-1}\eta\tau_-. 
\end{align}

(iii) The  map $\ze\,$ from (\ref{cc-a-new}) is compatible
with the following map of the polynomial 
representations: 
\begin{align}\label{cctoapol}
\sH_{q^2,u,v}\circlearrowright \v_{{}_{\!\X}}\ni F(\X)
\,\overset{\ze}\longmapsto\,
\tau_-^{-1}\bigl(F(X)\bigr)(1)\in 
\v_{{}_{\!X}}\circlearrowleft \HH^{A_1}_{q,t}.
\end{align}
For instance, the 
polynomials $\e_n(\X)$ upon the substitution
$$
q^{1/4}\mapsto q^{1/2},\ \, 
u_0^{1/2},u_1^{1/2}\mapsto t^{1/2},\ \,
v_0^{1/2},v_1^{1/2}\mapsto 1
$$ 
coincide with $E_n(X)$ (of type $A_1$) for any $n$ 
and, accordingly, the evaluation of $\e_n$ at $(u_1 v_1)^{-1/2}$
becomes that for $E_n$ at $t^{-1/2}$.
\end{theorem}

{\it Proof.} It is known that the substitution 
$u_1=1=v_0=v_1$ makes $\sH_{q,u,v}$ the DAHA
of type $A_1$ (claimed in Part $(i)$). However we need some
details here. We rewrite the identity
\begin{align*}
(TY^{-1})(q^{1/4}YT^{-1}X^{-1})(XT)= q^{1/4}T =
\pi\tilde{\pi}\overline{\pi}
\end{align*}
in $\HH^{A_1}_{q,t}$ 
as follows:\, $q^{1/4} \overline{\pi}\tilde{\pi}\pi T=1$.
Then we use that $\pi, \tilde{\pi}, \overline{\pi}$ are 
involutions, and this readily gives the map $\xi$
from $(i)$.
The comparison of the automorphisms involved is
straightforward (and is known); DAHA of type
$C^\vee C_1$ was actually introduced as a 
generalization of that for $A_1$. The compatibility
of $\xi$ with the definition of the polynomial
representation (here and below)
results from its interpretation as 
an induced module from (\ref{ccpolind}).
\smallskip

Let us focus on $(ii,iii)$. Now the identity 
in $\HH^{A_1}_{q,t}$ we need is
\begin{align*}
T_0T=Y^2=q^{-1/2}\tilde{\pi}\breve{\pi}=
q^{-1/2}(q^{1/4}YT^{-1}X^{-1})(q^{1/4}XTY),
\end{align*}
which can be rewritten as $q^{1/2}\breve{\pi}\tilde{\pi}T_0T=1$
and readily the $\ze$\~images of
the generators of $\sH_{q^2,u,v}$ from (\ref{cc-a-new}).
This proves $(iii)$.

Let us calculate the $\ze$\~image of $\tau_+^2$.
Note that $\ze(\tau_+)$ is not well defined since $\tau_+$
transposes $u_0$ and $v_0$. One has :
\begin{align*}
&\tau_+(U_1)\!=\!U_1,\, \tau_+(\X)\!=\!\X,\, 
\tau_+^2(\Y)\!=\!\tau_+(q^{-1/4}V_1\Y^{-1}U_1)\!=\!
V_1^{-1} U_1^{-1}\Y V_1 U_1,\\
&\ze(\tau_+)\,:\ 
T\ \mapsto\ T,\ \ \tau_-^{-1}(X)\ \mapsto\ \tau_-^{-1}(X),
\ \ 
Y^2\ \mapsto\ T^{-1}X^{-2}T\\
&=\ze(V_1^{-1} U_1^{-1}\Y V_1 U_1)
=\breve{\pi}T^{-1}Y^2\breve{\pi}T=
(\breve{\pi}Y^{-1})(Y T^{-1}Y) (Y\breve{\pi})\\
&=\!(\breve{\pi}Y^{-1})T (Y\breve{\pi})T
=(q^{\frac{1}{4}}XT)T (q^{-\frac{1}{4}}
T^{-1}X^{-1})T=T^{-1}X^{-2}T.
\end{align*}
On the other hand,
\begin{align*}
&\tau_-^{-1}\tau_+\tau_-(Y^2)=
\tau_-^{-1}\tau_+\tau_-^{-1}\bigl(\tau_-^2(Y^2)\bigr)=\si(Y^2)=
T^{-1}X^{-2}T,\\
&\tau_-^{-1}\tau_+\tau_-\bigl(\tau_-^{-1}(X)\bigr)=
\tau_-^{-1}(X),\ \tau_-^{-1}\tau_+\tau_-(T)=T.
\end{align*}
This justifies that $\ze(\tau_+)=\tau_-^{-1}\tau_+\tau_-.$
\smallskip

The $\ze$\~image of $\tau_-$ is as follows:
\begin{align*}
&\tau_-(U_1) = U_1,\,\ \tau_-(\Y) = \Y,\,\ 
\ze(\tau_-(\X)) = \ze(V_0^{-1}U_0^{-1})=\tilde{\pi}T Y^{-2}
\\
&\ze(\tau_-):\ 
T\ \mapsto\ T,\ \ Y^2\ \mapsto\ Y^2,
\ \ 
\ze(X)=\tau_-^{-1}X\ \mapsto\ \ze\bigl(\tau_-(X)\bigr)\\
&= \ze(V_0^{-1} U_1^{-1})
\!=\!\tilde{\pi}T^{-1}\!=\!q^{1/4}YT^{-1}X^{-1}T^{-1}\!=
\!q^{1/4}YX\!=\!
\tau_-(X).
\end{align*}
We conclude that $\ze(\tau_-)=\tau_-^2$.
\smallskip

Similarly, the relation $\ze(\eta)=\tau_-^{-1}\eta\tau_-$ 
results from: 
\begin{align*}
&\ze(\eta(\X))\,=\,\ze(X^{-1})=\tau_-^{-1}(X^{-1})\, =\,
\tau_-^{-1}\eta\tau_-\bigl(\ze(X)\bigr),\\
&\ze(\eta(\Y^2))\!=\! \ze(U_1 \Y^{-2} U_1^{-1})
\!=\!\eta(Y^2),\ 
\ze(\eta(U_1))\!=\!\eta(T).
\end{align*}
The fact that $\ze$ is the group automorphism of
$\Ga^\wedge_{\,0}(2)$ follows from the calculation above
(and irreducibility of $\v$). This can be readily
checked directly, since any relation between $\tau_+^2$
and $\tau_-$ holds for $\tau_+$ and $\tau_-^2$.
The conjugation by $\begin{pmatrix}1 & 0\\ 0 & 2\end{pmatrix}$
proves this fact in the corresponding matrices,
which is of course connected with our formula $\ze(\Y)=Y^2$.

\subsection{\bf DAHA-Jones polynomials for 
\texorpdfstring{{\mathversion{bold}$C^\vee C_1$}}{C-check-C-one}}
We define the evaluation functional for $H\in \sH_{q,u,v}$
as $\{H\}_{ev}\equal \hat{H}(1)((u_1v_1)^{-1/2})$ for the action
of $H$ in $\v$. Equivalently,
one can use the PBW-decomposition from (\ref{ccpbw}) and 
set 
\begin{align}\label{ccevf}
\bigl\{\,\sum_{n,\ep,m}C_{n,\ep,m}\X^n\, V_1^\ep\, 
\Y^m\,\bigr\}_{ev} =
\!\sum_{n,\ep,m}C_{n,\ep,m} (u_1v_1)^{-n}\, 
u_1^\ep\, (u_1u_0)^{m}. 
\end{align}
By construction, $\{\, \vph(H) \,\}_{ev}=\{H\}_{ev}$; recall that
$\vph$ transposes $u_0$ and $v_1$ fixing $u_1$.
Also,  $\{\, \eta(H) \,\}_{ev}=\{H\}^\star$, where $\star$
is the inversion of all parameters (and the standard generators), 
and $\{\, \si^2(H) \,\}_{ev}=\{H\}$. Thus the connections between
$\{\cdot\}_{ev}$ from (\ref{evsym}) all hold for $C^\vee C_1$.
However the DAHA-Jones polynomials must now carefully
address the action of $GL_{\,2}^\wedge(\Z)$ on the parameters
$\{u_i,v_i\}$. 
\smallskip

Given a pair \,$(\rr,\ss)$ satisfying the condition
\,gcd$(\rr,\ss)=1$,  we associate 
with it $\,\ga=\ga_{\rr,\ss}\in 
\Ga_{\,0}^\wedge(2)\equal 
\bigl\{\,
\begin{pmatrix}\al & \be\\ \ga & \de\\ \end{pmatrix}\,\mid\,
\be\in 2\Z \,
\bigr\}$ such that its first column $\,(\al,\ga)^{tr}$
is $\,(\rr,\ss)^{tr}$. Here $\ga\,$  is unique modulo right 
multiplication by 
matrices from $\begin{pmatrix} 1 & 2\Z\\ 0 & 1\\ \end{pmatrix}$. 
We lift such $\ga$ to $\hat{\ga}\in GL^\wedge_{\,2}(\Z)$.

Recall that $\hat{\ga}$ acts on the parameters $u_1,u_0,v_1,v_0$;
for the group $PSL^\wedge_{\,2}(\Z)$ generated by $\tau_{\pm}$
this action is via
$\varpi(\ga)$. One has:
\begin{align}\label{hgauv}
&\tau_{+}(u_1,\!\underline{u_0},\!v_1,\!\underline{v_0})\!=\! 
(u_1,\!v_0,\!v_1,\!u_0),\ 
\tau_{-}(u_1,\!u_0,\!\underline{v_1,\!v_0})\!=\!
(u_1,\!u_0,\!v_0,\!v_1),\\
\eta(u_1,& u_0,\!v_1,\!v_0) =\, 
(u_1^{-1}\!,\!u_0^{-1}\!,\!v_1^{-1}\!,\!v_0^{-1}),\ \, 
\vph(u_1,\underline{u_0,\!v_1},\!v_0) =\, 
(u_1,\!v_1,\!u_0,\!v_0),\notag
\end{align}
and $\eta(q)=q^{-1}$ (the other three preserve $q$). These
transformations are naturally extended to the fractional
powers of the parameters. 
\smallskip

Let $\varpi(\ss)=0,1$ be the parity of $\ss$; we also 
denote $\ga$ modulo $(2)$ by $\varpi(\ga)$, which is
$\begin{pmatrix}1 &  0\\ \varpi(\ss) & 1\\ \end{pmatrix}$
for $\ga=\ga_{\rr,\ss}$. We also set $\de(\ga)=\det(\ga)=\pm1$.

We will sometimes use the notation 
$\bar{\Ga}_0(2)$ and 
$\bar{\Ga}_{\,0}^\wedge(2)$ for the intersection of
$\Ga_0(2)$ and 
$\Ga_{\,0}^\wedge(2)$ with  $PSL_2(\Z)$, considered
as a subgroup of $GL_2(\Z)$, and 
$PSL^\wedge_{\,2}(\Z)$. Any admissible $(\rr,\ss)$ can be lifted
to $\ga\in \bar{\Ga}_0(2)$ and $\hat{\ga}\in 
\bar{\Ga}_{\,0}^\wedge(2)$; 
then only $\tau_{\pm}$
will occur in $\hat{\ga}$ ($\eta$ will not be involved).  
\smallskip

The next theorem and its proof
follow Theorem \ref{MAINTHM}.
% we use $\,\thicksim\,$ for the coincidence of Laurent polynomials
%up to a factor that is a
%monomial in terms of $q^{1/4},\{u_i^{1/2},v_i^{1/2}\}$.

\begin{theorem} \label{CCTHM} 
(i) Given an admissible pair $(\rr,\ss)$, i.e. satisfying
the condition \,gcd$(\rr,\ss)=1$, and $m\in \N$, we set 
\begin{align}\label{jones-dc}
\j\!D_{\rr,\ss}(m\,;\,q,u,v)\equal 
\bigl\{\,\hat{\ga}_{\rr,\ss}(\p_m)\bigr\}_{ev}\,/\,
\eta^{\de(\ga_{\rr,\ss})}
\tau_-^{\varpi(\ss)}\bigr(\bigl\{\p_m\bigr\}_{ev}\bigl),
\end{align}
where $\hat{\ga}_{\rr,\ss}$ is a lift of $\ga_{\rr,\ss}
\in \Ga_0(2)$ with
the first column $(\rr,\ss)^{tr}$.
Then 
$\j\!D_{\rr,\ss}(m\,;\,q,u,v)$ does not depend on the 
particular choice of $\ga_{\rr,\ss}\in \Ga_0(2)$
and $\hat{\ga}_{\rr,\ss}\in GL_{\,2}^{\wedge}(\Z)$ 
and is a Laurent polynomial in terms of $q^{1/4},u_i^{1/2},
v_i^{1/2}$ for $i=0,1$. 
One can switch to the $\e$\~polynomials here:
\begin{align}\label{jones-nonsymc} 
&\j\!D_{\rr,\ss}(m\,;\,q,u,v)\!=\!  
\bigl\{\,\hat{\ga}(\e_{\pm m})\bigr\}_{ev}\,/\,
\eta^{\de(\ga)}\tau_-^{\varpi(\ss)}
\bigl(\bigl\{\e_{\pm m}\bigr\}_{ev}\bigr),\
 \ga\!=\!\ga_{\rr,\ss}.
\end{align}

(ii) One has: $\j\!D_{1,\ss}(m\,;\,q,u,v)\!=\!
(q^{\frac{m}{2}}u_1u_0)^{-\ss m/2}$, 
\begin{align}\label{jones-symc} 
%\thicksim 1,\ \
&\j\!D_{\rr,-\ss}(m\,;\,q,u,v)=\bigr(\j\!D_{\rr,\ss}
(m\,;\,q,u,v)\bigl)^\star \hbox{\, for any \,}\rr,\ss,
\notag\\
&\j\!D_{\rr,\ss}(m\,;\,q,u,v)\,=\ 
\j\!D_{\ss,\rr}(m\,;\,q,u,v)\hbox{\,\, if \,\,}
\rr\ss \hbox{\,\, is even\,.}
\end{align}
Moreover, for any\,  
$\rr,\ss, v_1,v_0$\, and \,$m,n\in \Z_+$, 
\begin{align}\label{bcw-relc}
q^{\,(n+m)/2\,} u_0u_1=1 
\ \,\Rightarrow\ \,
\j\!D_{\rr,\ss}(m\,;\,q,t)=
\j\!D_{\rr,\ss}(n\,;\,q,t).
\end{align}
If $q=1$, then $\,\e_{-m}\!=\!\p_{m}$,
$\,\p_{m+n}\!=\!\p_m \p_n\,$ for any $\,m,n\in \Z_+\,$ and
\begin{align}\label{jones-evalc}
\j\!D_{\rr,\ss}(m\,;\,q=1,u,v)=%\thicksim
\j\!D_{\rr,\ss}(1\,;\,q=1,u,v)^m
\hbox{\, for any \,} \rr,\ss.
\end{align}

(iii) Let us use the sign-automorphisms 
$\,\varsigma_x,\varsigma_y,\varsigma_\ast,\varsigma_q\,$
of $\sH_{q,t,u}$ from (\ref{chngsign});\ \,
$GL_{\,2}^\wedge(\Z)$ naturally
acts in the group generated by these $\varsigma$. 
Namely, it fixes $\,\varsigma_\ast, \varsigma_q\,$
and acts on $\varsigma_x,\varsigma_y$ via formulas
(\ref{hgauv}). Then relations (\ref{varsixy}) and
(\ref{varsen}) give that
\begin{align*}
\varsigma\{ \ga_{\rr,\ss}(\e_m) \}_{ev}/
\{ \ga_{\rr,\ss}(\e_m) \}_{ev} =
\varsigma(\e_m)/\e_m \hbox{\, for any such\, }\varsigma,\ m\in \Z.
\end{align*}
Since the normalization point in 
(\ref{jones-dc}) depends on $\varpi(\ss)$,  
the square roots of $\,q^{1/2},u_i,v_i\,$ can occur in
$\j\!D_{\rr,\ss}(m\,;\,q,u,v)$ only as follows:
\begin{align}\label{varsQ}
&\j\!D_{\rr,\ss}(m\,;\,q,u,v)\ =
\a + 
(q^{\frac{1}{2}}u_1 u_0 v_1 v_0)^{\frac{1}{2}}\b
\hbox{\, for}\\ 
&\a,\b\,\, \in\,\, (u_0u_1)^{\frac{m\, \rr\, \ss}{2}}\,
\Z\,[\,q^{\pm\frac{1}{2}}, u_1^{\pm1},
u_0^{\pm1},v_1^{\pm1},v_0^{\pm1}\,]. \notag
\end{align}
\end{theorem}
\vskip -1.25cm 
\sq
\vskip 0.2cm

Obviously, the broken symmetry $(\rr,\ss)\mapsto (\ss,\rr)$
and the fact that 
$\j\!D_{2\pp+1,1}$ are nontrivial seem non-topological. 
However exactly the 
polynomials $\j\!D_{2\pp+1,1}$ appear connected with   
the {\, superpolynomials\, } of knots 
$T(2\pp+1,2)$, so they do have topological-geometric meaning.
\medskip

\subsection{\bf Reductions of parameters}
We will begin with the following case, which is a 
counterpart of the case $t=1$ for $A_1$ (where all
DAHA-Jones polynomials are tilde-trivial).

\begin{proposition}\label{U0U1ONE}
Let us impose the relations
$v_0^{1/2}\!=\!1\!=\!v_1^{1/2}$ and also
$u_1^{1/2}\!=\!q^{1/2}\!=\!u_0^{-1/2}$. 
We assume that $\rr$ is odd using  the symmetry (\ref{jones-symc})
if necessary; let $\,\rr=2\pp+1$. Then
\begin{align}\label{juouone}
\j\!D_{\rr,\ss}(m\,;\,q,u_1\!=\!q,u_0\!=\!q^{-1},1,1) =
q^{-m(\rr\, \ss\, m/2+\,\pp)}\, \frac{1+q^{m\,\rr}}{1+q^m}
\end{align}
for any admissible $\rr,\ss, m$.\sq
\end{proposition}

Here we check that $\j\!D_{\rr,\ss}(m\,;\,q,u,v)$ depends
only on $\,\pp\,$ up to a monomial factor
and then calculate $\,\j\!D\,$
for the pair $(2\pp+1,1)$.
\smallskip
 
Theorem \ref{CCTOONE} results in the 
following reduction formulas.

\begin{theorem}\label{JCCTOONE}
(i) Let $u_1^{1/2}\mapsto t^{1/2}$ and 
$u_0^{1/2}=1=v_0^{1/2}=v_1^{1/2}$. Then 
\begin{align}\label{jcc-aone}
\j\!D_{\rr,\ss}(m\,;\,q,u,v) =
J\!D^{A_1}_{\rr,\ss}(m\,;\,q,t)
\equal J\!D^{A_1}_{\rr,\ss}(m\om_1\,;\,q,t).
\end{align}

(ii) We take now $\,u_1^{1/2}=t^{1/2}=u_0^{1/2}$,
and replace $q^{1/4}$ by $\,q^{1/2}$. We also
assume that $\,\rr$ is odd, which is possible
due to (\ref{jones-symc}). Then 
\begin{align}\label{jcc-atwo}
\j\!D_{\rr,\ss}(m\,;\,q^2,t,\!t,1,\!1) =
J\!D^{A_1}_{\rr,2\ss}(m\,;\,q,t)
\end{align}
for any admissible $\rr, \ss, m$. \sq
\end{theorem}

The reduction formulas here
do not require the tilde-normalization
(up to a monomial factor) of $J\!D^{A_1}$.
This normalization is actually
unnecessary in the $A_1$\~case due
to the formulas from \cite{CJ}, Section 2.5\,:
$$
\tilde{J\!D}\raisebox{0.2ex}{\hbox{${}^{A_1}$}}_
{\kern -12pt \rr,\ss}\,(m\,;\,q,t)=
q^{m^2 \rr\ss/4}t^{m(\rr+\ss-1)/2}
J\!D^{A_1}_{\rr,\ss}(m\,;\,q,t).
$$
Accordingly, we have the following leading term 
in the case of $C^\vee C_1$: 
\begin{align}\label{cc-tilde}
\j\!D_{\rr,\ss}(m\,;\,q^2,t,\!t,1,\!1)\!=\!
q^{-m^2\frac{\rr\ss}{4}}t^{-m\frac{\rr+\ss-1}{2}}
\!\sum_{i,j\!\ge 0} C_{ij}q^i t^j \hbox{\, for\, } C_{00}=1.
\end{align}

Let us provide the simplest example of the
colored  $\j\!D$\~polynomial
and its reductions to $A_1$ (apart from the pairs
$(2\pp+1,1)$ to be discussed next).
Setting $\,A=(1-v_0),\ B=(1-u_0),\ \breve{q}=q^{-2}\,$,
\vskip 0.1cm
 
\centerline{$\j\!D_{3,2}(2\,;\,q,u,v)\,=\,
\j\!D_{2,3}(2\,;\,q,u,v)\,=$}

\vskip 0.1cm

\renewcommand{\baselinestretch}{0.5} 
{\small

\comment{
\noindent
\begin{align*}
&\frac{A B q^{15}}{{u_0}^{3/2} {u_1}^{5/2} 
\sqrt{{v_0}} {v_1}^{3/2}}+
\frac{B^2 q^{10}}{{u_0}^2 {u_1}^2 {v_1}}-
\frac{A^2 q^8}{{u_1}^2 {v_0} {v_1}}\\
&+\frac{1}{{u_0}^{3/2} {u_1}^{7/2} \sqrt{{v_0}} 
{v_1}^{3/2}}A q^3 \bigl(q^{14}+q^{16}+q^{10} {u_1}-
q^8 {u_0} {u_1}-2 q^{10} {u_0} {u_1}\\
&-q^4 {u_0} {u_1}^2-2 q^6 {u_0} {u_1}^2-q^8 {u_0} {u_1}^2+
q^4 {u_0}^2 {u_1}^2+q^6 {u_0}^2 {u_1}^2+{u_0}^2 {u_1}^3\\
&+q^2 {u_0}^2 {u_1}^3-q^9 \sqrt{{u_0}} \sqrt{{u_1}} 
\sqrt{{v_0}} \sqrt{{v_1}}+q^6 {u_1} {v_1}+q^8 {u_1} {v_1}-
q^2 {u_0} {u_1}^2 {v_1}\\
&-q^4 {u_0} {u_1}^2 {v_1}\bigr)
+\frac{1}{{u_0}^2 {u_1}^3 {v_1}^2}B q^2 \bigl(q^{16}+q^{18}+
q^{12} {u_1}-q^{12} {u_0} {u_1}\\
&-q^8 {u_0} {u_1}^2+q^{12} {v_1}+
q^{10} {u_1} {v_1}+q^4 {u_0} {u_1} {v_1}-q^4 {u_0} {u_1}^2 {v_1}-
{u_0}^2 {u_1}^2 {v_1}\bigr)\\
&+\frac{1}{{u_0}^2 {u_1}^4 {v_0} {v_1}^2}q^2 \bigl(q^{22} {v_0}-
q^{10} {u_0} {u_1}^2 {v_0}-q^{14} {u_0} {u_1}^2 {v_0}-
q^6 {u_0} {u_1}^3 {v_0}\\
&+q^6 {u_0}^2 {u_1}^3 {v_0}+
q^2 {u_0}^2 {u_1}^4 {v_0}+q^{10} {u_0} {u_1} {v_1}-
q^4 {u_0}^2 {u_1}^2 {v_1}+{u_0}^3 {u_1}^3 {v_1}\\
&+q^{14} {u_1} {v_0} {v_1}+2 q^4 {u_0}^2 {u_1}^2 {v_0} {v_1}-
q^6 {u_0} {u_1}^3 {v_0} {v_1}-2 {u_0}^3 {u_1}^3 {v_0} {v_1}\\
&-q^2 {u_0}^3 {u_1}^3 {v_0} {v_1}-q^4 {u_0}^2 {u_1}^2 {v_0}^2 {v_1}+
{u_0}^3 {u_1}^3 {v_0}^2 {v_1}+q^6 {u_1}^2 {v_0} {v_1}^2\bigr)
\end{align*}
\begin{align*}
&\frac{q^2}{{u_0}^2 {u_1}^4 {v_0} {v_1}^2}
\Bigl(\,A B\, q^{13}u_1(u_0 u_1 v_0 v_1)^{1/2}+
B^2\, q^{8}{u_1}^2 {v_0} {v_1}-
A^2\, q^{6}{u_0}^2 {u_1}^2 {v_1}\\
&+A\, q ({u_0} {u_1} {v_0} {v_1})^{1/2} 
\bigl(q^{14}+q^{16}+q^{10} {u_1}-
q^8 {u_0} {u_1}-2 q^{10} {u_0} {u_1}
-q^4 {u_0} {u_1}^2\\
&-2 q^6 {u_0} {u_1}^2-q^8 {u_0} {u_1}^2+
q^4 {u_0}^2 {u_1}^2+q^6 {u_0}^2 {u_1}^2+{u_0}^2 {u_1}^3
+q^2 {u_0}^2 {u_1}^3\\
&-q^9 ({u_0} {u_1} 
{v_0} {v_1})^{1/2}+q^6 {u_1} {v_1}+q^8 {u_1} {v_1}-
q^2 {u_0} {u_1}^2 {v_1}-q^4 {u_0} {u_1}^2 {v_1}\bigr)\\
&+ B\,{u_1} {v_0} \bigl(q^{16}+q^{18}+
q^{12} {u_1}-q^{12} {u_0} {u_1}-q^8 {u_0} {u_1}^2+q^{12} {v_1}
+ q^{10} {u_1} {v_1}\\
&+q^4 {u_0} {u_1} {v_1}-q^4 {u_0} {u_1}^2 {v_1}-
{u_0}^2 {u_1}^2 {v_1}\bigr) +  \bigl(q^{22} {v_0}-
q^{10} {u_0} {u_1}^2 {v_0}-q^{14} {u_0} {u_1}^2 {v_0}\\
&-q^6 {u_0} {u_1}^3 {v_0}+q^6 {u_0}^2 {u_1}^3 {v_0}+
q^2 {u_0}^2 {u_1}^4 {v_0}+q^{10} {u_0} {u_1} {v_1}
-q^4 {u_0}^2 {u_1}^2 {v_1}\\
&+{u_0}^3 {u_1}^3 {v_1}+q^{14} {u_1} {v_0} {v_1}+
2 q^4 {u_0}^2 {u_1}^2 {v_0} {v_1}-q^6 {u_0} {u_1}^3 {v_0} {v_1}
-2 {u_0}^3 {u_1}^3 {v_0} {v_1}\\
&-q^2 {u_0}^3 {u_1}^3 {v_0} {v_1}-q^4 {u_0}^2 {u_1}^2 {v_0}^2 {v_1}+
{u_0}^3 {u_1}^3 {v_0}^2 {v_1}+q^6 {u_1}^2 {v_0} {v_1}^2\bigr)\Bigr).
\end{align*}
}
\begin{align*}
&\frac{\breve{q}}{u_0^2 u_1^4 v_0 v_1^2}
\Bigl(\,A B\, \breve{q}^{6}u_1(\breve{q} u_0 u_1 v_0 v_1)^{1/2}+
B^2\, \breve{q}^{4}u_1^2 v_0 v_1-
A^2\, \breve{q}^{3}u_0^2 u_1^2 v_1\\
&+A\, (\breve{q} u_0 u_1 v_0 v_1)^{1/2} 
\bigl(\breve{q}^{7}+\breve{q}^{8}+\breve{q}^{5} u_1-
\breve{q}^4 u_0 u_1-2 \breve{q}^{5} u_0 u_1
-\breve{q}^2 u_0 u_1^2\\
&-2 \breve{q}^3 u_0 u_1^2-\breve{q}^4 u_0 u_1^2+
\breve{q}^2 u_0^2 u_1^2+\breve{q}^3 u_0^2 u_1^2
+u_0^2 u_1^3+\breve{q} u_0^2 u_1^3\\
&-\breve{q}^4 (\breve{q}u_0 u_1 
v_0 v_1)^{1/2}+\breve{q}^3 u_1 v_1+\breve{q}^4 u_1 v_1-
\breve{q} u_0 u_1^2 v_1-\breve{q}^2 u_0 u_1^2 v_1\bigr)\\
&+ B\,u_1 v_0 \bigl(\breve{q}^{8}+\breve{q}^{9}+
\breve{q}^{6} u_1-\breve{q}^{6} u_0 u_1-
\breve{q}^4 u_0 u_1^2+\breve{q}^{6} v_1
+ \breve{q}^{5} u_1 v_1\\
&+\breve{q}^2 u_0 u_1 v_1-\breve{q}^2 u_0 u_1^2 v_1-
u_0^2 u_1^2 v_1\bigr) +  \bigl(\breve{q}^{11} v_0-
\breve{q}^{5} u_0 u_1^2 v_0-\breve{q}^{7} u_0 u_1^2 v_0\\
&-\breve{q}^3 u_0 u_1^3 v_0+\breve{q}^3 u_0^2 u_1^3 v_0+
\breve{q} u_0^2 u_1^4 v_0+\breve{q}^{5} u_0 u_1 v_1
-\breve{q}^2 u_0^2 u_1^2 v_1\\
&+u_0^3 u_1^3 v_1+\breve{q}^{7} u_1 v_0 v_1+
2 \breve{q}^2 u_0^2 u_1^2 v_0 v_1-
\breve{q}^3 u_0 u_1^3 v_0 v_1
-2 u_0^3 u_1^3 v_0 v_1\\
&-\breve{q} u_0^3 u_1^3 v_0 v_1-
\breve{q}^2 u_0^2 u_1^2 v_0^2 v_1
+u_0^3 u_1^3 v_0^2 v_1
+\breve{q}^3 u_1^2 v_0 v_1^2\bigr)\Bigr).
\end{align*}

}
\renewcommand{\baselinestretch}{1.2} 

Let us consider the reduction  $u_1=t$ and $u_0=1=v_0=v_1$. 
Then $A=0=B$ and the substitution $\breve{q}\mapsto q^{-1/2}$ 
readily results in 
\begin{align*}
&J\!D^{A_1}_{3,2}(2\,;\,q,t)\!=\!\frac{q^{-6}}{{t}^4}
\bigl(1\!-\!q^{3} {t}^2\!-\!q^{2} {t}^2 
\!+\!q^{5} {t}^4\!+\!q^{3} {t}\!+\!q^{2} {t} 
\!-\!q^{4} {t}^3\!-\!q^{5} {t}^3\!+\!q^{4} {t}^2 \bigr).
\end{align*}

Making now $u_1=t=u_0$, \ $v_1=1=v_0$ and $\,\breve{q}\mapsto
q^{-1}\,$, we obtain
\begin{align*}
&J\!D^{A_1}_{4,3}(2\,;\,q,t)\!=\!\frac{q^{-12}}{{t}^6}\bigl(
 1 + q^2 t + q^3 t + q^4 t + q^5 t - q^2 t^2 - q^3 t^2 + 2 q^6 t^2 \\
&+ q^7 t^2 + q^8 t^2 - q^4 t^3 - 2 q^5 t^3 - 2 q^6 t^3 - 2 q^7 t^3 
+ q^9 t^3 + q^5 t^4 - 2 q^8 t^4\\
&- 2 q^9 t^4 + q^7 t^5 + q^8 t^5 + q^9 t^5 - q^{11} t^5 + 
q^{11} t^6\bigr).
\end{align*}
\medskip

\subsection{\bf Superpolynomials for  
\texorpdfstring{{\mathversion{bold}$T(2\rr+1,2)$}}{T(r,s)}}
Applying the previous theory to torus knots $T(2\pp+1,2)$,
we obtain that
\begin{align}\label{2p+1-rel}
&\j\!D_{2\pp\!+\!1,1}(m;\,q^2,u_1\!=\!t,u_0\!=\!t,1,1)=
J\!D^{A_1}_{2\pp\!+\!1,2}(m;\,q,t) \hbox{\,\, and}\\ 
\j\!D_{2\pp\!+\!1,1}(&m;\,q^2,u_1\!=\!t,u_0\!=\!1,1,1)\!=\!
J\!D^{A_1}_{2\pp\!+\!1,1}(m;\,q^2,t)\!=\!
(q^{\frac{m}{2}}t)^{-m\frac{2\pp+1}{2}}.
%\thicksim 1.
\notag
\end{align} 
On the other hand, the $a$\~degree of
$\h_{2\pp+1,2}(\om_1;\,q,t,a)$ is $1$ and it is completely
determined by the following two relations from
(\ref{jones-sup}):
\begin{align}\label{sup-rel}
\h_{\rr,\ss}(\om_1;\,q,t,a\!=\!-t)\!=\! 1,\,\,
\h_{\rr,\ss}(\om_1;\,q,t,a\!=\!-t^2)\!=\!
\tilde{J\!D}_{\rr,\ss}^{A_{1}}(1;\,q,t),
\end{align}
where $\rr=2\pp+1,\ss=2$.

Therefore the DAHA-superpolynomial 
$\h_{2\pp+1,2}(\om_1;\,q,t,a)$ from Theorem \ref{STABILIZ}
must be tilde-proportional (up to fractional powers of
$q,t$)  to 
$\imath a^{1/2}\j\!D_{2\pp+1,1}
(m;\,q^2,u_1\!=\!t,u_0\!=\!-at^{-1},1,1)$.

The imaginary unit $\imath$ occurs here due to
square roots of the parameters $u_1,u_0,v_1,v_0$ 
in the formulas for $\j\!D$\~polynomials; the coefficient
of proportionality $\imath a^{1/2}$ corresponds to the
choice $\imath a^{1/2}=(u_1u_0)^{1/2}$.

Indeed, $\a=0$ in (\ref{varsQ}) for $(2\pp+1,1)$
and any $\,m$, which means that
$(q^{\frac{1}{2}}u_0u_1)^{\frac{1}{2}}
\j\!D_{2\pp+1,1}(m;\,q,u,v)$ 
does not involve square roots of $q^{1/2}$, $u_i, v_i$. 
This is because $\j\!D_{2\pp+1,1}(m;\,q,u,v)$
does not depend on $v_1$. 
\smallskip

%Involving the recurrence relations
%in terms of $m$ for DAHA- superpolynomials and 
%DAHA-Jones polynomials of type $C^\vee C_1$, which will be
%not discussed in the present paper, we extend this proportionality
%to all $m\in \N$. Here the nonsymmetric polynomials are
%used; see (\ref{jones-nonsym})
%and (\ref{jones-nonsymc}). There is no less computational
%proof of the proposition below at the moment.

\begin{theorem}\label{SUPER-CC}
Let us substitute
\begin{align}\label{a-subst}
&q^{1/2}\!\mapsto\! q,\ \,u_1^{1/2}\!\mapsto\! t^{1/2},\ \, 
u_0^{1/2}\!\mapsto\! \imath a^{1/2} t^{-1/2},\ \,
v_0^{1/2}=1=v_1^{1/2}
\end{align}
in the formula for $\j\!D_{2\pp+1,1}(m;\,q,u,v)$.
Then the resulting polynomial and the DAHA-superpolynomial 
$\h_{2\pp+1,2}(m\om_1;\,q,t,a)$ from Theorem \ref{STABILIZ}
are proportional to each other. Namely,
\begin{align}\label{prop-a-h}
&q^{m^2\,\rr/2}t^{m(\rr-1)/2}\,
 \j\!D_{2\pp+1,1}(m;\,q^2,u_1\!=\!t,u_0\!=\!-\!at^{-1},v_i\!=\!1)\\
=&(\imath^m/a^{\frac{m}{2}})\, \h_{2\pp+1,2}(m\om_1;\,q,t,a) 
\hbox{\,\, for any\,\, } 
\pp\in \Z_+,\, m\in \N.\notag
\end{align}
\end{theorem}
%\vskip -0.2cm
%\sq

{\it Proof.} The key is using the
color exchange. We have the following.

\begin{lemma}\label{LEMCC}
For an arbitrary admissible pair $(\rr,\ss)$
and $N=m,m+1,\ldots, 2m-1$, correspondingly,
$n\equal N-m=0,1,\ldots, m-1$ the 
following $\,m\,$ color-exchange relations
from (\ref{bcw-relc}) hold: 
\begin{align}\label{2p+1-rs}
&\j\!D_{\rr,\ss}(m;\,q^2,u_1=t,u_0=q^{-N}t^{-1},v_1,v_0)\\
=&\j\!D_{\rr,\ss}\,(n;\,q^2,u_1=t,u_0=q^{-N}t^{-1},v_1,v_0).\notag
\end{align}
On the other hand, using (\ref{jdiota}) and then
applying the $q\leftrightarrow t^{-1}$\~duality one has:
\begin{align}\label{super-rs}
&\h_{\rr,\ss}(\om_m;\,q,t,a=-t^{N})\,\, = \,\,
\h_{\rr,\ss}(\om_n;\,q,t,a=-t^{N})\ \Rightarrow\notag\\
&\h_{\rr,\ss}(m\om_1;\,q,t,a=-q^{-N})=
\h_{\rr,\ss}(n\om_1;\,q,t,a=-q^{-N}).
\end{align}
\end{lemma}
\vskip -1.3cm \sq
\vskip 0.2cm

Thus we have obtain that (\ref{prop-a-h}) holds
for $a\!=\!-q^{-m},\ldots, a\!=\!-q^{1\!-\!2m}$. 
Since the $a$\~degrees of both expressions are
$m$, we need $m+1\,$ coinciding independent 
specializations. Relations  (\ref{2p+1-rel}) 
and (\ref{sup-rel})
provide them at $a\!=\!-t$ and 
$a\!=\!-t^2$. This concludes the proof of the theorem. \sq 
%$\j\!D_{2\pp+1,1}(m;\,q^2,t,-\!at^{-1},1,1)$ and
%$\h_{2\pp+1,2}(m\om_1;\,q,t,a)$.  

\smallskip

The theorem is actually very surprising, since the $\h$
require all $A_n$ for their definition, while the
$\j\!D$ are obtained by a rank-one calculation
(for $(2\pp+1,2))$. 
%This somehow resembles using
%{\em rational DAHA\,} for such knots, but we do not 
%see any connection.  
In process of proving the theorem
and the lemma, we obtain the following two formulas:
\begin{align}\label{2p+1form}
&\h_{2\pp+1,2}(m\om_1;\,q,t,a) =
\frac{(q;q)_m}{(-a;q)_m\,(1-t)}\sum_{k=0}^m (-1)^{m-k}
(qt)^{\frac{m-k}{2}}\\
\times\!\Bigl(
(q^{\frac{m(m\!+\!1)}{2}}&\!\!-\!q^{\frac{k(k\!+\!1)}{2}})
\bigl\{\!\frac{\hbox{\small\em t}}{\hbox{\small\em q}}\bigr\}
^{\!\!\frac{m\!-\!k}{2}}\Bigr)^{2\pp\!+1}
\,\frac{(t;q)_k\,(-a;q)_{m\!+\!k}\,
(-a/t;q)_{m\!-\!k}\,(1\!-\!q^{2k}t) }
{(q;q)_k(qt;q)_{m\!+\!k}\,(q;q)_{m\!-\!k}};\notag\\
\h_{3,2}(m\om_1)&\! =\!\!
\sum_{k=0}^m q^{mk}t^k \frac{(q;q)_m(-a/t;q)_k}
{(q;q)_k(q;q)_{m-k}},\ (x;q)_n\!=\!(1\!-\!x)\cdots 
(1\!-\!x q^{n-1}).\notag
\end{align}

Note that a simple calculation with the
first summation from (\ref{2p+1form}) gives that 
its value at $\,a=-q^{-m}\,$ is $\,q^{\pp\,m^2}t^{m\,\pp}$.
Therefore, 
\begin{align*}
&\h_{2\pp+1,2}(m\om_1;\,q,t,-q^{-m})= 
\,q^{\pp\,m^2}t^{m\,\pp},\ 
\j\!D_{2\pp+1,1}(m;\,q^2,t,q^{-m}t^{-1},1,1)\\
=&q^{-m^2(\pp+1/2)}t^{-m(\rr-1)/2}\,
(\imath^m/a^{\frac{m}{2}})\, \h_{2\pp+1,2}(m\om_1;\,q,t,a)
\hbox{\, for\, }
a^{1/2}\mapsto \imath q^{-m/2}.
\end{align*}
The latter is $1$, which matches formula
(\ref{2p+1-rs}) for $N=m,n=0$.

\medskip
The first formula in (\ref{2p+1form}) was conjectured in 
(229) of \cite{DMMSS} based on extensive calculations there 
with superpolynomials of torus and similar knots.
It was conjectured in \cite{FGS} too as well as
the second (trefoil) formula; see (2.23), (2.32) there.
The second one is a generalization of Habiro's formula 
(around 2000) for $\pp\!=\!1, a\!=\!-t^2, t\!=\!q$.

The approach from \cite{FGS} is based on  
(conjectural) symmetries of superpolynomials. The method 
from \cite{DMMSS} is somewhat different (and more
understandable algebraically); the authors assume 
the split of the formulas for superpolynomials into the 
components reflecting respectively the knot itself 
and the contribution of the corresponding Macdonald polynomials. 
Such a split can be readily seen in 
(\ref{2p+1form}). Both papers do not provide details;
we mention that (2.12) from \cite{CJ} links our 
parameters to those used in \cite{DMMSS}. 
\smallskip

Comparing the first formula for $\pp=1$ and the second one,
we come to a certain $q$\~identity, which is not too
difficult to check; see (\ref{Aznmp}) below.
To clarify its meaning, let us perform the reduction 
$q^{1/2}\mapsto 1$; here we will double check the output using 
the evaluation formula
(\ref{jones-seval}). The first formula becomes
\begin{align*}
&\frac{1}{(1\!-\!t)}\!\sum_{k=0}^m 
t^{k(\pp\!+\!1)}(\!-\!1)^k\binom{m}{k}(1\!+\!a)^{m-k}
\left(1\!+\!\frac{a}{t}\right)^k\!=
\!\Bigl(\frac{1\!-\!t^{\pp\!+\!1}}{1-t}\!+\!
\frac{1\!-\!t^{\pp}}{1\!-\!t}a
\Bigr)^m,
\end{align*}
which does coincide for $\pp=1$ with the second (as  
$q^{1/2}\mapsto 1$):
\begin{align*}
&\sum_{k=0}^m t^k\binom{m}{k}\left(1+\frac{a}{t}\right)^{k}
=\bigl(1+t+a\bigr)^m.
\end{align*}
\smallskip

Let us expand the first formula from
(\ref{2p+1form}) in terms of $(-a/t;q)_m$ with the 
coefficients from the second (trefoil) one. One has:
\begin{align}\label{Aznmp}
&\h_{2\pp+1,2}(m\om_1;\,q,t,a)\!\equal\sum_{n=0}^m q^{mn}t^n
A(n;\,m,\pp)\,
\frac{(q;q)_m\,(-\frac{a}{t};q)_m}{(q;q)_{m-n}\,(q;q)_n},\\
A(n;m&+1,\pp)\,=\, \Xi\bigl[\,A(n;\,m,\pp)\,\bigr],\ \,
\Xi\bigl[\,\sum C_{ij}q^i t^j\,\bigr]\equal
\sum C_{ij}q^{i+2j}t^j,
\notag\\
&\hbox{and\, } A(n;\, m,\pp\to\!\infty)\,=\,
\prod_{i=2m-n}^{2m-1}\frac{1}{1\!-\!q^i t}\hbox{\,\, for fixed\, } 
0\le n\le m.
\label{Azstab}
\end{align}
Relations for  $A(n;\,m,\pp)$ result in the second 
formula from (\ref{2p+1form}). They follow from a variant of 
the color exchange for $C^\vee C_1$:
\begin{align*}%\label{Azcolor}
\h_{2\pp\!+\!1,2}(m\om_1;\,q,t,a\!=\!-tq^{-n}) =
\h_{2\pp\!+\!1,2}(n\om_1;\,q,t,a\!=\!-tq^{-m})
\mid_{t\mapsto t q^{2(m\!-\!n)}}
\end{align*}
for any $n,m>0$; take $n=m-1$ for (\ref{Aznmp}).
We omit the justification.
% Thus it suffices to focus on the leading 
%coefficients $A(n;\,n,\pp)$ for $n=0,1,\ldots\,$.   

By the limit in (\ref{Azstab}), we mean that for any 
given $0\le n\le m$, 
each $A(n;\, m,\pp)$ is a subsum (with positive coefficients) 
of $A(n;\, m,\infty)$ and eventually will include any term of
the latter for sufficiently large $\pp$. It is quite likely
that the coefficients of 
$A(n;\, m,\pp)$ monotonically increase
in terms of $\,\pp\,$ for any fixed $n,m$ and $q^at^b$. 
This may be somehow connected with Conjecture 7.4 from
\cite{GGS}.
%\medskip

\subsection{\bf Generalized Verlinde algebras}
We will extend {\em perfect $\HH^{A_1}_{q,t}$\~modules\,}
at roots of unity, also called {\em generalized nonsymmetric 
Verlinde algebras\,}, from $A_1$ to $C^\vee C_1$  
and discuss the action of the absolute Galois group there.
Following \cite{C101}, such a module is 
a quotient of the polynomial representations by
the radical of the evaluation pairing provided 
its irreducibility and the projective action of
$PSL^\wedge_{\,2}(\Z)$ there.   We will only extend 
the {\em main series\,} of such modules (which
directly generalize the classical Verlinde algebras).
%We set $u_i^{1/2}=q^{k_i/2}$ and $v_i^{1/2}=q^{l_i/2}$.

The evaluation pairing for $C^\vee C_1$ is as follows:
$$
\{f,g\}_{ev}=\{\vph(f)g\}_{ev}=
\vph\bigl(f(\X)\bigr)(g)(u_1v_1^{-1/2}).
$$  
It is defined on $f\times g\in \v^\vph\times \v$, where 
$\v^\vph$ is the image of $\v$ under $\vph$ with the 
action of algebra  $\sH^\vph_{q,u,v}$ obtained from $\sH_{q,u,v}$
by the transformation $u_0\leftrightarrow v_1$; 
see (\ref{vphcc}). 
We will set $u_i^{1/2}=q^{k_i/2}$ and $v_i^{1/2}=q^{l_i/2}$ 
for $i=0,1$ in the following theorem.

\begin{theorem}\label{VERCC}
Let $q^{1/4}$ be a primitive 
$(4\!N)${\small th} root of unity and $\,u_1^{1/2}=q^{k_1/2}$
for $2k_1\in \Z$ subject to the inequality
$0< k_1< N/2$. We set $\,M\!=\!N\!-\!2k_1$.
If \,$2k_0\in \Z$, we will assume that $\,0\le k_0\le k_1$,
which is a $\tau_-$\~invariant condition.
We also impose the condition
\begin{align}\label{vergncc}
&1-\ep\, q^{\frac{j}{2}+\frac{1}{4}}
(u_1u_0v_1v_0^{\,\ep})^{\frac{1}{2}}
\neq 0 \for \ep=\pm1 \and 0\le j \le M.
\end{align} 
If $\,u_0,v_1,v_0\in \Z/2$, then
$\frac{1}{2}\!+\!k_1\!+\!k_0\!+\!l_1\!+\!l_0\not\in \Z$ 
is sufficient, which is a $GL_{\,2}^\wedge(\Z)$\~invariant condition.

(i) Then the quotient $\V$ of $\v$ by the (right) radical $Rad$ of
the evaluation pairing
$\{f,g\}_{ev}$ is irreducible of dimension $2M=2\!N\!-\!4k_1$ with 
semisimple action of $\Y$, which spectrum is simple. 
More explicitly,\, the polynomial $\e_{-M}$ is well defined,
$\{\e_{-M}\}_{ev}=0$ and $\V=\v/(\e_{-M})$\, as an 
$\X$\~module. Furthermore, the polynomials $\,\e_{m}\,$ are well 
defined for $\,-M<m\le M\,$ and form a basis of $\V$; 
their evaluations are all nonzero for such $m$.

(ii) Assuming that the parameters $\,u_0,v_1,v_0\,$ are formal 
variables, the action of the group $\Ga_{\,0}^\wedge$ can be
represented by skew-linear maps from $\V$ to its images 
under the corresponding transformation of the 
parameters. For instance, the subgroup $\bar{\Ga}^\wedge(2)$ 
of $PSL_{\,2}^\wedge(\Z)$ generated by $\{\tau_{\pm}^2\}$ 
preserves $\V$ for any admissible values of the parameters
and its action is linear in this module. Here 
$PSL_{\,2}^\wedge(\Z)$ is the
span of $\{\tau_{\pm}\}$ in $GL_{\,2}^{\wedge}(\Z)$,\, 
$\bar{\Ga}^\wedge(2)$ corresponds to the group 
$\bar{\Ga}(2)\subset PSL_{\,2}(\Z)$ of matrices identical 
modulo $(2)$.

(iii) Continuing to assume $\,M\!=\!N\!-\!2k_1\in \N$,
let 
\begin{align}\label{tilde-perfect}
&\,k_0,l_1,l_0\in \Z_+/2,\,\ k_0\!\le\! k_1\! \hbox{\,\, and\,\, }
1/2+k_1+k_0+l_1+l_0\in \Z,
\end{align}
which violates (\ref{vergncc}). Then the 
claims from (i,ii) still hold if $\,k_1\!-\!k_0\!>\!l_1\!+\!l_0$. 
However, the opposite inequality $\,k_1\!-\!k_0\!<\!l_1\!+\!l_0\,$ 
(the strict equality is impossible here due to (\ref{tilde-perfect})) 
and the additional condition 
$\,M'\!\equal\!N\!+\!1/2\! -\!(k_1\!+\!k_0\!+\!l_1\!+\!l_0
)>0\,$ give that $\,\v/Rad\equal\V'=\v/(\e_{M'})\,$ 
is irreducible of dimension $\,2\!M'\!\!-\!1\!<\!2M\,$ with 
simple spectrum of $\Y$.

(iv) We return to imposing (\ref{vergncc}). 
For instance, one may assume that
\begin{align}\label{rigid-cc}
k_i,l_i\!\in\!\Z_+/2,\  \,k_0\!\le\! k_1\!<\! N/2
\hbox{\,\, and\,\, }
\frac{1}{2}\!+\!k_1\!+\!k_0\!+\!l_1\!+\!l_0\not\in \Z.
\end{align}
Then $\V$ is rigid; more exactly, it is a unique 
$\sH_{q,u,v}$\~module of dimension $\,2\!M$ up to isomorphisms. 
For instance, this results in the existence of the projective 
action of $\bar{\Ga}^\wedge(2)$ in $\V$, claimed in (ii).
\end{theorem}
\smallskip

{\it Proof.} Part $(i)$ is straightforward using
the evaluation formulas (\ref{cceval-1}), 
(\ref{cceval-2}) and the spectrum of $\Y\,$ from
formula (\ref{Eccheck}). 
One can use formula (\ref{tauminpol}) for the action of 
$\tau_-$ on the $\e$\~polynomials and the relation 
$\vph \tau_-\vph=\tau_+$ to verify $(ii)$.

Another proof of $(ii)$ is as follows. We readily establish
$(iv)$ for $\,u_0,v_1,v_0\,$ in a neighborhood of
the specialization  $\,u_0^{1/2}=1=v_1^{1/2}=v_0^{1/2}$,
since it holds for $A_1$. See case $(\al)$
of Theorem 1.4 from \cite{CG}. This gives that for {\em all\,}
parameters, there exists a skew-symmetric
isomorphism $\V\to\hat{\ga}(\V)$, where the parameters
are conjugated by $\hat{\ga}$. Here 
$\ga\in\bar{\Ga}^\wedge(2)$, assuming that 
the parameters $u_0,v_1,v_0$ are 
{\em admissible\,}, i.e. such that the existence, 
irreducibility and semisimplicity of $\V$ hold.
\smallskip
  
Claim $(iii)$ is straightforward. The last claim uses
the technique of intertwiners; formulas (\ref{Sinter})
and counterparts of (\ref{Sintere}) in arbitrary modules
are needed. The proof is parallel to
the classification of the irreducible modules for $A_1$
from \cite{C101}. However $(iv)$ and its 
counterpart from Theorem 1.4 of \cite{CG} do not require
knowing the whole classification of irreducible modules;
this is a simpler problem.
\sq
\smallskip

We mention here paper \cite{ObS} devoted to
the rigid modules for $C^\vee C_1$ apart from the
roots of unity, the corresponding Deligne-Simpson
problem and Crawley-Boevey's results.  Their rigid modules 
remain rigid for sufficiently general $N$ and can be used too. 
However a systematic theory of perfect and rigid modules at 
roots of unity for $C^\vee C_1$ is needed.
\medskip

As it was already mentioned in the proof of $(ii)$,
the module $\V\,$ becomes $V_{2\!N\!-\!4k}\,$ from Part $(i)$
of Theorem 1.2 from \cite{CG} for $k_1=k, k_0=0=
l_1=l_0$. The second reduction to $A_1$ from
Theorem \ref{CCTOONE} above, is
$k_1=k=k_0$ and $q\mapsto q^2$. This corresponds to
the so-called {\em Little Verlinde algebra} from
\cite{C101}, defined for $Y^2$ instead of $Y$; see
also Section 1.4 from \cite{CG}.

Claim $(ii)$ of the theorem
is directly related to DAHA-Jones polynomials.
For the sake of transparency,
we will restrict ourselves to the subgroup 
$\bar{\Ga}^\wedge(2)\subset \Ga_{\,0}^\wedge(2)$ in
the definition from (\ref{jones-dc}) of 
Theorem \ref{CCTHM}.

We use that the evaluation
$\{f\}_{ev}$ factors through $\V$\,; the notation will be
$\{f\}_{ev}^\V$.
Let $\tga\in \hbox{Aut}(\V)$ be the matrix 
representing $\hat{\ga}\in \bar{\Ga}^\wedge(2)$
for $\ga=\ga_{\rr,\ss}$.

\begin{corollary}\label{JDVER}
(i) In the setting from $(i,ii)$ of 
Theorem \ref{VERCC},
\begin{align}\label{jones-ver}
&\j\!D_{\rr,\ss}(m'\,;\,q,u,v)\ =\  
\{\tga\,\e_m\,\tga^{-1}\}_{ev}^{\V}\,/\,
\{\e_m\}_{ev}^{\V},\\
&\hbox{where\, }
0\le m'\!=\!m\!\!\mod(2\!N)  
\hbox{\, for\, } 0\le m\le M. \notag
\end{align}
\end{corollary}
\vskip -1.2cm \sq
\vskip 0.2cm

This corollary establishes a link to the 
approach to Jones polynomials of torus knots
and their superpolynomials based on generalized 
Verlinde algebras. See \cite{AS}, which
triggered \cite{CJ}; the projective
action of $PSL_2(\Z)$ used in \cite{AS} is
from \cite{Ki}.
One can see which kind of {\em refined Chern-Simons
theory\, } can be needed to match the
DAHA-Jones polynomials. It must involve 
$3$ essentially free parameters ($u_0,v_1,v_0$)
and a half-integral parameter $1\le\! k_1\le\! N/2\,$ from
the relation $u_1=q^{k_1}$. Note that the 
$C^\vee C_1$\~theory provides Selberg-type integrals 
with $4$ parameters (and $q$). Similar 
integrals did occur in conformal field theory,
but their relevance to our work is not clear.
\smallskip

{\sf Absolute Galois group.}
We will conclude this paper with an
arithmetic extension of Corollary \ref{JDVER}.
The key point is that practically any automorphisms
of the spherical DAHA can be used there. Following \cite{CG},
we show that the {\em absolute Galois group\,} can 
substitute for  $GL^\wedge_{\,2}(\Z)$.

We continue to assume that $q^{1/4}$ is a primitive
$(4\!N)${\small th} root of unity and assume that 
$u_i^{1/2},v^{1/2}\in \Z[q^{1/4}]$ for $i=0,1$.
Then we fix a prime number $p\,$ and a number (conductor) 
$n\in \N$, extending $\,(p)\,$ to a prime ideal $\mathfrak{p}$
in $\Z[q^{1/4}]$. We need the counterparts $\V_{\mathfrak{p},n}$
of module $\V$ from Theorem \ref{VERCC} defined over the 
rings $\Z_{\mathfrak{p},n}\equal \Z[q^{1/4}]/(\mathfrak{p}^n)$
and the $\mathfrak{p}$\~adic limit $\V_{\mathfrak{p}}$
over $\Z_{\mathfrak{p}}=\varprojlim_n\Z_{\mathfrak{p},n}$.
Here we assume that $\V_{\mathfrak{p},n}$ and $\V_{\mathfrak{p}}$
exist, are irreducible,
rigid and $\Y$\~semisimple for such $\mathfrak{p}$. Due to 
Hensel's Lemma, this must be established for $n=1$.  
For generic $u_0,v_1,v_0$, the condition \,gcd$(p,2\!N)=1$
is sufficient.
\smallskip

Let us consider the projective line $P^1_{\C}$ 
punctured at pairwise distinct points $O_1,O_2,O_3,O_4$.
The fundamental group $\Pi_1=\pi_1(P^1_{\C}\setminus \{O_j\};o)$
(for a base point $o\not\in  \{O_j\}$)
is generated by the standard counterclockwise loops $A_j$
around these points satisfying the relation $A_1A_2A_3A_4\!=\!1$.
Let $\a_{\mathfrak{p},n}$ be the image of $\Pi_1$ in 
Aut${}_{\Z_{\mathfrak{p},n}}\V_{\mathfrak{p},n}$, 
where $A_1,A_2,A_3,A_4$ map to  
$q^{1/4}U_1,U_0,V_0,V_1$  acting in $\V_{\mathfrak{p},n}$
(in this order, matching that in (\ref{CCdoub})).

We use the Riemann Existence Theorem to construct  
the cover $F\to P^1_{\C}$ ramified only at $\{O_j\}$
such that Aut$(F\to P^1)=\a_{\mathfrak{p},n}$ and
the images of $A_i$ generate the (cyclic) subgroups of 
the elements fixing certain ramification points over $O_i$.
See \cite{CG} here and below. Assuming that $O_1,\ldots,O_4$ 
are defined over $\Q(q^{1/4})$ and that $u_i,v_i\in \Z[q^{1/4}]$,
we claim that $F$ can be actually defined over $\Q(q^{1/4})$
and the {\em absolute Galois group\,}
$\g=$Gal$\bigl(\overline{\Q}/\Q(q^{1/4})\bigr)$
acts in Aut$(F/P^1)$ by outer automorphisms.
We use the rigidity of $\V_{\mathfrak{p},n}$, which
also gives that $\g$ projectively acts there by linear 
automorphisms. For $g\in \g$, let $\phi_g$ be the corresponding
matrix (unique up to proportionality). This is 
a modification of Belyi's theory enriched
by Katz' {\em linear rigidity\, } (M. Dettweiler and others).

\begin{conjecture}\label{CONJGALOIS}
Following (\ref{jones-ver}), we set 
\begin{align}\label{jones-gal}
\j\!D_{g}(m)\!= \varprojlim_n
\{\phi_g\,\,\e_m\,\phi_g^{-1}\}_{ev}^{\V}\,/\,
\{\e_m\}_{ev}^{\V} \for \V=\V_{\mathfrak{p},n},\, g\in \g,
\end{align}
where $1\le m\le M=N\!-\!2k_1$. 
We conjecture that there
exists $\varrho_{m}(g)\in \Z[q^{1/4}]$ such that 
$\varrho_m(g)=\j\!D_{g}(m)$ for all places $p,\mathfrak{p}$
of good reduction\,,
where $\V_{\mathfrak{p}}$ is irreducible,
rigid and $\Y$\~semisimple.\sq
\end{conjecture}

The  conjecture is not directly related to the standard 
$\ell$\~adic conjectures for Tate modules in Number Theory 
(see below), though we follow the same lines. It can be 
easily verified if the image $\a$ of $\Pi_1$ in the 
whole Aut$\V$ is finite; then restricting the ring of 
coefficients to $\Z_{\mathfrak{p},n}$ is not needed. 
We described all such finite $\a$ in the case of  
$A_1$ in \cite{CG}. This seems significantly more ramified
for $C^\vee C_1$, though hopefully doable under
the conditions from (\ref{rigid-cc}).
\smallskip
  
The absence of denominators in  $\varrho_m(g)$
at the places  $p,\mathfrak{p}\,$ of 
{\em bad reduction\,}
is a direct counterpart of the polynomiality of 
$\j\!D_{\rr,\ss}$. Here one can actually drop
the inequality $m\!<\!N\!-\!4k_1$. Also, we can take 
$\g=$Gal$(\overline{\Q}/\Q)$ provided that $O_1$ and
the triple $\{O_2,O_3,O_4\}$ are defined over $\Q$.
Then $\phi_g$ will be coupled with the restriction $g'$ 
of $g\in\g$ to  $\Q[q^{1/4}]$ and 
$\j\!D_{g}(m)$ will be conjugated by $g'$. 
We modify (\ref{jones-gal}) 
following (\ref{jones-dc}) if $\,O_2,O_3,O_4\,$ 
are permuted by $\,g'$; cf. \cite{CG}.
\smallskip

The classical (unramified) {\em Tate modules\,}
require considering the stabilization in terms of $N\!=\!\ell^e$
as $e\to \!\infty$ for a given prime $\ell$, which
corresponds to $k_1\!=\!0\!=\!k_0\!=\!l_1\!=\!l_0$ in our ramified
theory ($\ell$ must be odd). Then the group $\a\,$ is finite 
and we do not need using $p$. Furthermore, $\g\,$ acts in $\V\,$ via 
$\,\psi_g\in GL_2(\Z)$ modulo $(\ell^e)$ and 
$\,\psi_g\,$ becomes a matrix in 
$GL_2(\Z_{\ell})$ upon the $\ell$\~adic limit $e\to \infty$.

Classically, the eigenvalues of $\psi_g$ are considered,
which depend only on the conjugacy class of $\,g\in\g$.
Our $\j\!D_{g}(m)$ is not conjugacy invariant. It equals
$q^{-m^2\rr\,\ss\,\kappa/4}$ in this case for
$\rr,\ss,\kappa\in \Z_{\ell}$, where $(\rr,\ss)^{tr}$ is the 
first column of the matrix $\psi_g$ and 
$\,g(q)=q^{\kappa}$. We represent $q=1+\ep$ and tend 
$e\to\infty$, expanding $\j\!D_{g}(m)$ in terms of $\ep$.
\smallskip

We do not generally expect the groups $\a_{\mathfrak{p},n}$ to
act via $GL_{\,2}^\wedge(\Z)$ in $\V_{\mathfrak{p},n}$.
Actually one of the objectives of
this ``motivic" direction is to enlarge the group 
$GL_{\,2}^\wedge(\Z)$ of the standard DAHA-automorphisms
using the absolute Galois group.
We suggest some tools for managing the stabilization
in terms of $N$ in \cite{CG} based on 
{\em $q$\~deformed nonsymmetric Verlinde algebras\,},
but this is generally an open problem. 
%\smallskip

We mention that in contrast to the rest of the paper,
the arithmetic direction is restricted to the
rank one so far, since we use the Riemann Existence 
Theorem. We expect that $\,\g\,$ always acts in
$\mathfrak{p}$\~adic rigid DAHA modules at roots of 
unity for proper $(p)\subset \mathfrak{p}$ and that the 
monodromy of the 
{\em Knizhnik-Zamolodchikov-Bernard\,} equation associated 
with DAHA can be used for the justification.
However {\em KZB\,} is not involved here 
so far even in the rank-one case.
\smallskip

% \bigl\{\overset{t}{\underset{q}{-}}\bigr\}
% \raisebox{-1.5ex}{\hbox{\ $\overset{w}\longmapsto$\ }}

{\bf Acknowledgements.}
The paper was mainly written at RIMS;
the author thanks Hiraku Nakajima and RIMS for
the invitation and hospitality and
the Simons Foundation (which made this visit possible).

%\vskip -5cm
%\medskip
\bibliographystyle{unsrt}

\begin{thebibliography} {ABCD}
%\vfil
%\medskip

\bibitem [AS] {AS}
{M.~Aganagic}, and {S.~Shakirov},
{\em Knot homology from
refined Chern-Simons theory},
Preprint arXiv:1105.5117v1 [hep-th], 2011.

\bibitem [AW] {AW}
{R.~Askey}, and {J.~Wilson},
{\em  Some basic hypergeometric orthogonal polynomials
that generalize Jacobi polynomials},
Memoirs AMS {  319} (1985), 1--55.

\bibitem [BS] {BeS}
{Yu.~Berest}, and {P.~Samuelson},
{\em Double affine Hecke algebras and generalized Jones polynomials},
Preprint arXiv:1402.6032v2 [math.QA] (2014).

\bibitem [Bo] {Bo}
{N.~Bourbaki},
{\em Groupes et alg\`ebres de Lie}, Ch. { 4--6},
Hermann, Paris (1969).

\bibitem [Ch1] {C101}
{I.~Cherednik},
{\em Double affine Hecke algebras},
London Mathematical Society Lecture
Note Series, { 319}, Cambridge University Press, Cambridge, 2006.
%%\smallskip

\bibitem [Ch2] {C103}
\bysame
{\em Nonsemisimple Macdonald polynomials},
Selecta Mathematica { 14}: 3-4 (2009), 427--569.

\bibitem [Ch3] {C100}
\bysame,
{\em Irreducibility of perfect representations of double affine 
Hecke algebras},
Studies in Lie theory, Progr. Math. { 243}, 79--95,
Birkh\"auser Boston, Boston, MA, 2006.

\bibitem [Ch4] {CJ}
\bysame
{\em Jones polynomials of torus knots via DAHA},
Int. Math. Res. Notices, { 2013}:23 (2013), 5366--5425.

\bibitem [Ch5] {CG}
\bysame
{\em On Galois action in rigid DAHA modules},
Preprint arXiv:1310.2581v3 [math.QA] (2013).

\bibitem [DMS] {DMMSS}
{P.~Dunin-Barkowski},
and {A.~Mironov}, and {A.~Morozov}, and
{A.~Sleptsov}, and  {A.~Smirnov},
{\em Superpolynomials for toric knots from evolution 
induced by cut-and-join operators},
Preprint arXiv:1106.4305v2 [hep-th] (2012).

\bibitem [DGR] {DGR}
{N.~Dunfield}, and {S.~Gukov}, and {J.~Rasmussen},
{\em The superpolynomial for knot homologies},
Experimental Mathematics, { 15}:2 (2006), 129--159.
%arXiv:math/0505662v2 [math.GT].

\bibitem [FGS] {FGS}
{H.~Fuji}, and {S.~Gukov}, and {P.~Sulkowski},
{\em Super-$A$-polynomial for knots and BPS states},
Preprint arXiv:1205.1515v2 [hep-th] (2012).

\bibitem [GN] {GoN}
{E.~Gorsky}, and {A.~Negut},
{\em Refined knot invariants and Hilbert schemes},
Preprint arXiv:1304.3328v2 [math.RT] (2013).

\bibitem [GORS] {GORS}
{E.~Gorsky}, and {A.~Oblomkov},  
and {J.~Rasmussen}, and {V.~Shende},
{\em Torus knots and rational DAHA},
arXiv:1207.4523 (2012).

\bibitem [GGS] {GGS}
{E.~Gorsky}, and {S.~Gukov}, and {M.~Stosic},
{\em Quadruply-graded colored homology of knots},
Preprint arXiv:1304.3481 [math.QA] (2013).

\bibitem [GS] {GS}
{S.~Gukov}, and {M.~Stosic},
{\em Homological algebra of knots and BPS states},
Preprint arXiv:1112.0030v1 [hep-th] (2011).

\bibitem [Hi] {Hi}
{K.~Hikami},
{\em $q$-series and $L$-functions related to half-derivatives 
of the Andrews--Gordon identity}
Ramanujan J. { 11} (2006), 175--197.

%\bibitem [Jon] {Jon}
%{V.~F.~R.~Jones},
%{\em A polynomial invariant for knots via von 
%Neumann algebras}, 
%Bull. Amer. Math. Soc. { 12} (1985), 103--111.

\bibitem [KhR1] {KhR1}
{M.~Khovanov}, and  {L.~Rozansky}, 
{\em Matrix factorizations and link homology}, 
Fundamenta Mathematicae, { 199} (2008), 1--91.
%Preprint arxiv: math.QA/0401268.

\bibitem [KhR2] {KhR2}
\bysame, and \bysame,
{\em Matrix factorizations and link homology II}, 
Geometry and Topology, { 12} (2008), 1387--1425.
%Preprint arxiv: math.QA/0505056.

\vfill\eject %%%%%%%%%%%%%%%%%%%%%%%%%

\bibitem [Ki] {Ki}
{A.~Kirillov, Jr.},
{\em On inner product in modular tensor categories. I},
Jour. of AMS { 9} (1996), 1135--1170.


\bibitem [Ko] {Ko}
{T.~Koornwinder}, {\em Askey-Wilson polynomials for 
root systems of type $BC$}, Contemp. Math. { 138} (1992), 
189-–204.

\bibitem [Mac] {Mac}
{I.~Macdonald},
{\em Affine Hecke algebras and orthogonal polynomials}, 
Cambridge University Press (2003).

\bibitem [NS] {NS}
{M.~Noumi}, and {J.V.~Stokman}, 
{\em Askey-Wilson polynomials: 
an affine Hecke algebra approach}, 
in: Laredo Lectures on Orthogonal Polynomials and 
Special Functions, pp. 111--144, 
Nova Sci. Publ., Hauppauge, NY, 2004; 
Preprint arXiv:math/0001033v1. 

\bibitem [Ob] {Ob}
{A.~Oblomkov},
{\em Double affine Hecke algebras of rank 
$1$ and affine cubic surfaces}, 
Int. Math. Res. Not. { 2004}:18 (2004), 877--912. 

\bibitem [OS] {ObS}
{A.~Oblomkov}, and {E.~Stoica},
{\em
Finite dimensional representations of the double
affine Hecke algebra of rank 1},
Journal of Pure and Applied Algebra, { 213}:5 (2009), 766-–771.

\bibitem [Ras] {Ras}
{J.~Rasmussen},\ \ \ 
{\em Some differentials on Khovanov-Rozansky homology},
Preprint arXiv:math.GT/0607544 (2006).

\bibitem [RJ] {RJ}
{M.~Rosso}, and {V.~F.~R.~Jones},
{\em  On the invariants of torus knots derived 
from quantum groups},
Journal of Knot Theory and its Ramifications, { 2} (1993), 
97--112.

\bibitem [Rou] {Rou}
{R.~Rouquier},
{\em Khovanov-Rozansky homology and 2-braid groups},
Preprint arXiv:1203.5065 [math.RT] (2012).

\bibitem [Sa] {Sa}
{S.~Sahi}, 
{\em Nonsymmetric Koornwinder polynomials and duality}, 
Ann. of Math. (2) { 150}:1
(1999), 267-–282.

\bibitem [SV] {SV}
{O.~Schiffmann}, and {E.~Vasserot},
{\em The elliptic Hall algebra, Cherednik Hecke
algebras and Macdonald polynomials}, 
Compos. Math. { 147} (2011), 188--234.

\bibitem [Ste] {Ste}
{S.~Stevan},
{\em Chern-Simons invariants of torus links},
Annales Henri Poincar\'e { 11} (2010), 1201--1224.

\bibitem [Sto] {Sto}
{J.V.~Stokman}, 
{\em Difference Fourier transforms for 
non\-reduced root sys\-tems},
Sel. Math., New ser. { 9}:3 (2003), 409--494.

\bibitem [Web] {Web}
{B.~Webster},
{\em Knot invariants and higher representation theory},
Preprint arXiv:1309.3796 [math.GT] (2013).


%\vfill
\end{thebibliography}

\end{document}